\documentclass[ijoc,nonblindrev]{informs3}
\DoubleSpacedXI

\usepackage{xcolor}

\usepackage{natbib}
\usepackage{enumitem}
 \bibpunct[, ]{(}{)}{,}{a}{}{,}%
 \def\bibfont{\small}%
\TheoremsNumberedThrough     % Preferred (Theorem 1, Lemma 1, Theorem 2)
\ECRepeatTheorems

\EquationsNumberedThrough    % Default: (1), (2), ...
\usepackage[toc,page]{appendix}
\usepackage{subcaption}
\usepackage{mathtools}
\usepackage{bigstrut}
\usepackage{booktabs}
\usepackage{mathrsfs}

\usepackage[utf8x]{inputenc}
\usepackage{multicol}

\usepackage{multirow}

\begin{document}
\graphicspath{{figures/}}
%%%%%%%%%%%%%%%%

% Outcomment only when entries are known. Otherwise leave as is and
%   default values will be used.
%\setcounter{page}{1}
%\VOLUME{00}%
%\NO{0}%
%\MONTH{Xxxxx}% (month or a similar seasonal id)
%\YEAR{0000}% e.g., 2005
%\FIRSTPAGE{000}%
%\LASTPAGE{000}%
%\SHORTYEAR{00}% shortened year (two-digit)
%\ISSUE{0000} %
%\LONGFIRSTPAGE{0001} %
%\DOI{10.1287/xxxx.0000.0000}%

% Author's names for the running heads
% Sample depending on the number of authors;
% \RUNAUTHOR{Jones}
% \RUNAUTHOR{Jones and Wilson}
%\RUNAUTHOR{Zhang, Lu, and Shen}
% \RUNAUTHOR{Jones et al.} % for four or more authors
% Enter authors following the given pattern:
%\RUNAUTHOR{}

% Title or shortened title suitable for running heads. Sample:
% \RUNTITLE{Bundling Information Goods of Decreasing Value}
% Enter the (shortened) title:
	\RUNTITLE{Distributionally Robust Chance-Constrained Building Load Control}
	
	% Full title. Sample:
	% \TITLE{Bundling Information Goods of Decreasing Value}
	% Enter the full title:
	\TITLE{{Building Load Control using Distributionally Robust Chance-Constrained Programs with Right-Hand Side Uncertainty and the Risk-Adjustable Variants}}
	
	% Block of authors and their affiliations starts here:
	% NOTE: Authors with same affiliation, if the order of authors allows, 
	%   should be entered in ONE field, separated by a comma. 
	%   \EMAIL field can be repeated if more than one author
	\ARTICLEAUTHORS{%
		\AUTHOR{Yiling Zhang}
		\AFF{Department of Industrial and Systems Engineering, University of Minnesota \EMAIL{yiling@umn.edu}}
		\AUTHOR{Jin Dong\footnotemark[\value{footnote}]}
		\AFF{Electrification and Energy Infrastructures Division, Oak Ridge National
			Laboratory, \EMAIL{dongj@ornl.gov}}
		% Enter all authors
		
	}% end of the block
	
	\ABSTRACT{%
		Aggregation of heating, ventilation, and air conditioning (HVAC) loads can provide reserves to absorb volatile renewable energy, especially solar photo-voltaic (PV) generation. 
	%		However, the time-varying   PV generation is not perfectly known when the system operator decides the HVAC control schedules. 
	In this paper, we decide HVAC control schedules under uncertain PV generation, using 
	%		To consider the unknown uncertain PV generation, in this paper, we formulate 
	a distributionally robust chance-constrained (DRCC)  building load control model under two typical ambiguity sets: the moment-based and Wasserstein ambiguity sets. We derive mixed integer linear programming (MILP) reformulations for DRCC problems under both sets. Especially, for  the Wasserstein ambiguity set, we utilize the right-hand side (RHS) uncertainty to derive a more compact MILP reformulation than the commonly known MILP reformulations with big-M constants. All the results  also apply to  \textcolor{black}{general individual  chance constraints with RHS uncertainty}.
	Furthermore, we propose an adjustable chance-constrained variant to achieve trade-off between the operational risk and costs. We derive \textcolor{black}{MILP reformulations} under the Wasserstein ambiguity set and second-order conic programming (SOCP) reformulations under the moment-based set. Using real-world data, we conduct computational studies to demonstrate the efficiency of the solution approaches and the effectiveness of the solutions.
		
	}%
	
	% Sample 
	%\KEYWORDS{deterministic inventory theory; infinite linear programming duality; 
	%  existence of optimal policies; semi-Markov decision process; cyclic schedule}
	
	% Fill in data. If unknown, outcomment the field
% 	\HISTORY{}
	
	\maketitle
	%%%%%%%%%%%%%%%%%%%%%%%%%%%%%%%%%%%%%%%%%%%%%%%%%%%%%%%%%%%%%%%%%%%%%%
With growing environmental consciousness and government regulations, renewable energy sources (RESs) are expected to account for 29\% of the total electricity consumption by 2040 \citep{conti2016international}. 
% Higher penetration of renewable energy is challenging, due to potential two-way power flow issues, variability in output and stress on electricity grids' balance, e.g., network frequency and voltage stability \cite{teodorescu2011grid}.
Given that the renewable energy generally cannot adjust their output to reflect changes of demand, higher penetration of renewable energy  may cause electrical supply and demand imbalance issues and can be challenging with variability in output and stress on electricity grids' balance, e.g., network frequency and voltage stability \citep{teodorescu2011grid}.

With the advanced development of smart sensing and control technologies, one solution is utilizing heating, ventilation, and air conditioning (HVAC) systems as grid-responsive flexible load resource, \textcolor{black}{i.e., demand response}.
Given their large amount of power consumption, enormous thermal mass and considerable resistances, the flexible HVAC loads can be employed as virtual storage resources to compensate high frequent fluctuations in renewable energy such as solar photo-voltaic (PV) \citep{yin2016quantifying}. 

In \citet{dong2017adaptive,dong2018model},  they study the problem of using aggregated HVAC systems to absorb solar PV generation. However, the inherent uncertainties of the problem are ignored in their deterministic models, such as uncertainties of the thermal controlled loads (TCLs) and renewable resources, which are mainly determined by factors of weather and consumer behavior \citep{zhang2018data}. The uncertainties can be further intensified by missing samples and low resolution information in HVAC data collection \citep{wijayasekara2015data,vzavcekova2014towards}. In this paper, we take the uncertainties into consideration by employing distributionally robust chance-constrained (DRCC) programs. 

%In this paper, we consider using HVAC systems
%to compensate fluctuations in solar PV generation,
%as previously studied by \cite{dong2017adaptive,jin2018MPC},  which develop deterministic models to track PV generation. However, the inherent uncertainties are ignored in the deterministic models, such as uncertainties of the thermal controlled loads (TCLs) and renewable resources, which are mainly determined by factors of weather and consumer behavior \citep{zhang2018data}. The uncertainties can be further intensified by missing samples and low resolution information in HVAC data collection \citep{wijayasekara2015data,vzavcekova2014towards}.  

% due to potential two-way power flow issues,  
%  Higher penetration of renewable energy can cause the imbalance issues between electrical supply and demand. Because the renewable energy generally cannot adjust their output to reflect changes on the demand side and thus can hardly provide consistent and predictable power generation to the grid. 

\subsection{Relevant Literature}
The aggregated HVAC loads as a virtual storage have become a key player in providing grid \textcolor{black}{demand-responsive} services including load balancing
%The virtual storage concept has been widely studied  and the aggregated HVAC loads have become a key player in providing grid services
\citep[see, e.g.,][]{lu2012evaluation,dong2017adaptive,barooah2019virtual,wang2020operating}. \citet{hao2014aggregate} provide a virtual storage model to characterize aggregate energy flexibility from building loads.
It is followed by virtual storage model identification and flexibility quantification in \citet{hughes2015virtual} and \citet{stinner2016quantifying}, respectively. \textcolor{black}{With the support of home energy management systems (HEMSs), distribution system operators (DSOs) can connect with customers {(e.g., using virtual storage model via aggregated HVAC loads)} to realize system-wide control objectives, e.g., demand response.}
Deterministic optimization models have been proposed to orchestrate the aggregated virtual storage devices \citep[see, e.g.,][]{hao2017optimal,dong2018model}.
% Moreover, various optimization techniques have been proposed in \citet{hao2017optimal,jin2018MPC} to orchestrate the aggregated virtual storage devices. 
\textcolor{black}{In addition to the scalability and privacy issues, residential demand response programs also confront the challenges of handling uncertain parameters, e.g., uncertainties of weather and consumer behavior. To consider the impacts of forecasting errors,}
both robust optimization and stochastic programming techniques have  been introduced to account for modeling  disturbance uncertainties \citep[see, e.g.,][]{chen2012real,nguyen2014risk,zhang2019distributionally,kocaman2020stochastic}. \textcolor{black}{For example, \citet{diekerhof2017hierarchical} schedule flexible devices to reduce peak loads and customers bills under weather and occupancy uncertainties. They propose a distributed robust optimization framework to hedge against the worst-case scenario.
	\citet{lu2020economic} apply robust optimization to maintain thermal comfort in heat and electricity integrated energy systems.  \citet{nguyen2014coordinated} propose a two-stage stochastic program which uses HVAC loads to smooth out the power fluctuation of a wind farm and/or a solar farm. \citet{alhaider2016benders} further take  sizing decisions of PV systems and decisions of battery energy storage systems into account, and formulate a two-stage stochastic integer linear program. } 

{Recently, distributionally robust optimization (DRO) techniques have gained wide interest. Instead of assuming a specific probability distribution of the system uncertainties as in stochastic programming, the DRO approaches consider a family of probability distributions with prior knowledge of the uncertainties, termed as ambiguity set.  The DRO approaches have been applied to  many problems in power systems, such as energy storage operation \citep[see, e.g.,][]{yang2019data}, %congestion line management \citep[see, e.g.,][]{qiu2015distributionally}, economic dispatch \citep[see, e.g.][]{wei2016distributionally}, 
	optimal power flow \citep[see, e.g.,][]{zhang2016distributionally,duan2018distributionally}, and unit commitment \citep[see, e.g.,][]{zhao2018distributionally}.
	Two typical groups of the ambiguity sets employed in DRO are \emph{moment-based} \citep[e.g.,][]{delage2010distributionally} and \emph{distance-based} \citep[e.g., Wasserstein metric][]{esfahani2018data} ambiguity sets.
	%
	%Recently, the distributionally robust optimization (DRO) techniques have been applied to virtual storage. 
	%
	\textcolor{black}{
		For example, considering a moment-based ambiguity set, \citet{zhang2019distributionally} employ a  DRCC program to enable more effective use of uncertain renewables with HVAC systems. 
		A similar formulation has been proposed in \citet{guo2020optimal} to solve optimal pump coordination in water distribution networks under uncertain water demand.  They consider a distributionally robust two-stage stochastic program under a Wasserstein ambiguity set.}
	The stochastic model predictive control \textcolor{black}{(MPC)}, an approach for energy efficiency in HVAC units \citep{dong2018model}, has been considered with distributionally robust chance constraints in \citet{mark2020stochastic}. 
	%\citet{mark2020stochastic} consider stochastic model predictive control with disitributionally robust chance constraints, which can be a control approach for energy efficiency in HVAC units \citep{dong2018model}.
	They use conditional value-at-risk (CVaR) to approximate the chance constraints.  
}

\textcolor{black}{
	%	
	%	Adjustable chance constraint: \citet{shen2014using}: individual chance constraint with right-hand side uncertainty, p-sufficient points
	%%
	%
	Another stream of research relevant to this paper is on adjustable chance constraints. Instead of considering a predetermined (fixed) risk level for chance constraints, decision makers can be interested in balancing between the risk level and operational costs by varying the risk level. The adjustable chance constrained models treat the risk levels as decision variables, which have been applied to various problems, such as \textcolor{black}{metal} melting \citep{evers1967new}, flexible ramping capacity \citep{wang2018adjustable}, power dispatch \citep{qiu2016data,ma2019distributionally}, portfolio optimization \citep{lejeune2016multi}, and humanitarian relief network design \citep{elcci2018chance}.
	\citet{wang2018adjustable} and \citet{evers1967new}  assume that the inverse of the cumulative distribution is known, which, however, is not always accessible. \citet{qiu2016data} approximate the chance constraint by the sample average approximation and transform the problem to a mixed-integer program.
	%  relied on big-M coefficients.
	% , which poses a formidable challenge when solving large instances. 
	%
	When the uncertainty only happens on the right-hand side, assuming a discrete distribution, \citet{shen2014using}  proposes a mixed integer linear programming (MILP) reformulation based on $p$-efficient point using a special ordered set of type 1 (SOS1) constraint. Along the same vein, \citet{elcci2018chance}  propose an alternative MILP using a knapsack inequality which yields \textcolor{black}{an equivalent  linear programming relaxation} as the one in \citet{shen2014using}. All the research above works on individual chance constraint. 
	For joint chance constraints, Bonferroni approximation is one classical approximation, which enforces individual chance constraints with variable risk levels and bounds on the sum of the risk levels. %Equivalently, the risk level of each single chance constraint is treated as a variable.
	\citet{xie2019optimized} study the Bonferroni approximation of distributionally robust joint chance constraints under a moment-based ambiguity set which matches the exact mean and covariance. \citet{ma2019distributionally} apply a similar joint chance constraint, which further requires unimodality, to a power dispatch problem with do-not-exceed limits. }

\subsection{Summary of Main Contributions}

%\citet{zhang2019distributionally} employ distributionally robust chance-constrained  (DRCC) program considering uncertain renewables to enable more-effectively use of the renewable generation under a moment-based ambiguity set.
In this paper, in addition to the moment-based ambiguity set used in the prior work of \citet{zhang2019distributionally}, we further consider DRCC programs under the Wasserstein ambiguity set. In particular, by exploiting the right-hand-side (RHS) uncertainty, we derive an MILP reformulation based on the conditional value-at-risk (CVaR) interpretation for DR chance constraints pointed out by \citet{xie2018distributionally, chen2018data}. 
\textcolor{black}{Recently, based on the CVaR (primal) interpretation,  \citet{ho2020distributionally}  provide an MILP reformulations for joint DRCC programs with RHS uncertainty.}
%	Different from their work, our results are derived from the dual perspective of CVaR.
\textcolor{black}{From the primal perspective, CVaR is a conditional expectation at the tail of a distribution; while the dual representation further indicates that CVaR is a weighted sum of the least favorable outcomes. Our results are built based on the dual perspective by deriving the weights of the least favorable outcomes.}
%In this paper, we extend our prior work of \citet{zhang2019distributionally} to Wasserstein ambiguity set.
%The results can be applied to general individual \textcolor{blue}{binary} chance constraints with RHS uncertainty (Regarding improved formulations and valid inequalities for joint distributionally robust chance constraints with RHS uncertainty, we refer the readers to \citet{ho2020distributionally}.).
Moreover, to better balance the operational cost and PV utilization, we propose an adjustable chance-constrained formulation that treats the risk level of chance constraints as a decision variable rather than a given (fixed) parameter in the DRCC formulations.
%   We derive MILP reformulations for the Wasserstein set and solve for the moment-based ambiguity set using two second-order cone programs (SOCPs).
%
%
% 
% However, low resolution information and missing samples are common problems in HVAC data collection \citep{wijayasekara2015data,vzavcekova2014towards}. (ADD REFERENCE: review emilio Maddalena) yet the abundance of measuring instruments has rather
% been regarded as an opportunity to develop data-driven solutions to
% existing problems in the area
%
%
% In \cite{dong2017adaptive,jin2018MPC}, they focus on minimizing the difference between total power consumption of HVAC units and PV signal plus state deviation, while assuming that PV generation is perfectly known or predictable. 
% Such formulation can be further enhanced with energy storage systems to provide better building-to-grid (B2G) services (e.g., %\cite{razmara2017enabling,
% \cite{allen2018supervisory,8431619}).  
% Although utilizing building loads for regulation services has been an important research area for years, the above deterministic solutions are not able to handle uncertainties associated with TCL models and intermittent renewable resources, which are ultimately shaped by outside factors such as weather and consumer behavior \cite{zhang2018data}.
% 
% 
We summarize our contributions as follows:
\begin{enumerate}
	\item We formulate the building load control (BLC) problem of HVAC units using DRCC optimization under two types of ambiguity sets: moment-based and Wasserstein ambiguity sets in Section \ref{sec:DRCC}. We also propose their variants with adjustable chance constraints to balance operational cost and performance \textcolor{black}{in Section \ref{sec:adjustable-cc}}.
	
	\item We provide exact reformulations for the DRCC and the adjustable DRCC formulations \textcolor{black}{(with binary decision variables)}. %, which can also be applied to general individual \textcolor{blue}{binary} chance constraints with RHS uncertainty.
	\textcolor{black}{Specifically, under the Wasserstein ambiguity set, we derive exact MILP reformulations for both  DRCC and its adjustable variant. Under the moment-based ambiguity set, DRCC yields an MILP reformulation, while solving the adjustable DRCC is equivalent to solving two second-order conic programs (SOCPs). \textcolor{black}{All the results for DRCC models hold for \emph{general} individual DR chance constraints with RHS uncertainty (even with continuous decision variables), while the results for the adjustable variants hold  for general adjustable individual \emph{binary}  DR chance constraint with RHS uncertainty.}}
	
	\item We conduct computational tests on various instances and demonstrate the efficiency and effectiveness of the proposed reformulations via real-world data \textcolor{black}{in Sections \ref{sec:comp-setup}--\ref{sec:comp-adjustable}}.
\end{enumerate}

%The reformulation of DRCC under the Wasserstein ambiguity set is  based on the CVaR interpretation and quantile strengthening. 
% \citet{ho2020distributionally} also exploit the CVaR 
% interpretation and the quantile strengthening to derive an MILP reformulation for joint DR chance constraint with RHS uncertainty over Wasserstein ambiguity set. Our formulation is derived from a dual perspective.  

The remainder of the paper is organized as follows. Section \ref{sec:modeling} presents the mathematical formulations of the deterministic and stochastic chance-constrained BLC problems. \textcolor{black}{In Sections \ref{sec:DRCC}--\ref{sec:adjustable-cc}, we present the DRCC models and their adjustable variants, and derive exact reformulations under both the moment-based and Wasserstein ambiguity sets.
	% In Section XXX, we derive the reformulations of the DRCC models.
	%In Section \ref{sec:adjustable-cc}, we further extend the DRCC models under both ambiguity sets by considering adjustable chance constraints. We derive their tractable reformulations. 
	In Sections \ref{sec:comp-setup}--\ref{sec:comp-adjustable}, we conduct extensive numerical studies on both non-adjustable DRCC models and adjustable chance-constrained models. Finally, we draw conclusions in Section \ref{sec:conclusions}. }

% \cite{hao2015potentials}).

% However, the introduction of smart sensing, metering, and control technology enables grid-responsive building load. Specifically, the flexible building loads can be used as virtual storage resources to compensate fluctuations in renewable energy such as solar photovoltaic and thus integrate more distributed energy resources into the grid. For the building owners, the load flexibility generates a new value stream for them by more efficient use of electricity. Research shows that by 8\% reduction from only water heating and air conditioning, which still maintains comfort and service quality, results in 10\% - 40 \% on customer bills and 13 billion bill savings per year for grid costs.

% One solution is utilizing flexibility in heating, ventilation, and air conditioning
% (HVAC) systems, a type of thermostatically controlled loads,
% TCLs, providing a great flexible resource given their large
% amount of power consumption, enormous thermal storage and
% considerable resistances. Performance and effects of demand
% response using HVAC loads have been widely studied (e.g., \cite{nguyen2014joint}).
% Moreover, HVAC systems can serve as means to provide
% ancillary services (see %\cite{hao2014ancillary,hao2015potentials,balandat2014contract}).
% \cite{hao2015potentials}).

\section{Model Formulation}

% \section{DRCC Models}
\label{sec:modeling}

The BLC problem utilizes an aggregated HVAC load of $N_{\text{HVAC}}$ units, e.g., buildings, to  absorb the solar PV generation (collected from $N_\text{PV}$ PV panels) locally while delicately maintaining desired indoor temperature for each unit throughout the day. We discretize the day-time duration into $N_p$  periods with a time interval of $\Delta t$ (e.g., 10 minutes). For each period $t = 1,\ldots,N_p$, denote a binary decision variable $u_{t,\ell}\in \{0,1\}$ for HVAC unit $\ell$  to indicate its scheduled mode: if $u_{t,\ell} = 1$, ON; otherwise $u_{t,\ell} = 0$, OFF.

% In the system, we aim to utilize a fleet of $N_{\text{HVAC}}$ HVAC units to consume PV power generation {\color{black}(collected from $N_\text{PV}$ PV panels)} with a desired probability during day time of a single day. We divide the day-time duration into $N_p$ periods. At period $t = 1,\ldots,N_p$, we  decide for each HVAC unit $j = 1,\ldots, N_{\text{HVAC}}$ whether in ON or OFF mode, denoted by a binary variable $u_{t,j} \in \{ 0,1\}$, equal to  $1$ if ON, and $0$ otherwise.  
For each HVAC unit, to characterize the dynamics of room temperature and outdoor temperature,
we consider a widely used building thermal model \citep[see e.g.,][]{mathieu2013state}, where
the system state is the room  temperature $T$, the system input is HVAC ON/OFF status, $Mode_\text{HVAC}$, and the system disturbances include outdoor temperature $T_\text{out}$ and solar irradiance $Q_\text{out}$. Based on the building thermal model, a continuous-time linear time invariant (LTI) system for each HVAC unit in the state-space form is as follows. \begin{equation}\label{eq:LTI}
	\dot{T} = \frac{1}{RC} T_\text{out} + \frac{1}{RC} T + \frac{1}{C}Q_\text{out} + \frac{Mode_\text{HVAC}}{C}Q_\text{HVAC},
\end{equation}
where $R$ is the building's thermal resistance, $C$ is the building's thermal capacity, and $Q_\text{HVAC}$ is  cooling capacity of the building. \textcolor{black}{These parameters can be estimated following standard building energy modeling techniques,  Resistance-Capacitance (RC) model \citep{belic2016thermal, cui2019hybrid}, in particular. Such continuous-time LTI model is further converted into a discrete-time model using various techniques including the Zero-Order Hold (ZOH) method.} In this paper, for \textcolor{black}{each individual} building $\ell=1,\ldots,N_{\text{HVAC}}$, we consider a discrete-time building thermal model with a sampling interval of $\Delta t$ as
\begin{equation}\label{eq:x-function}
	x_{t,\ell} = A_\ell x_{t-1,\ell} + B_\ell u_{t,\ell} + G_\ell v_\ell,
\end{equation}
where $x_{t,\ell}$ is the room temperature of period $t$, $u_{t,\ell}$ is the binary mode decision variable, $v_\ell$ is the system disturbance, and the parameters $A_\ell, B_\ell, G_\ell$ can be computed from the continuous-time model \eqref{eq:LTI}.

In this study, we focus on the single-period models, where we solve for the optimal ON/OFF mode decisions at the current period given an initial room temperature resulted from previous periods' decisions. To solve the control problem for all $N_p$ periods, we sequentially solve $N_p$ small optimization problems, each for one period.
The problem can also be formulated as a multi-period model, which solves one monolithic  optimization  problem  for all  the  mode  decision  variables over $N_p$ periods. In our preliminary results of \citet{zhang2019distributionally}, via extensive computational studies, 
%we compare the performance of the multi-period and the single-period models using DRCC approaches under the moment-based ambiguity set  
we demonstrate that the multi-period models  have similar  solutions as those of single-period models \textcolor{black}{using DRCC approaches under the moment-based ambiguity set}. However, the multi-period models can suffer from computational difficulty  especially for large-sized problems. Therefore, in this paper, we focus on the single-period models. For the rest of this section, we present two optimization formulations of the single-period models: a deterministic formulation, which assumes that the solar PV generation is deterministic and perfectly known, and a chance-constrained formulation, which assumes the PV generation is stochastic.

\subsection{Deterministic Formulation}
\label{sec:determ}

At period $t$, we denote $x_t = (x_{t,\ell},\ \ell=1,\ldots,N_{\text{HVAC}})^\top$ the \textcolor{black}{ auxiliary decision vector} of room temperature of $N_{\text{HVAC}}$ buildings, and denote $u_t = (u_{t,\ell}, \ \ell=1,\ldots,N_{\text{HVAC}})^\top$ the mode vector of ON/OFF decisions of $N_{\text{HVAC}}$ buildings. \textcolor{black}{For  building $\ell$ in the ON mode, the energy consumption of one period is $P_\ell$.}
% In the deterministic model, we assume that PV power output at period $t$ is known perfectly as $P_{\text{PV},t} \in \mathbb R^{N_\text{PV}}$. Given the initial room temperature $x_{t-1,j}$, the MILP formulation is as follows.
We assume that the solar PV generation at period $t$ is perfectly known as $P_{\text{PV},t} \in \mathbb R^{N_\text{PV}}$. At time $t$,  given an initial room temperature $x_{t-1,\ell}, \ \ell = 1,\ldots,N_{\text{HVAC}}$, the BLC problem can be formulated as the following  MILP \citep{dong2018model,zhang2019distributionally}.
\begin{subequations}
	\label{eq:determ}
	\begin{eqnarray}
		% \textbf{\text{Determ}$(t,x_{t-1},P_\text{PV})$}\quad
		\label{eq:determ-obj}\min_{\color{black}{u_t, \eta, \beta_{t,\ell}, x_{t,\ell}}} &&    c_\text{sys}\sum_{\ell=1}^{N_{\text{HVAC}}}\beta_{t,\ell}+ c_\text{switch}\sum_{\ell=1}^{N_{\text{HVAC}}}u_{t,\ell} +c_\text{PV}\eta\\
		\mbox{s.t.} && \eqref{eq:x-function} \nonumber \\
		&& -\eta \le \sum_{\ell=1}^{N_{\text{HVAC}}} P_\ell u_{t,\ell} - \sum_{i=1}^{N_\text{PV}} P_{\text{PV},t,i} \le \eta\\  
		\label{eq:absolute}&& -\beta_{t,\ell} \le x_{t,\ell} - x_\text{ref} \le \beta_{t,\ell}, \ \ell =1,\ldots,N_{\text{HVAC}}\\
		\label{eq:determ-x-bounds} && x_\text{min} \le x_{t,\ell} \le x_\text{max},\ \ell = 1,\ldots,N_{\text{HVAC}}\\
		\label{eq:determ-binary}&& u_t \in \{0,1\}^{N_{\text{HVAC}}},
	\end{eqnarray}
\end{subequations}
\textcolor{black}{where, in the objective, auxiliary variable $\beta_{t,\ell} = |x_{t,\ell} - x_\text{ref}|$ denotes the room temperature deviation from the set-point $x_\text{ref}$ and $\eta = |\sum_{\ell=1}^{N_{\text{HVAC}}} P_\ell u_{t,\ell} - \sum_{i=1}^{N_\text{PV}} P_{\text{PV},t,i}|$ denotes the signal deviation between the total control signal (of all $N_\text{HVAC}$ buildings) and the total PV signal (of all $N_\text{PV}$ PV panels).}
%where $|\cdot|$ takes the absolute value of a real number. 
The objective \eqref{eq:determ-obj} minimizes the total penalty cost of i) discomfort (indicated by the room temperature deviation), ii) switching cycles, and iii) PV tracking error, with unit cost parameters $c_\text{sys}, c_\text{switch}, c_\text{PV}$. Constraints \eqref{eq:determ-x-bounds} require that the room temperature $x_{t,\ell}$ is maintained in the comfort band $[x_\text{min}, x_\text{max}]$ for building $\ell=1,\ldots,N_{\text{HVAC}}$. The last constraint \eqref{eq:determ-binary} enforces binary decision $u_t$.

% An accurate estimate of the solar PV output $P_{\text{PV},t}$ is critical to the deterministic formulation
% \eqref{eq:determ}. However, in  practice, a good prediction of the PV output may not be available and can thus be uncertain due to the fluctuating nature of solar energy which is  introduced by cloud shadows, wind speed, and other factors.
% % In practice, an accurate estimate of the solar PV output may not be available due to many weather factors such as 
% % When the solar PV output $\widetilde{P}_{\text{PV},t}\in\mathbb R^{N_\text{PV}}$ is not accurately estimated, 
% We employ the Monte-Carlo sampling approach to generate $M$ i.i.d. scenarios of the uncertain PV output, denoted by $P_{\text{PV},t}^1,\ldots, P_{\text{PV},t}^M$. Each scenario $P_{\text{PV},t}^m$ is associated with a probability $p_{t,m}\ge 0$, such that $\sum_{m=1}^M p_{t,m}=1$. The stochastic programming formulation, via
%  the sample average approximation approach \cite{kleywegt2002sample}, is formulated as  follows.

%  \begin{eqnarray}
% \label{eq:stoch}
% && \min_{u_t} \left\{   c_\text{sys}\sum_{j=1}^{N_{\text{HVAC}}}|x_{t,j}-x_\text{ref}|+ c_\text{switch} \sum_{j=1}^{N_{\text{HVAC}}}u_{t,j} + c_\text{PV}\sum_{m=1}^{M} p_{t,m}\left| \sum_{j=1}^{N_{\text{HVAC}}}P_j u_{t,j} - \sum_{i=1}^{N_\text{PV}}P_{\text{PV},t,i}^m\right|: \ \eqref{eq:determ-x-bounds}-\eqref{eq:determ-binary} \right\}.\hspace{3mm}%\nonumber
% \end{eqnarray}

\subsection{Chance-Constrained Formulation}
\label{sec:cc}

An accurate prediction of the solar PV output $P_{\text{PV},t}$ is critical to the performance of the deterministic formulation
\eqref{eq:determ}. However, in  practice, a good prediction of the PV output may not be available,  due to the fluctuating nature of solar energy, which can be  introduced by cloud shadows, wind speed, and other factors, and  thus can be uncertain. In this section, we introduce a chance-constrained formulation to take into account uncertain PV output.

% In the last term of the objective function \eqref{eq:determ-obj}, we penalize the signal deviation  with cost parameter $c_\text{PV}$ to ensure most PV generation can be consumed locally by the HVAC fleet. 
Instead of penalizing the signal deviation $\eta$ in the objective to enforce all PV generation being consumed,  the chance-constrained formulation \eqref{eq:cc} employs a soft constraint \eqref{eq:cc-cc}, the chance constraint to ensures that, with a high probability, the solar PV output is consumed by the HVAC fleet. The chance-constrained formulation is 
% In the deterministic model, we assume that the PV power output can be perfectly known ahead of time, which does not always hold since the PV power is uncertain due to many random weather factors such as outdoor temperature, cloud shadows, wind speed, changes in composition of the clear
% atmosphere (e.g., dust, smoke, humidity), etc. In the following section, we present a stochastic programming formulation which is a classic method  used to model system uncertainties.
\begin{subequations}
	\label{eq:cc}
	\begin{eqnarray}
		\label{eq:cc-obj} \min_{u_t,  \beta_{t,\ell}, x_{t,\ell}} && c_\text{sys}\sum_{\ell=1}^{N_{\text{HVAC}}}\beta_{t,\ell} + c_\text{switch}\sum_{\ell=1}^{N_{\text{HVAC}}}u_{t,\ell}\\
		%\left( x - x_{\text{ref}}\right)^{\top}Q\left(x  - x_{\text{ref}}\right)\\
		%\min_{u,z}&& z \\
		\mbox{s.t.}
		%&& z \ge \left( x - x_{\text{ref}}\right)^{\top}Q\left(x  - x_{\text{ref}}\right) \\
		\label{eq:cc-cc} &&  \mathbb P\left( \sum_{\ell=1}^{N_{\text{HVAC}}} P_\ell u_{t,\ell} - \sum_{i=1}^{N_\text{PV}}\widetilde{P}_{\text{PV},t,i} \ge 0\right)  \ge 1-\alpha_t \hspace{7mm}\\
		&& \eqref{eq:x-function}, \ \eqref{eq:absolute}-\eqref{eq:determ-binary}, \nonumber
	\end{eqnarray}
\end{subequations}
where $\widetilde{P}_{\text{PV},t,i}$ denotes the uncertain PV generation at period $t$ of panel $i$.
Constraint \eqref{eq:cc-cc} ensures that the PV generation is absorbed by the HVAC fleet with probability $1-\alpha_t$. The risk level $\alpha_t \in (0,1)$ is  predefined and reflects  system operator's risk preference, usually a small number.

To solve the chance-constrained model, we employ the Sample Average Approximation (SAA) approach \citep[see, e.g.,][]{luedtke2008sample} to derive bounds and obtain feasible solutions. Using the Monte Carlo sampling method, we generate a set of finite samples of the uncertainty $\widetilde{P}_{\text{PV},t}$. 
and enforce $ \sum_{\ell=1}^{N_{\text{HVAC}}} P_\ell u_{t,\ell} - \sum_{i=1}^{N_\text{PV}}\widetilde{P}_{\text{PV},t,i} \ge 0$ for sufficiently many samples.

Specifically, we generate
$N$ i.i.d. scenarios of the uncertain PV output $\widetilde{P}_{\text{PV},t}$, denoted by $P_{\text{PV},t}^1,\ldots, P_{\text{PV},t}^N$. Each scenario $P_{\text{PV},t}^n$ is associated with a probability $p_{t,n}\ge 0$, such that $\sum_{n=1}^N p_{t,n}=1$. For each scenario $n$, we associate a binary variable $\rho_n$ such that $\rho_n = 0$ indicates 
\begin{equation}\label{eq:cc-inner} \sum_{\ell=1}^{N_{\text{HVAC}}} P_\ell u_{t,\ell} - \sum_{i=1}^{N_\text{PV}}{P}^n_{\text{PV},t,i} \ge 0\end{equation} 
and when $\rho_n = 1$, constraint \eqref{eq:cc-inner}  \textcolor{black}{is relaxed and can} be violated. The chance constraint \eqref{eq:cc-cc} is approximated \textcolor{black}{\citep[see, e.g., ][]{ruszczynski2002probabilistic}} by
\begin{subequations}
	\begin{eqnarray}
		\label{eq:cc-milp1}&&\sum_{\ell=1}^{N_{\text{HVAC}}} P_\ell u_{t,\ell} - \sum_{i=1}^{N_\text{PV}}{P}^n_{\text{PV},t,i} \ge -M \rho_n, \ n =1,\ldots,N\\
		\label{eq:cc-milp2}&&\sum_{n=1}^N \textcolor{black}{p_n}\rho_n  \le \alpha_t\\
		\label{eq:cc-milp3}&&\rho_n \in \{0,1\}, \ n = 1,\ldots,N,
	\end{eqnarray}
\end{subequations}
where $M$ is a big-M coefficient. Constraint \eqref{eq:cc-milp2} ensures that the probability of violating \eqref{eq:cc-inner} is no more than $\alpha_t$.  By replacing the chance constraint \eqref{eq:cc-cc} with \eqref{eq:cc-milp1}--\eqref{eq:cc-milp3}, we obtain an MILP approximation of the chance-constrained model \eqref{eq:cc}.

% Note that constraint \eqref{eq:cc-cc} allows the constraint $ \sum_{j=1}^{N_{\text{HVAC}}} P_ju_{t,j} - \sum_{i=1}^{N_\text{PV}}\widetilde{P}_{\text{PV},t,i} \ge 0$ to be violated sometimes. 
% In other words, constraint $ \sum_{j=1}^{N_{\text{HVAC}}} P_ju_{t,j} - \sum_{i=1}^{N_\text{PV}}\widetilde{P}_{\text{PV},t,i} \ge 0$ is not satisfied always. the PV generation can be completely absorbed by controllable HVAC loads locally instead of flowing back to the distribution grid with at least probability $1-\alpha_t$.
% That is, the difference of the HVAC consumption and the PV generation can be positive or negative, which aligns with the deterministic and stochastic formulations.

\section{DRCC Formulations}
\label{sec:DRCC}

In the stochastic chance-constrained formulation \eqref{eq:cc}, full knowledge of the PV's probability distribution is required. However, an accurate probability distribution can be challenging to obtain especially when the underlying distribution (while ambiguous) is time-varying. As a consequence, the solution obtained from the chance-constrained model might be sensitive to the choice of probability distribution and thus results in poor performance.
This phenomenon  is called the
optimizer's curse \citep{smith2006optimizer} of solving stochastic programs.
% In the  stochastic chance-constrained formulation \eqref{eq:cc}, we require full knowledge of the PV's probability distribution, which, however, is often challenge to obtain. On the other hand, optimal solutions might be  sensitive to the choice of probability distribution. An inaccurate estimate can lead to poor performance, i.e., low probability of absorbing the PV output. 
To address the curse, a natural way is to employ a set of plausible probability distributions, denoted as $\mathcal D_t$, rather than assuming a specific probability distribution. Specifically,  we consider the DRCC formulation as follows
% In this section, we consider the DRCC formulation \eqref{eq:drcc} which does not assume a specific probability distribution but a set of plausible distributions, denoted by $\mathcal D_t$.
%
\begin{subequations}
	\label{eq:drcc}
	\begin{eqnarray}
		\label{eq:drcc-obj} \min_{u_t, \beta_{t,\ell}, x_{t,\ell}} && c_\text{sys}\sum_{\ell=1}^{N_{\text{HVAC}}}\beta_{t,\ell} + c_\text{switch}\sum_{\ell=1}^{N_{\text{HVAC}}}u_{t,\ell}\\
		%\left( x - x_{\text{ref}}\right)^{\top}Q\left(x  - x_{\text{ref}}\right)\\
		%\min_{u,z}&& z \\
		\mbox{s.t.} %&& z \ge \left( x - x_{\text{ref}}\right)^{\top}Q\left(x  - x_{\text{ref}}\right) \\
		\label{eq:drcc-cc} && \inf_{f \in\mathcal D_t} \mathbb P\left( \sum_{\ell=1}^{N_{\text{HVAC}}} P_\ell u_{t,\ell} - \sum_{i=1}^{N_\text{PV}}\widetilde{P}_{\text{PV},t,i} \ge 0\right)  \ge 1-\alpha_t \hspace{7mm}\\
		&& \eqref{eq:x-function}, \ \eqref{eq:absolute}-\eqref{eq:determ-binary}. \nonumber
		%&& x = \frac{1}{RC} T_\text{out} - \frac{1}{RC} x_0
		%+ \frac{1}{C} Q_\text{out} + \frac{u}{C} Q_\text{HVAC}\\
		%  && x_\text{min} \le x_k \le x_\text{max},\ k=1,\ldots,N_{\text{HVAC}}\\
		%  && u \in \{0,1\}^{N_{\text{HVAC}}}
	\end{eqnarray}
\end{subequations}
Constraint \eqref{eq:drcc-cc} ensures that the probability of absorbing the PV generation locally by the HVAC fleet is guaranteed at least $1-\alpha_t$ for any probability distribution $f\in \mathcal D_t$. That is, for all probability distributions in  $D_t$, the worst-case probability of coordinating the HVAC fleet to consume the PV generation is no less than $1-\alpha_t$. We note that constraint \eqref{eq:drcc-cc} is an individual DR chance constraint with RHS  uncertainty.

\subsection{Ambiguity Sets}

One critical question of the DRCC formulation is how to choose the ambiguity set $\mathcal D_t$. A good choice of $\mathcal D_t$ should take into account the characteristics of the underlying probability distribution and the tractability of the DRCC formulation.
% The choice of the ambiguity set 
% The ambiguity set is the key ingredient to the DRCC model which characterizes the candidate distributions that can be the true underlying distribution of the PV output. 
Two types of ambiguity sets have been widely studied: (i) moment-based  and (ii) distance-based ambiguity sets. In this paper, we consider a moment-based ambiguity set containing moment constraints on the first- and second-order moments \citep[see, e.g.,][]{delage2010distributionally} and a distance-based ambiguity set using Wasserstein metric \citep[e.g.,][]{esfahani2018data}. 
% Let $\widetilde{P}_{\text{PV},t}\in\mathbb R^{N_\text{PV}}$ denote the vector of $\widetilde{P}_{\text{PV},t,i}, \ i = 1,\ldots,N_\text{PV} $. 
Given a series of independent samples, $\{{P}_{\text{PV},t}^n\}_{n=1}^N$, sampled from the true underlying distribution of the PV generation,  we consider the following two distributional ambiguity sets. 

\paragraph{(i) Moment-based Ambiguity Set.}
The empirical mean and covariance matrix can be calculated as $\mu_t = \frac{1}{N}\sum_{n=1}^N {P}_{\text{PV},t}^n$, $\Sigma_t = \frac{1}{N}\sum_{n=1}^N ({P}_{\text{PV},t}^n - \mu_t)({P}_{\text{PV},t}^n - \mu_t)^\top$. 
The ambiguity set based on the two moment estimates $\mu_t$, $\Sigma_t$,  first  proposed by \citet{delage2010distributionally}, is as follows.
\begin{equation*}
	% \scriptstyle
	\mathcal{D}_t^1
	% (\mu_t,\sigma_t,\gamma_1,\gamma_2) 
	\ = \ \left\{f :
	% \in {\color{black}\mathcal{P}(\mathbb{R}^{N_\text{PV}})}:
	\
	\begin{array}{l} 
		\mathbb P_f(\widetilde{P}_{\text{PV},t} \in \mathbb R^{N_\text{PV}}) = 1\\
		(\mathbb{E}_f[\widetilde{P}_{\text{PV},t}] - \mu_t)^{\top} \Sigma_t^{-1} (\mathbb{E}_f[\widetilde{P}_{\text{PV},t}] - \mu_t) \ \leq \ \gamma_1, \\%[0.3cm]
		\mathbb{E}_f\bigl[(\widetilde{P}_{\text{PV},t} - \mu_t) (\widetilde{P}_{\text{PV},t} - \mu_t)^{\top}\bigr] \ \preceq \ \gamma_2 \Sigma_t, 
	\end{array}
	\right\},
\end{equation*}
where $\gamma_1 \ge 0$ and $\gamma_2 \ge \max\{\gamma_1,1\}$. The three constraints guarantee that (1) the true mean of $\widetilde{P}_{\text{PV},t}$ lies in an ellipsoid centered at $\mu_t$ and (2) the true covariance of $\widetilde{P}_{\text{PV},t}$ is bounded above by $\gamma_2 \Sigma_t$.
% Here, $\mu_t$ and $\Sigma_t$ are estimates of the mean and coveriance of $\widetilde{P}_{\text{PV},t}$, which can be empirical mean and covariance matrix, as
% $\mu_t = \frac{1}{N}\sum_{n=1}^N {P}_{\text{PV},t}^n$ and $\Sigma_t = \frac{1}{N}\sum_{n=1}^N ({P}_{\text{PV},t}^n - \mu_t)({P}_{\text{PV},t}^n - \mu_t)^\top$. 
The two parameters $\gamma_1$ and $\gamma_2$ reflect the system operator's tolerance of the moment and distributional ambiguity: the larger the two parameters are, the more the tolerance towards ambiguity and  the  more robustness of the optimal solutions are. The values of $\gamma_1$ and $\gamma_2$   depend on the samples size, support size, and confidence level \citep[See more details in Definition 2 in][]{delage2010distributionally}.

\paragraph{(ii) {Wasserstein Ambiguity Set}.}
Given a positive radius $\delta_t>0$, the Wasserstein ambiguity set defines a ball around the discrete empirical distribution based on the $N$ samples, 
$
\mathbb P_{\widetilde{P}^N_{\text{PV},t}}\left[\widetilde{P}_{\text{PV},t} = P_{\text{PV},t}^n \right] = {{1}/{N}},
$
in the space of probability distributions as follows. 
% Given a positive $\delta_t$, Wasserstein ambiguity set is
% for time period $t$ of the total PV output $\widetilde{P}_{\text{total},t} = \sum_{i=1}^{N_\text{PV}}\widetilde{P}_{\text{PV},t,i}$ is
\begin{equation*}
	\mathcal D_t^2 = \left\{ f: \mathbb P_f\{\widetilde{P}_{\text{PV}, t} \in\mathbb R^{N_\text{PV}} \} = 1,\ W(\mathbb P_f, \mathbb P^N_{\widetilde{{P}}_{\text{PV},t}})\le \delta_t
	\right\},
\end{equation*}
where the Wasserstein distance is defined as 
\begin{equation*}
	W(\mathbb P_1, \mathbb P_2) = \inf_{\mathbb Q} \left\{ \int_{\mathbb R^{N_\text{PV}}\times \mathbb R^{N_\text{PV}}} \|{P}_{\text{PV},t}^1 - {P}_{\text{PV},t}^2\|\mathbb Q(\mathrm{d}{P}_{\text{PV},t}^1, \mathrm{d}{P}_{\text{PV},t}^2): \  \parbox{0.3\textwidth}{%
		$\mathbb Q$ is a joint distribution of $\widetilde{P}_{\text{PV},t}^1$ and $\widetilde{P}_{\text{PV},t}^2$ with marginals $\mathbb P_1$ and $\mathbb P_2$, respectively
	} \right\}.
\end{equation*}
The radius of the ambiguity set controls the degree of the conservatism of the DRCC model. If we set $\delta_t = 0$, the ambiguity set $\mathcal D_t^2$ only contains the empirical distribution and we can recover a chance-constrained model.

\subsection{DRCC Reformulation under Moment-based Ambiguity Set $\mathcal D_t^1$}
\label{sec:drcc-ref-moment}
% 
%In this section, we present an MILP reformulation of the DRCC model \eqref{eq:drcc} under  the moment-based ambiguity set $\mathcal D_t = \mathcal D_t^1$.
Let $\theta_t = \mathbf{1}^\top \mu_t$ and $\sigma_t = \mathbf{1}^\top \Sigma_t \mathbf{1}$, where $\mathbf{1}\in\mathbb{R}^{N_\text{PV}}$ is a vector with all ones.
We follow \citet{zhang2018ambiguous} to rewrite the DR chance constraint \eqref{eq:drcc-cc} under the moment-based ambiguity set $\mathcal D_t^1$ as a linear constraint.
%
%According to Theorem 3.2 in \citet{zhang2018ambiguous}, under the ambiguity set $\mathcal D_t^1$, the DR chance constraint \eqref{eq:drcc-cc} is equivalent to a linear constraint as
%
\textcolor{black}{\begin{proposition}[Adapted from Theorem 3.2 of \citet{zhang2018ambiguous}]
		The  DR chance constraint \eqref{eq:drcc-cc} under $\mathcal D_t = \mathcal D_t^1$ is equivalent to
		\begin{equation}\label{eq:linear-eq}
			\sum_{\ell=1}^{N_{\text{HVAC}}}P_\ell u_{t,\ell} \ge \theta_t + \Omega_t \sigma_t, \text{ 	where  $\Omega_t = \left\{ \begin{array}{ll}
					\sqrt{\gamma_1} + \sqrt{{(1-\alpha_t)}(\gamma_2 - \gamma_1)/{\alpha_t} }, & {\gamma_1}/{\gamma_2} \le \alpha_t\\
					\sqrt{{\gamma_2}/{\alpha_t}}, & {\gamma_1}/{\gamma_2} > \alpha_t.
				\end{array}\right.$}
		\end{equation}
		%	where  $\Omega_t = \left\{ \begin{array}{ll}
		%		\sqrt{\gamma_1} + \sqrt{{(1-\alpha_t)}(\gamma_2 - \gamma_1)/{\alpha_t} }, & {\gamma_1}/{\gamma_2} \le \alpha_t\\
		%		\sqrt{{\gamma_2}/{\alpha_t}}, & {\gamma_1}/{\gamma_2} > \alpha_t
		%	\end{array}\right.$. 
	\end{proposition}
	Instead of enforcing the DR chance constraint with the RHS uncertainty of PV generation, constraint \eqref{eq:linear-eq} requires the total HVAC load no less than the nominal PV generation $\theta_t$ plus the product of $\Omega_t$ and its standard deviation $\sigma_t$, where the value of $\Omega_t$ is determined by the relationship of the ambiguity set parameters $\gamma$'s and the risk level $\alpha_t$.}
Then the DRCC formulation \eqref{eq:drcc} under $\mathcal D_t^1$ is equivalent to the following  MILP problem.
\begin{equation}\label{eq:lp}
	\min_{u_t, \beta_{t,\ell}, x_{t,\ell}}\left\{ c_\text{sys}\sum_{\ell=1}^{N_{\text{HVAC}}}\beta_{t,\ell} + c_\text{switch}\sum_{\ell=1}^{N_{\text{HVAC}}}u_{t,\ell}:   \ \eqref{eq:x-function}, \ \eqref{eq:absolute}-\eqref{eq:determ-binary}, \ \eqref{eq:linear-eq} \right\}.
	\nonumber
\end{equation}

\subsection{DRCC Reformulations under Wasserstein Ambiguity Set $\mathcal D_t^2$}
\label{sec:drcc-ref-wasserstein}

% \paragraph{(ii) Wasserstein Ambiguity Set $\mathcal D_t = \mathcal D^2_t$.}
% 
Denote the total PV output ${P}^n_{\text{total},t} = \mathbf{1}^\top P^n_{\text{PV},t}$.
According to Corollary 2 in \citet{xie2018distributionally} \citep[Theorem 3 in][]{chen2018data}, the DR chance constraint \eqref{eq:drcc-cc} under the Wasserstein ambiguity set $\mathcal D_t^2$ is feasible if and only if the following constraints, \textcolor{black}{with auxiliary variables $\gamma$ and $z_n, \ n=1,\ldots, N$,} are satisfied
% According to Corollary 2 in \citet{xie2018distributionally}, the feasible region described by the DR chance constraint \eqref{eq:drcc-cc} is 
\begin{subequations}
	\label{eq:z}
	\begin{eqnarray}
		% Z =\bigg\{ u_{t}:
		\label{eq:z-1}&& \delta_t  - \alpha_t \gamma \le \frac{1}{N}\sum_{n = 1}^N z_n\\
		\label{eq:z-2}&& -\max\left[\sum_{\ell=1}^{N_{\text{HVAC}}}P_\ell u_{t,\ell}  - P_{\text{total},t}^n, \ 0 \right] \le -z_n - \gamma, \ n=1,\ldots,N \hspace{3mm}\\
		\label{eq:z-3}&& z_n \le 0, \ n = 1,\ldots,N\\
		\label{eq:z-4}&&\gamma \ge 0.
		% \bigg\}.
	\end{eqnarray}
\end{subequations}
\textcolor{black}{We remark that the reformulation \eqref{eq:z} admits a CVaR interpretation \citep{xie2018distributionally}, i.e., \begin{equation}\label{eq:cvar}
		\frac{\delta}{\alpha_t} + \text{CVaR}_{1-\alpha_t}\left[ -\max\left\{\sum_{\ell=1}^{N_\text{HVAC}} P_\ell u_\ell - \widetilde{P}_\text{PV}, 0 \right\} \right] \le 0,
	\end{equation} 
	where $\text{CVaR}_{1-\alpha_t}\left[-\max\left\{\sum_{\ell=1}^{N_\text{HVAC}} P_\ell u_\ell - \widetilde{P}_\text{PV}, 0 \right\} \right] = \min_\gamma \left\{ \gamma + \frac{1}{\alpha_t} \mathbb E_{\mathbb P_{\tilde{P}_{\text{PV},t}}^N} \left[-\max\left\{\sum_{\ell=1}^{N_\text{HVAC}} P_\ell u_\ell - \widetilde{P}_\text{PV}, 0 \right\} - \gamma \right]_+ \right\}$. } 
To linearize the the nonlinear constraints \eqref{eq:z-2}, we introduce big-M coefficients for $n=1,\ldots,N$, $M^1_n = \max_{u_t} \left\{ \left|\sum_{\ell=1}^{N_{\text{HVAC}}}P_\ell u_{t,\ell}  - P_{\text{total},t}^n \right| \right\} = \max\left\{\left|\sum_{\ell=1}^{N_{\text{HVAC}}}P_\ell  - P_{\text{total},t}^n\right|, P_{\text{total},t}^n\right\}$. \textcolor{black}{We also introduce an auxiliary variable $s_n = \max\left[\sum_{\ell=1}^{N_{\text{HVAC}}}P_\ell u_{t,\ell}  - P_{\text{total},t}^n, \ 0 \right]$ and a binary indicator variable $y_n$ for $n=1,\ldots, N$.} The constraints \eqref{eq:z-2} are equivalent to
% The constraints \eqref{eq:z1-1}
% -- \eqref{eq:z1-4} are equivalent to
% Therefore, the DR chance constraint \eqref{eq:drcc-cc} under the Wasserstein ambiguity set $\mathcal D_t^2$ is equivalent to  
the following MILP constraints.
% 
% 
% According to Corollary 2 in \citet{xie2018distributionally}, the DR chance constraint \eqref{eq:drcc-cc} under the Wasserstein ambiguity set $\mathcal D_t^2$ is equivalent to  the MILP constraints
% 
\begin{subequations}
	\label{eq:was}
	\begin{eqnarray}
		%  \label{eq:was-1}&& \delta_t  - \alpha_t \gamma \le \frac{1}{N}\sum_{n = 1}^N z_n\\
		\label{eq:was-2}&&  z_n + \gamma \le s_n, \ n = 1,\ldots, N\\
		\label{eq:was-3}&& s_n \le \sum_{\ell=1}^{N_{\text{HVAC}}}P_\ell u_{t,\ell}  - P_{\text{total},t}^n + M^1_n (1-y_n), \ n = 1,\ldots,N\\
		\label{eq:was-4}&& s_n \le M^1_n y_n, \ n = 1,\ldots, N\\
		% \label{eq:z-2}&& -\max\left[\sum_{j=1}^{N_{\text{HVAC}}}P_j u_{t,j}  - P_{\text{total},t}^\ell, \ 0 \right] \le -z_\ell - \gamma, \ \ell=1,\ldots,N \hspace{3mm}\\
		% \label{eq:was-5}&& z_n \le 0, \ n = 1,\ldots,N\\
		&& y_n \in \{0,1\}, \ n = 1,\ldots,N\\
		\label{eq:was-5}&& s_n \ge 0 , \ n=1,\ldots,N.
		% \label{eq:was-6}&&\gamma \ge 0 , 
	\end{eqnarray}
\end{subequations}
%where $ s_n, y_n$ are auxiliary variables. 
% The big-M coefficent $M^1_n = \max_{u_t} \left\{ \left|\sum_{j=1}^{N_{\text{HVAC}}}P_j u_{t,j}  - P_{\text{total},t}^n \right| \right\} = \max\left\{\left|\sum_{j=1}^{N_{\text{HVAC}}}P_j  - P_{\text{total},t}^n\right|, P_{\text{total},t}^n\right\}$.
Therefore, the DRCC formulation \eqref{eq:drcc} under $\mathcal D_t^2$ is reformulated as an MILP problem as follows.
% 
%\small
\begin{equation}\label{eq:lp-wasserstein}
	\text{\textbf{MILP1: }}
	\min_{u_t, \beta_{t,\ell}, x_{\ell,t}}\left\{ c_\text{sys}\sum_{\ell=1}^{N_{\text{HVAC}}}\beta_{t,\ell} + c_\text{switch}\sum_{\ell=1}^{N_{\text{HVAC}}}u_{t,\ell}:  \ \eqref{eq:x-function},\ \eqref{eq:absolute}-\eqref{eq:determ-binary},  \ \eqref{eq:z-1}, \ \eqref{eq:z-3}-\eqref{eq:z-4}, \ \eqref{eq:was-2}-\eqref{eq:was-5} \right\}.
	\nonumber
\end{equation}
\normalsize%
\begin{remark}Reformulation \eqref{eq:was}  involves 
	$3N$ constraints and $N$ binary variables, which may pose computational challenges as $N$ grows large. \end{remark}
Next, by exploiting the dual of the CVaR interpretation, we provide a more compact  MILP reformulation for the DR chance constraint with only   $2\lfloor \alpha_t N \rfloor + 2 $ constraints and $\lfloor \alpha_t N \rfloor + 1 $ binary variables. The problem size can be significantly reduced when $\alpha_t$ is small. %We note that the MILP reformulation can also be applied to any individual DR chance constraint with RHS uncertainty under Wasserstein ambiguity set $\mathcal D_t^2$.
%\textcolor{red}{ IS THERE ANY EXAMPLE?}
% 

We first sort $\{P^n_{\text{total},t}\}_{n=1}^N$ such that $P_{\text{total},t}^{(1)} \ge P_{\text{total}, t}^{(2)} \ge \cdots \ge P_{\text{total},t}^{(N)}$ and obtain the non-increasing permutation $\left\{(1), (2), \ldots, (n)\right\}$ of $\{1,2,\ldots,N \}$. Denote $P_{\text{total},t}^{(0)} = \sum_{\ell = 1}^{N_{\text{HVAC}}} P_\ell $, the maximum load provided by turning on all HVAC units. %Here, we assume that $P_{\text{total},t}^{(1)} > 0$. Otherwise, 
We make the following assumption.% for the more compact MILP reformulation.
\begin{assumption}\label{assump1}
	$ P_{\text{total},t}^{(0)}
	% \sum_{\ell = 1}^{N_{\text{HVAC}}} P_\ell 
	> P_{\text{total},t}^{(N)}$.
\end{assumption}
\noindent This is a mild assumption that requires  the smallest value of the solar PV generation realization  smaller than the maximum HVAC load provided when all the HVAC units are ON.% That is, the smallest PV generation can be absorbed fully when all the HVAC units are ON.
% \noindent (A1) there exists a solution $u_t$ such that $a_{(N)} > 0$.\\
%\textcolor{red}{REWRITE THIS THEOREM IN A MORE GENERAL FORM}
\textcolor{black}{\begin{theorem}\label{thm:Wasserstein-LP}
		Under Assumption \ref{assump1}, % For any solution $u_t$, assume that\\
		% Given a solution $u_t$ and $\alpha_t$, assume that\\
		% \noindent (A1) $a_{(N)} > 0$.\\
		%  \noindent (A2) there exist index $j^*$ such that $a_{(j^*)} > 0 \ge a_{(j^*-1)}$.\\
		% the DR chance constraint \eqref{eq:drcc-cc} under the Wasserstein ambiguity set $\mathcal D_t = \mathcal D_t^2$ is feasible if and only if the following MILP constraints are feasible
		% there exist feasible values of $\Delta_{jk}, \varepsilon_{jk}, o_{\ell j k}$ such that 
		the DR chance constraint \eqref{eq:drcc-cc} under the Wasserstein ambiguity set $\mathcal D_t^2$ is   feasible if and only if the following  linear constraints are feasible with auxiliary variables $a_n, h_n \in \{0,1\}, \ n =1,\ldots,k+1$, %\rho_{n\ell}\in\mathbb R, \ $,  
		%\ \ell = 1,\ldots,N_{\text{HVAC}}$,
		where \textcolor{black}{$k = \lfloor \alpha_t N \rfloor$}.
		\begin{subequations}\label{eq:drcc-wass-lp}
			\begin{eqnarray}
				&&\label{eq:drcc-wass-lp1}	\frac{1}{N}\sum_{n=1}^{{k}} {a}_{n} + ({\alpha}_t - \frac{k}{N}) {a}_{{k}+1} \ge \delta\\
				%	&& \frac{1}{N}\sum_{n=1}^{{k}} {a}_{(n)} + ({\alpha}_t - \frac{k}{N}) {a}_{({k}+1)} \ge \delta,
				&& \label{eq:drcc-wass-lp2} a_{n} \le \sum_{\ell=1}^{N_{\text{HVAC}}} P_\ell {u}_{t,\ell} - P_{\text{total},t}^{(n)} + M_n^2(1-h_n), \ n = 1,\ldots,k+1\\
				&& \label{eq:drcc-wass-lp3} a_{n} \le M_n h_n , \ n = 1,\ldots,k+1\\
				&&\label{eq:drcc-wass-lp4}  h_n\in\{0,1\}, \ a_n \ge 0 , \ n=1,\ldots,k+1,
			\end{eqnarray}
		\end{subequations}
		where $M_n$ and $M_n^2$ are  sufficiently large big-M constants.
		%${{k}}/{N} \le {\alpha}_t < {({k}+1)}/{N} $.
		%MILP constraints with auxiliary variables $\Delta_{jk}\in \{0,1\}, \varepsilon_{jk}\in\mathbb R, \tau_{\ell jk} \in \mathbb R, o_{\ell j k} \in \mathbb R $ are feasible. 
		% 
		%\begin{subequations}
		%\begin{eqnarray}
		%&& \label{eq:drcc-wass-lp1}
		% -\frac{1}{N} \sum_{n = 1}^{{k}} \sum_{i=1}^n x_i \left(P_{\text{total},t}^{({k}+1)} - P_{\text{total},t}^{(n)}\right) + \left({\alpha}_t - \frac{\sum_{n=1}^{k+1} n h_n-1}{N}\right) P_{\text{total},t}^{({k}+1)}- \left(\alpha_t + \frac{1}{N}\right) \sum_{\ell=1}^{N_{\text{HVAC}}}P_\ell u_{t,\ell}\nonumber\\
		% &&  + \frac{1}{N}\sum_{n=1}^{k+1}n\sum_{\ell=1}^{N_{\text{HVAC}}}P_\ell\rho_{n\ell}  \le - \delta\hspace{7mm}\\
		% && P_{\text{total},t}^{(n)} - P_{\text{total},t}^{(0)}(1-h_n)\le \sum_{\ell=1}^{N_{\text{HVAC}}}P_\ell {u}_{t,\ell} \le P_{\text{total},t}^{(n-1)} + P_{\text{total},t}^{(0)}(1-h_n), \ n = 1,\ldots, k+1\\
		% && \sum_{n=1}^{k+1} h_n = 1\\
		% \label{eq:drcc-wass-mcc1}&& \rho_{n\ell} \ge h_n + u_{t,\ell} - 1, \ n = 1,\ldots,k+1, \ \ell = 1,\ldots,N_{\text{HVAC}}\\
		% &&\rho_{n\ell} \le u_{t,\ell}, \ n = 1,\ldots,k+1, \ \ell = 1,\ldots,N_{\text{HVAC}}\\
		% && \rho_{n\ell} \le h_n, \ n = 1,\ldots,k+1, \ \ell = 1,\ldots,N_{\text{HVAC}}\\
		% \label{eq:drcc-wass-mcc4}&& \rho_{n\ell} \ge 0, \ n = 1,\ldots,k+1, \ \ell = 1,\ldots,N_{\text{HVAC}}\\
		% && \label{eq:drcc-wass-lp4} h_n \in \{0,1 \}, \ n = 1,\ldots, k+1.
		%\end{eqnarray}
		%\end{subequations}
		%
	\end{theorem}
}
\proof{Proof:}%
According to \eqref{eq:z}, the DR chance constraint \eqref{eq:drcc-cc} with $\mathcal D_t = \mathcal D_t^2$ is satisfied if and only if the optimal value of the following linear program is no more than $-\delta_t$.
\begin{subequations}\label{eq:z-opt}
	\begin{eqnarray}\label{eq:z-opt1}
		\min && -\frac{1}{N} \sum_{n=1}^{N} z_n - {\alpha}_t \gamma \\ % + \delta_t \le 0\\
		\label{eq:z-opt2}\mbox{s.t.} && -{a}_n \le - z_n - \gamma, \ n = 1,\ldots,N\\
		&& z_n \le 0, \ n = 1,\ldots,N\\
		\label{eq:z-opt3}&& \gamma \ge 0,
	\end{eqnarray}
\end{subequations}
where 
\begin{equation}\label{eq:a-def}
	{a}_n = \max \left[ \sum_{\ell=1}^{N_{\text{HVAC}}} P_{\ell} {u}_{t,\ell} - P_{\text{total},t}^{(n)}, 0 \right], \ n = 1,\ldots,N.
\end{equation} 
We associate dual variables $\pi_n \ge 0,\ n=1,\ldots,N$ with constraints in \eqref{eq:z-opt2} and obtain the dual problem as follows.
\begin{subequations}
	\label{eq:z-opt-dual}
	\begin{eqnarray}
		\label{eq:z-opt-dual-obj}\max && -\sum_{n=1}^N \pi_n {a}_n \\ \label{eq:z-opt-dual-1}\mbox{s.t.} && \pi_n \le \frac{1}{N}, \ n = 1,\ldots,  N\\
		\label{eq:z-opt-dual-2}&& \sum_{n=1}^N \pi_n \ge {\alpha}_t\\
		\label{eq:z-opt-dual-3}&& \pi_n \ge 0, \ n= 1,\ldots,N.
	\end{eqnarray}
\end{subequations}
Due to strong duality, the optimal value of \eqref{eq:z-opt-dual} equals to 
that of \eqref{eq:z-opt}, which is $\le -\delta_t$. \textcolor{black}{Note that \eqref{eq:z-opt-dual} is the dual interpretation of the CVaR \eqref{eq:cvar} scaled by $\alpha_t$.}
% The problem \eqref{eq:z-opt-dual-obj}--\eqref{eq:z-opt-dual-3} is a resource allocation problem, which  allocates $\alpha_t$ to $\pi$ such that $-\sum_{n=1}^{N} \pi_n a_n$ is maximized. 
The problem \eqref{eq:z-opt-dual}  is always feasible as $0<{\alpha_t}< 1$. %Based on Assumption (A1), we have $\hat{a}_{(N)} >0$. 
If Assumption \ref{assump1} does not hold, the DR chance constraint \eqref{eq:drcc-cc} is infeasible. Because when the Assumption \ref{assump1} does not hold,
${a}_{n} = 0,\ n=1,\ldots,N$, for any solution of ${u}_t$ and thus the optimal value of the dual problem \eqref{eq:z-opt-dual}
is zero. As the optimal value is more  than $-\delta_t$, the DR chance constraint \eqref{eq:drcc-cc} is then infeasible and so is the DRCC problem \eqref{eq:drcc}.
% If the problem \eqref{eq:z-opt-dual-obj}--\eqref{eq:z-opt-dual-3} 
% is feasible,

% We assume that $k/N \le \alpha_t < (k+1)/N$. 
The problem \eqref{eq:z-opt-dual} can be converted to a relaxed knapsack problem and one optimal solution that can be obtained by the greedy algorithm    is
\begin{equation}\nonumber\label{eq:dual-sol}
	\pi_{n} = \begin{cases} \frac{1}{N}, & n = 1,\ldots,{k}\\ {\alpha}_t - \frac{{k}}{N}, & n={k}+1\\ 0,& n={k}+2,\ldots N.\end{cases}
	%  \pi_n = \frac{1}{N}, \ n = 1,\ldots,k, \ \pi_{k+1} = \alpha_t - \frac{k}{N}, \ \pi_{n} = 0, \ n = k+2,\ldots,n.
\end{equation}
% where $\hat{k}$ is an index such that 
% \begin{equation}\label{eq:alpha-con}\frac{\hat{k}}{N} \le {\alpha}_t < {(\hat{k}+1)}{N} .\end{equation}  %$a_{(k-1)} \le 0 < a_{(k)}$. 
\textcolor{black}{We note that $\pi/\alpha_t$ is the conditional (discrete) probability distribution of CVaR in \eqref{eq:cvar}.} Recall that $\left\{ (1), (2), \ldots, (n)\right\}$ is a permutation such that $P_{\text{total},t}^{(1)} \ge P_{\text{total}, t}^{(2)} \ge \cdots \ge P_{\text{total},t}^{(N)}$, or equivalently, $0 \le {a}_{1} \le {a}_{2} \le \cdots \le {a}_{N}$ (because for any $n \le m$, $\sum_{\ell=1}^{N_{\text{HVAC}}} P_\ell {u}_{t,\ell} - P_{\text{total},t}^{(n)} \le \sum_{\ell=1}^{N_{\text{HVAC}}} P_\ell {u}_{t,\ell} - P_{\text{total},t}^{(m)} $ and thus
${a}_{n} \le {a}_{m}$) and $k = \lfloor \alpha_t N \rfloor$ is an index such that  ${{k}}/{N} \le {\alpha}_t < {({k}+1)}/{N}$.
%\begin{equation}\label{eq:k} 
%{{k}}/{N} \le {\alpha}_t < {({k}+1)}/{N}
%\end{equation}%$a_{(k-1)} \le 0 < a_{(k)}$.
{In the case where  ${k}=0$, the optimal solution is $\pi_{1} = {\alpha}_t$ and $\pi_{n} = 0, \ n= 2,\ldots,N$.}

%To see the optimality of the solution \eqref{eq:dual-sol}, we notice that the dual problem \eqref{eq:z-opt-dual} can be reformulated as a relaxed knapsack problem by replacing the decision variable $\pi_n$ with its negative  $\pi_n^\prime = -\pi_n$.
%It is well known that the relaxed knapsack problem can be effectively solved by the greedy algorithm by assigning the variable with the largest objective coefficient :
%\begin{enumerate}
%    \item  Order the decision variables $\pi^\prime_n, \ n=1,\ldots,N$ in a non-decreasing order of their objective coefficients ${a}_{n}$  (Note that the non-decreasing order is provided by  the permutation  $\left\{ (1), (2), \ldots, (n)\right\}$, i.e., 
%    % $P_{\text{total},t}^{(1)} \ge P_{\text{total}, t}^{(2)} \ge \cdots \ge P_{\text{total},t}^{(N)}$. We have
%    % 
%    $0 \le a_{(1)} \le {a}_{(2)} \le \cdots \le {a}_{(N)}$.
%);
%\item Find the  index ${k}$ based on the condition \eqref{eq:k};
%\item Let $\pi_{(n)}^\prime = -1/N, \ n = 1,\ldots,{k}$, $\pi_{(k+1)}^\prime = {\alpha}_t - {k}/N$, and $\pi_{(n)}^\prime = 0, \ n = {k}+2,\ldots,N$.
%\end{enumerate}
%% 
%

\textcolor{black}{Now, we obtain an equivalent reformulation of \eqref{eq:z} as %$	-1/N\sum_{n=1}^{{k}} {a}_{n} - ({\alpha}_t - k/N) {a}_{{k}+1} \le -\delta \text{ and constraint \eqref{eq:a-def}}.$
	\begin{equation}
		-\frac{1}{N}\sum_{n=1}^{{k}} {a}_{(n)} - ({\alpha}_t - \frac{k}{N}) {a}_{({k}+1)} \le -\delta_t \text{ and constraint \eqref{eq:a-def}}.
	\end{equation}
	We obtain\eqref{eq:drcc-wass-lp1} by multiplying -1 on both sides of the first inequality above. }
%
% The solution obtained by the greedy algorithm is optimal to the relaxed knapsack problem. So is \eqref{eq:dual-sol} to the problem \eqref{eq:z-opt-dual-obj}--\eqref{eq:z-opt-dual-3} with optimal value $-1/N\sum_{n=1}^{\hat{k}} \hat{a}_{(n)} - (\hat{\alpha}_t - \hat{k}/N) \hat{a}_{(\hat{k}+1)} $.
% 
% 
% 
% Otherwise, the optimal value of the dual problem \eqref{eq:z-opt-dual-obj}--\eqref{eq:z-opt-dual-3} 
% is zero, which is greater than $-\delta$ and therefore the DRCC under the Wasserstein set $\mathcal D_t^2$ is infeasible. 
\textcolor{black}{To linearize \eqref{eq:a-def}, we introduce binary variables $h_n$ for $n=1,\ldots,k+1$ and big-M coefficients to obtain \eqref{eq:drcc-wass-lp2} and \eqref{eq:drcc-wass-lp3}.} Now, we complete the proof.
\Halmos
\endproof
The big-M coefficients in Theorem \ref{thm:Wasserstein-LP} can take values as $M_n^2 = P_{\text{total},t}^{(n)} - P_{\text{total},t}^{(k)}$ and $M_n = M_n^1$.
Now, we obtain the second  MILP reformulation of the DRCC formulation \eqref{eq:drcc} under the Wasserstein ambiguity set $\mathcal D_t^2$.
\begin{equation}\label{eq:lp-wasserstein2}
	\text{\textbf{MILP2: }}\min_{u_t, \beta_{t,\ell}, x_{t,\ell}}\left\{ c_\text{sys}\sum_{\ell=1}^{N_{\text{HVAC}}}\beta_{t,\ell} + c_\text{switch}\sum_{\ell=1}^{N_{\text{HVAC}}}u_{t,\ell}:  \ \eqref{eq:x-function},\ \eqref{eq:absolute}-\eqref{eq:determ-binary},  \ \eqref{eq:drcc-wass-lp1} - \eqref{eq:drcc-wass-lp4} \right\}.
	\nonumber
\end{equation}
MILP2 can be further strengthened by the following proposition.
%To strengthen MILP2, we present the following proposition based on the underlying problem structure.
\textcolor{black}{\begin{proposition} \label{prop:strengthen}
		Let $k = \lfloor \alpha_t N\rfloor$.
		\begin{enumerate} 
			\item[i.]If MILP2 is feasible,  then
			\begin{equation}\label{eq:k+1}a_{k+1} =\sum_{\ell=1}^{N_{\text{HVAC}}} P_\ell {u}_{t,\ell} - P_{\text{total},t}^{({k}+1)}.\end{equation}
			Thus constraints \eqref{eq:drcc-wass-lp2} and \eqref{eq:drcc-wass-lp3} indexed with $k+1$ and variable $h_{k+1}$ can be removed from MILP2.
			\item[ii.] The following inequalities are valid for MILP2.
			\begin{equation}\label{eq:ordering-valid}
				h_n \le h_{n+1}, \ n = 1,\ldots,k-1.
			\end{equation}
		\end{enumerate}
\end{proposition}}
\proof{Proof:}%
For a given solution  ${u}_t$, \textcolor{black}{there exists a critical index $j$ such that}
%we are interested in a critical index ${j}$ such that 
\begin{equation*}P_{\text{total},t}^{({j})} < \sum_{\ell=1}^{N_{\text{HVAC}}}P_\ell {u}_{t,\ell} \le  P_{\text{total},t}^{({j}-1)}\end{equation*} 
and thus $a_{1} = \ldots ={a}_{{j} -1 } = 0 < {a}_{{j}} \le \ldots\le a_{N} $. In the case where for all $n=1,\ldots,N$, ${a}_{n} > 0$, or equivalently $P_{\text{total},t}^{n} < \sum_{\ell = 1}^{N_{\text{HVAC}}} P_{\ell} {u}_{t,\ell}$, let $j=1$. \textcolor{black}{According to constraints \eqref{eq:drcc-wass-lp2}--\eqref{eq:drcc-wass-lp3}, binary variable $h_n = 1$ if only if $a_n > 0 $. Given a feasible solution $u_t$, $h_n = 0$ for $n=1,\ldots,j-1$ and $h_n = 1$ for $n=j,\ldots,k$. So inequalities \eqref{eq:ordering-valid} are valid.}

If MILP2 is feasible, then ${j} \le k+1$. Since, otherwise, the optimal value of \eqref{eq:z-opt-dual} is zero which is larger than $-\delta$ and result in infeasibility.   \textcolor{black}{ Therefore, \eqref{eq:k+1} holds. 
}
\Halmos
\endproof
\begin{remark}
	\textcolor{black}{ Theorem \ref{thm:Wasserstein-LP} and Proposition \ref{prop:strengthen} also apply to general DR individual chance constraint with RHS uncertainty under the Wasserstein ambiguity. In particular, Theorem \ref{thm:Wasserstein-LP} can be  extended to a joint chance constraint setting. For the sake of space, we leave this to future work.  
		% 
		% \citet{ho2020distributionally} also exploit the CVaR interpretation and the quantile strengthening to derive an MILP reformulation for joint DR chance constraint with RHS uncertainty over Wasserstein ambiguity set. 
		%\citet{ho2020distributionally} study DRCC problems with RHS uncertainty based on the CVaR interpretation as well. 
		Different from \citet{ho2020distributionally} based on the CVaR interpretation, 
		we exploit the dual perspective of the CVaR interpretation which results in fewer binary variables (\citet{ho2020distributionally} require $N$ binary variables in contrast with $\lfloor \alpha_t N\rfloor$ in MILP2). We also reveal the solution structure of the binary variables associated with samples in the individual chance constraint setting.  
		In Section \ref{sec:W-CPU}, we show the computational comparison between their reformulation and MILP2. }
\end{remark}

\section{Adjustable  Chance-Constrained Formulations}
\label{sec:adjustable-cc}

In the DRCC formulation \eqref{eq:drcc}, $\alpha_t$ is a pre-defined risk level, which is often chosen based on operators' experience and is usually a small number to guarantee high PV generation utilization. However, a too small $\alpha_t$ may result in infeasibility of the problem and can lead to high operational cost  \citep[see, e.g.,][]{ma2019distributionally}. It is challenging for the system operator to decide the value of $\alpha_t$  for optimally trading off between the risk and cost. Next, 
%following \cite{ma2019distributionally,qiu2016data,wang2018adjustable}, 
we consider $\alpha_t$ as an adjustable decision variable and modify the DRCC model \eqref{eq:drcc} into the following adjustable DRCC problem. 
\begin{subequations}
	\label{eq:adjustable-drcc}
	\begin{eqnarray}
		\label{eq:adjustable-drcc-obj}
		\min_{u_t,\beta_{t,\ell}, x_{t,\ell}, \alpha_t} && c_\text{sys}\sum_{\ell=1}^{N_{\text{HVAC}}}\beta_{t,\ell} + c_\text{switch}\sum_{\ell=1}^{N_{\text{HVAC}}}u_{t,\ell} + c_t \alpha_t \\
		%\left( x - x_{\text{ref}}\right)^{\top}Q\left(x  - x_{\text{ref}}\right)\\
		%\min_{u,z}&& z \\
		\label{eq:adjustable-drcc-cc}\mbox{s.t.} %&& z \ge \left( x - x_{\text{ref}}\right)^{\top}Q\left(x  - x_{\text{ref}}\right) \\
		&& \inf_{f \in\mathcal \mathcal D_t} \mathbb P\left( \sum_{\ell=1}^{N_{\text{HVAC}}} P_\ell u_{t,\ell} - \sum_{i=1}^{N_\text{PV}}\widetilde{P}_{\text{PV},t,i} \ge 0\right)  \ge 1-\alpha_t \hspace{7mm}\\
		% && x_\text{min} \le x_{t,j} \le x_\text{max},\ j = 1,\ldots,N_{\text{HVAC}}\\
		\label{eq:adjustable-alpha-bounds}&& 0 \le \alpha_t \le 1\\
		%  && u_t \in \{0,1\}^{N_{\text{HVAC}}}, \\
		&& \eqref{eq:x-function}, \ \eqref{eq:absolute}-\eqref{eq:determ-binary},\nonumber
		% && \eqref{eq:determ-x-bounds}-\eqref{eq:determ-binary}. \nonumber
		%&& x = \frac{1}{RC} T_\text{out} - \frac{1}{RC} x_0
		%+ \frac{1}{C} Q_\text{out} + \frac{u}{C} Q_\text{HVAC}\\
		%  && x_\text{min} \le x_k \le x_\text{max},\ k=1,\ldots,N_{\text{HVAC}}\\
		%  && u \in \{0,1\}^{N_{\text{HVAC}}}
	\end{eqnarray}
\end{subequations}
where $c_t > 0$ is the coefficient weight on the risk level $\alpha_t$. If $c_t$ is close to zero, the utilization of solar PV output is low. Section \ref{sec:wass-adjustable-alpha} presents the impact of $c_t$ on the  risk level  and operational cost.
% \paragraph{Moment-based Ambiguity Set $\mathcal D_t = \mathcal D_t^1$}

% % We present constraint \eqref{eq:linear-eq} in a more general form
% % \begin{equation}

% \end{equation}
\subsection{Adjustable DRCC Reformulations under Moment-based Ambiguity Set $\mathcal D_t^1$}
According to \eqref{eq:linear-eq}, the choice of the coefficient $\Omega_t$ depends on the values of $\alpha_t$ and $\gamma_1/\gamma_2$. In this section, we present 0-1 SOCP reformulations under the two cases: (i) $\gamma_1/\gamma_2 \le \alpha_t$ and (ii) $\gamma_1/\gamma_2 > \alpha_t$. 
%$\Omega_t = $ $\Omega_t = \left\{ \begin{array}{ll}
% \sqrt{\gamma_1} + \sqrt{{1-\alpha_t}/{\alpha_t} (\gamma_2 - \gamma_1)}, & {\gamma_1}/{\gamma_2} \le \alpha_t\\
%\sqrt{{\gamma_2}/{\alpha_t}}, & {\gamma_1}/{\gamma_2} > \alpha_t
%\end{array}\right.$. 
% For simplicity, we drop the time index of \eqref{eq:linear-eq}:
% \begin{equation}\label{eq:linear-eq-1}
%   \sum_{j=1}^{N_{\text{HVAC}}}P_ju_{j} \ge \theta + \Omega \sigma,
% \end{equation}
% where $\Omega = \left\{ \begin{array}{ll}
%  \sqrt{\gamma_1} + \sqrt{\left(\frac{1-\alpha}{\alpha} (\gamma_2 - \gamma_1)\right)}, & \frac{\gamma_1}{\gamma_2} \le \alpha\\
% \sqrt{\frac{\gamma_2}{\alpha}}, & \frac{\gamma_1}{\gamma_2} > \alpha
% \end{array}\right.$. Denote $u = (u_j,\ j = 1,\ldots,N_{\text{HVAC}})^\top$.
\begin{theorem}\label{thm:socp1}
	If ${\gamma_1}/{\gamma_2} \le \alpha_t$,  the adjustable DR chance constraint \eqref{eq:adjustable-drcc-cc} under $\mathcal D_t^1$, or equivalently, 
	%\begin{subequations}
	\begin{equation}\label{eq:linear-ref1}
		\sum_{\ell=1}^{N_{\text{HVAC}}}P_\ell u_{t,\ell} \ge \theta_t + \left(\sqrt{\gamma_1} + \sqrt{\frac{1-\alpha_t}{\alpha_t} (\gamma_2 - \gamma_1)}\right) \sigma_t,\end{equation}
	%\end{subequations} 
	is equivalent to the following 0-1 SOCP constraints
	\begin{subequations}
		\begin{eqnarray}
			\label{eq:thm-soc1}&& \left\| \begin{matrix} 2\sigma_t \sqrt{\gamma_2 - \gamma_1} \\ \alpha_t - d \end{matrix} \right\| \le \alpha_t + d\\
			\label{eq:thm-soc2}&& d \le P\cdot g - 2 \left(\theta_t + \sigma_t \sqrt{\gamma_1}\right)\sum_{\ell =1}^{N_{\text{HVAC}}} P_\ell u_{t,\ell}  + \theta_t^2 + 2 \theta_t\sigma_t \sqrt{\gamma_1} + \gamma_2\sigma_t^2 \\
			%(\theta_t + \sigma_t \sqrt{\gamma_1})^2 + \sigma_t^2 (\gamma_2 - \gamma_1)\\
			\label{eq:thm-soc7}&& \sum_{\ell = 1}^{N_{\text{HVAC}}} P_\ell u_{t,\ell} \ge \theta_t + \sigma_t \sqrt{\gamma_1}\\
			\label{eq:thm-soc3}&& g_{ij} \ge u_{t,i} + u_{t,j} - 1, \ g_{ij} \le u_{t,i}, \ g_{ij} \le u_{t,j},\ g_{ij} \ge 0, \ i,j = 1,\ldots, N_{\text{HVAC}},
			%
			%\label{eq:thm-soc3}&& g_{ij} \ge u_{t,i} + u_{t,j} - 1, \ i,j = 1,\ldots, N_{\text{HVAC}} \\ 
			%\label{eq:thm-soc4}&& g_{ij} \le u_{t,i}, \ i,j = 1,\ldots, N_{\text{HVAC}} \\
			%\label{eq:thm-soc5}&& g_{ij} \le u_{t,j}, \ i,j = 1,\ldots, N_{\text{HVAC}} \\
			%\label{eq:thm-soc6}&& g_{ij} \ge 0, \ i,j = 1,\ldots, N_{\text{HVAC}},
		\end{eqnarray}
	\end{subequations}
	where the operator $\cdot$  represents the Frobenius inner product, 
	$P\in \mathbb R^{N_{\text{HVAC}} \times N_{\text{HVAC}}}$ with $P_{ij} = P_i P_j$ and $g \in \mathbb R^{N_{\text{HVAC}} \times N_{\text{HVAC}}}$.
\end{theorem}
\proof{Proof:}
See Appendix \ref{sec:proof-thm-socp1} in the online supplement.
\Halmos\endproof

\begin{theorem}\label{thm:socp2}
	If ${\gamma_1}/{\gamma_2} > \alpha_t$,  the adjustable DR chance constraint \eqref{eq:adjustable-drcc-cc}, or equivalently,
	\begin{equation}\label{eq:linear-ref2}
		\sum_{\ell=1}^{N_{\text{HVAC}}}P_\ell u_{t,\ell} \ge \theta_t + \sqrt{\frac{\gamma_2}{\alpha_t}} \sigma_t,\end{equation}
	is equivalent to the following 0-1 SOCP constraints
	\begin{subequations}
		\begin{eqnarray}
			\label{eqn:soc-1} && \sum_{\ell=1}^{N_{\text{HVAC}}}P_\ell u_{t,\ell} \ge \theta_t + \phi \sigma_t\sqrt{\gamma_2}\\
			\label{eqn:soc-2}&& \alpha_t + \phi \ge \left\|\begin{matrix}\alpha_t-\phi\\ 2q\end{matrix}\right\|\\
			\label{eqn:soc-3}&& \phi \ge w^2\\
			\label{eqn:soc-4}&& q+w \ge \left\|\begin{matrix} q-w\\ 2
			\end{matrix}\right\|.
			% \label{eqn:ref-4+} && w \ge 0
		\end{eqnarray}
	\end{subequations}
\end{theorem}
\proof{Proof:}
See Appendix \ref{sec:proof-thm-socp2}  in the online supplement.
\Halmos\endproof

Therefore, to solve the adjustable DRCC model \eqref{eq:adjustable-drcc} under the moment-based ambiguity set $\mathcal D_t^1$, we solve the following two 0-1 SOCP problems, separately.
% \small
\begin{equation} \text{\textbf{SOCP1:}}\ \label{eq:socp-moment1}
	\min_{u_t,\beta_{t,\ell}, x_{t,\ell}, \alpha_t}\left\{ c_\text{sys}\sum_{\ell=1}^{N_{\text{HVAC}}}\beta_{t,\ell} + c_\text{switch}\sum_{\ell=1}^{N_{\text{HVAC}}}u_{t,\ell} + c_t \alpha_t:  \ 1 \ge \alpha_t \ge \frac{\gamma_1}{\gamma_2},\ \eqref{eq:x-function}, \ \eqref{eq:absolute}-\eqref{eq:determ-binary}, \ \eqref{eq:thm-soc1}-\eqref{eq:thm-soc3} \right\},
	\nonumber
\end{equation}
\begin{equation} \text{\textbf{SOCP2:}}\ \label{eq:socp-moment2}
	\min_{u_t,\beta_{t,\ell}, x_{t,\ell}, \alpha_t}\left\{ c_\text{sys}\sum_{\ell=1}^{N_{\text{HVAC}}}\beta_{t,\ell} + c_\text{switch}\sum_{\ell=1}^{N_{\text{HVAC}}}u_{t,\ell} + c_t \alpha_t:  \ 0\le \alpha_t \le \frac{\gamma_1}{\gamma_2},\ \eqref{eq:x-function} , \ \eqref{eq:absolute}-\eqref{eq:determ-binary}, \  \eqref{eqn:soc-1}-\eqref{eqn:soc-4} \right\}.
	\nonumber
\end{equation}
\normalsize
After obtaining the two optimal values,  we compare them and let  the solution with the higher optimal value be the optimal solution to the adjustable DRCC model \eqref{eq:adjustable-drcc}.
% 
%  
%  Therefore, the adjustable DRCC model \eqref{eq:adjustable-drcc} under the moment-based ambiguity set is equivalent to the 0-1 second-order conic program

%   \begin{equation}\label{eq:socp-moment}
% \min_{u_t}\left\{ c_\text{sys}\sum_{j=1}^{N_{\text{HVAC}}}|x_{t,j}-x_{\text{ref}}| + c_\text{switch}\sum_{j=1}^{N_{\text{HVAC}}}u_{t,j} + c_t \alpha_t:  \ \eqref{eqn:soc-1}-\eqref{eqn:soc-4}, \ \eqref{eq:determ-x-bounds}-\eqref{eq:determ-binary} \right\}.
% \nonumber
% \end{equation}
% 
% 
% 
Note that the reformulation SOCP1
% the reformulation \eqref{eq:thm-soc1}--\eqref{eq:thm-soc7} 
incorporates $N_{\text{HVAC}}^2$ auxiliary variables $g_{ij}, \ i,j = 1,\ldots,N_{\text{HVAC}}$, which can result in computational burden when $N_{\text{HVAC}}$ is large. Below, 
%inspired by 
%the proof of Theorem 14 in 
%\citet{xie2019optimized}, 
we present a more compact approximation that incorporates only four auxiliary variables if ${\gamma_1}/{\gamma_2} \le \alpha_t \le 0.75$.
\begin{theorem}\label{thm:moment-adjust}
	if $\gamma_1/\gamma_2 \le \alpha_t \le 0.75$, the adjustable DR chance constraint \eqref{eq:adjustable-drcc-cc} is outer approximated by the following 0-1 SOCP constraints
	\begin{subequations}
		\begin{eqnarray}
			\label{eq:thm-outer-qc1}&& \sum_{\ell=1}^{N_{\text{HVAC}}} P_\ell u_{t,\ell} \ge \theta_t + \left( \sqrt{\gamma_1} + r\sqrt{\gamma_2 - \gamma_1} \right)\sigma_t\\
			%		 \sigma_t\sqrt{\gamma_1} \ge \sigma_t \sqrt{\gamma_2 - \gamma_1} r\\ 
			\label{eq:thm-outer-socp1}
			&& 2 r \ge \phi\\
			&& \alpha_t + \phi \ge \left\| \begin{matrix} \alpha_t -\phi \\ 2q \end{matrix}
			\right\|\\
			&& q + w \ge \left\| \begin{matrix} q - w\\ 2
			\end{matrix}
			\right\|\\
			\label{eq:thm-outer-socp2}&& \phi \ge w^2.
		\end{eqnarray}
	\end{subequations}
\end{theorem}
\proof{Proof:}
See Appendix \ref{sec:proof-thm-moment-adjust}  in the online supplement.
\Halmos\endproof
Therefore,
when $\gamma_1/\gamma_2 \le \alpha_t \le 0.75$, to solve the adjustable DRCC model \eqref{eq:adjustable-drcc} under the moment-based ambiguity set, we can implement a branch-and-cut algorithm, which solves the following 0-1 SOCP problem.
\small
%\noindent\textbf{SOCP3:}
\begin{equation}
	\text{\textbf{SOCP3:}}\ 
	\label{eq:socp-moment}
	\min_{u_t,\beta_{t,\ell}, x_{t,\ell}, \alpha_t}\left\{ c_\text{sys}\sum_{\ell=1}^{N_{\text{HVAC}}}\beta_{t,\ell} + c_\text{switch}\sum_{\ell=1}^{N_{\text{HVAC}}}u_{t,\ell} + c_t \alpha_t:  \ 0.75 \ge \alpha_t \ge \frac{\gamma_1}{\gamma_2},\ \eqref{eq:x-function}, \ \eqref{eq:absolute}-\eqref{eq:determ-binary},  \
	\eqref{eq:thm-outer-qc1}-\eqref{eq:thm-outer-socp2}
	\right\}.
	\nonumber
\end{equation}
\normalsize
At each iteration, obtaining the current solution $(\hat{\alpha}_t, \hat{r}, \hat{u}_t)$, if \textcolor{black}{constraint \eqref{eq:thm-qc2} in the online appendix \ref{sec:proof-thm-moment-adjust}: $r \ge \sqrt{({{1-\alpha_t})/{\alpha_t}}}$,} is  satisfied, we claim that $(\hat{\alpha}_t, \hat{r}, \hat{u}_t)$  is optimal. Otherwise, we generate the following supporting hyperplane as a valid inequality.
\begin{equation*}
	r \ge \left( -\frac{1}{2} (1-\hat{\alpha}_t)^{-\frac{1}{2}} \hat{\alpha}_t^{-\frac{3}{2}} \right) \alpha_t + (1-\hat{\alpha}_t)^{-\frac{1}{2}} \hat{\alpha}_t^{-\frac{1}{2}} (\frac{3}{2} - \hat{\alpha}_t).
\end{equation*}    
If the decision maker chooses to only consider $\alpha_t$ such that $0 \le \alpha_t \le 0.75$ for the adjustable DRCC model, he first solves SOCP2 for $0\le \alpha_t \le \gamma_1/\gamma_2$ and implement the branch-and-cut algorithm for $\gamma_1/\gamma_2 \le \alpha_t \le 0.75$. Then comparing the two optimal values, let the solution with the higher optimal value  be the optimal solution to the adjustable DRCC model.

% \paragraph{Wasserstein Ambiguity Set $\mathcal D_t = \mathcal D_t^2$}
\subsection{Adjustable DRCC Reformulations under Wasserstein Ambiguity Set $\mathcal D_t^2$}\label{sec:adjustable-wasserstein}
% The formulation \eqref{eq:was} provides a 0-1 SOCP reformulation for the adjustable DRCC model under the Wasserstein ambiguity set. In the following we show that the DRCC model can be reformulated as a 0-1 MILP model.
In the section, we first show that the adjustable DRCC model under the Wasserstein ambiguity set $\mathcal D_t^2$ can be reformulated as a 0-1 MILP formulation with big-M coefficients based on the MILP1. Then we present a big-M free MILP reformulation by exploiting the RHS uncertainty. %But we note that the result apply to general individal chance constraint with the RHS uncertainty.

% According to Corollary 2 in \citet{xie2018distributionally}, the feasible region described by the DR chance constraint \eqref{eq:drcc-cc} is 
% \begin{subequations}
% \label{eq:z}
% \begin{eqnarray}
% Z =\bigg\{ u_{t}: \label{eq:z-1}&& \delta_t  - \alpha_t \gamma \le \frac{1}{N}\sum_{n = 1}^N z_n\\
% \label{eq:z-2}&& -\max\left[\sum_{j=1}^{N_{\text{HVAC}}}P_j u_{t,j}  - P_{\text{total},t}^n, \ 0 \right] \le -z_n - \gamma, \ n=1,\ldots,N \hspace{3mm}\\
% \label{eq: z-3}&& z_n \le 0, \ n = 1,\ldots,N\\
% \label{eq:z-4}&&\gamma \ge 0 \bigg\}.
% \end{eqnarray}
% \end{subequations}

We denote $Z = \{u_t: \ \eqref{eq:z-1}-\eqref{eq:z-4} \}$  the feasible region (see Section \ref{sec:drcc-ref-wasserstein}) described by the DR chance constraint \eqref{eq:drcc-cc} under Wasserstein ambiguity set $\mathcal D_t^2$. 
% In Section  \ref{sec:drcc-ref-wasserstein}, according to Corollary 2 in \citet{xie2018distributionally}, we denote $Z = \{u_t: \ \eqref{eq:z-1}-\eqref{eq:z-4} \}$ which is the feasible region described by the DR chance constraint \eqref{eq:drcc-cc} under Wasserstein ambiguity set $\mathcal D_t^2$. 
Now, we consider the risk level $\alpha_t$ as a decision variable in the adjustable formulation. The product $\alpha_t\gamma$ in \eqref{eq:z-1} becomes a bilinear term. In the following proposition, we derive an equivalent set to $Z$ by eliminating the bilinear term.

% But here the risk level $\alpha_t$ becomes a decision variable in the adjustable chance constraint \eqref{eq:adjustable-drcc-cc} instead of a known parameter in \eqref{eq:drcc-cc}. In the following proposition, we eliminate the bilinear term $\alpha_t\gamma$ in \eqref{eq:z-1}.

\begin{proposition}\label{prop:adjust}
	The set $Z$ is equivalent to the following set:
	\begin{subequations}
		\label{eq:z1}
		\begin{eqnarray}
			Z_1 = \bigg\{ u_{t}: \label{eq:z1-1}&& \delta_t \lambda - \alpha_t \le \frac{1}{N}\sum_{n=1}^N z_n\\
			\label{eq:z1-2}&& -\max\left[ \sum_{\ell=1}^{N_{\text{HVAC}}}P_\ell \lambda u_{t,\ell} - P_{\text{total},t}^n \lambda, \ 0 \right] \le -z_n - 1, \ n = 1,\ldots,N\hspace{7mm}\\
			\label{eq:z1-3}&& \lambda \ge 0,\ z_n \le 0, \ n = 1,\ldots,N
			%\label{eq:z1-4}&& 
			% && v_{t,j} \ge 0, \ j = 1,\ldots, N_{\text{HVAC}}\\
			% && v_{t,j} \ge  \lambda^U u_{t,j} + \lambda - \lambda^U, \ j = 1,\ldots, N_{\text{HVAC}}\\
			% && v_{t,j}\le \lambda,\ j = 1,\ldots, N_{\text{HVAC}} \\
			% && v_{t,j} \le \lambda^U u_{t,j}, \ j = 1,\ldots, N_{\text{HVAC}} 
			\bigg\}. 
		\end{eqnarray}
	\end{subequations}
\end{proposition}

\proof{Proof:}
See appendix \ref{sec:proof-prop-adjust}  in the online supplement.
\Halmos\endproof
By introducing suitable big-M coefficients $M_n^3, \ n = 1,\ldots,N$, the set $Z_1$ can be further reformulated as a mixed integer set below.
\begin{subequations}\label{eq:z1-mip}
	\begin{eqnarray}
		Z_1 = \bigg\{ u_{t}:\label{eq:z1-mip1}&& \delta_t\lambda - \alpha_t \le \frac{1}{N}\sum_{n = 1}^N z_n\\
		\label{eq:z1-mip2}&& z_n + 1 \le s_n,\ n = 1,\ldots,N\\
		\label{eq:z1-mip3}&& s_n \le \sum_{\ell=1}^{N_{\text{HVAC}}} P_\ell \lambda u_{t,\ell} - P_{total,t}^n \lambda  + M^3_n (1-y_n), \ n = 1,\ldots, N\\
		\label{eq:z1-mip4}&& s_n \le M^3_n y_n, \ n = 1,\ldots,N\\
		%\label{eq:z1-mip5}&& \lambda \ge 0\\
		\label{eq:z1-mip6}&& \lambda \ge 0, \ s_n \ge 0,\ z_n \le 0, \ y_n\in\{0,1\}, \ n = 1,\ldots,N\bigg\}.
	\end{eqnarray}
\end{subequations}
%where $M_n^3$ is a suitable big-M constant for each $n = 1,\ldots,N$. 
We denote $w_{t,\ell} := \lambda u_{t,\ell}$ and constraint 
\eqref{eq:z1-mip3} becomes
\begin{equation}\label{eq:z1-mip3-ref}
	s_n \le \sum_{\ell=1}^{N_{\text{HVAC}}} P_\ell w_{t,\ell} - P_{total,t}^n \lambda  + M^3_n (1-y_n), \ n = 1,\ldots, N.
\end{equation}
According to the McCormick inequalities, we introduce the following linear constraints 
% The nonlinear term $\lambda u_{t,j}$ in \eqref{eq:z1-mip3} can be replace by $w_{t,j}$ which satisfies McCormick inequalities as
% There remain  nonlinear terms $\lambda u_{t,j}$ in \eqref{eq:z1-mip3}. Let $w_{t,j} = \lambda u_{t,j}$ and we can present $w_{t,j}$ by McCormick inequalities as 
\begin{equation}
	\label{eq:mccormick} w_{t,\ell} \ge  0, \ w_{t,\ell} \ge \lambda - (1-u_{t,\ell})\lambda^\text{U}, \ w_{t,\ell}\le \lambda^\text{U}u_{t,\ell}, \ \ w_{t,\ell}\le \lambda,	
\end{equation}
%\begin{subequations}
%	\label{eq:mccormick}
%\begin{eqnarray}
%\label{eq:mccormick1}&& w_{t,\ell} \ge  0\\
%&& w_{t,\ell} \ge \lambda - (1-u_{t,\ell})\lambda^\text{U}\\
%&& w_{t,\ell}\le \lambda^\text{U}u_{t,\ell}\\
%\label{eq:mccormick4}&& w_{t,\ell}\le \lambda,
%\end{eqnarray}
%\end{subequations}
where $\lambda^\text{U}$ is an upper bound of $\lambda$.
Therefore, the adjustable DRCC model \eqref{eq:adjustable-drcc} under the Wasserstein ambiguity set is equivalent to the following 0-1 MILP formulation.

\small
\noindent\textbf{MILP3:}
\begin{equation}\label{eq:lp-wass}
	\min_{u_t,\beta_{t,\ell}, x_{t,\ell}, \alpha_t}\left\{ c_\text{sys}\sum_{\ell=1}^{N_{\text{HVAC}}}\beta_{t,\ell} + c_\text{switch}\sum_{\ell=1}^{N_{\text{HVAC}}}u_{t,\ell} + c_t \alpha_t: \ \eqref{eq:x-function}, \ \eqref{eq:absolute}-\eqref{eq:determ-binary}, \ \eqref{eq:z1-mip1},\ \eqref{eq:z1-mip2}, \ \eqref{eq:z1-mip4} - \eqref{eq:z1-mip6}, \ \eqref{eq:z1-mip3-ref},\ \eqref{eq:mccormick} ,\ \eqref{eq:adjustable-alpha-bounds}   \right\}.
	\nonumber
\end{equation}
\normalsize
MILP3 is, however, difficult to solve in certain cases, according to \citet{xie2018distributionally} and \citet{chen2018data}. A bad choice of too large values for the two big-M parameters ($M_n^3$ and $\lambda^\text{U}$) may lead to weak linear relaxations and thus can be detrimental to efficient computation.
Taking into account the RHS uncertainty exploited in Theorem \ref{thm:Wasserstein-LP} and Proposition \ref{prop:strengthen},  next, we derive  a big-M free reformulation. %based on Theorem \ref{thm:Wasserstein-LP}.
% 
%
%
%
% We sort $P_{\text{total},t}$ such that $P_{\text{total},t}^{(1)} \ge P_{\text{total}, t}^{(2)} \ge \cdots \ge P_{\text{total},t}^{(N)}$ and obtain the permutation $\left\{(1), (2), \ldots, (n)\right\}$ of $\{1,2,\ldots,N \}$. %Here, we assume that $P_{\text{total},t}^{(1)} > 0$. Otherwise, 
% We assume that \\
% \noindent \textbf{(A1)} $\sum_{\ell = 1}^{N_{\text{HVAC}}} P_\ell  > P_{\text{total},t}^{(N)}$.\\
% \noindent This is a mild assumption where we assume that the smallest value of the solar PV generation realization can be absorbed fully when all the HVAC units are ON, i.e.,  $u_{\ell, t} = 1, \ \ell = 1,\ldots,N_{\text{HVAC}}$.
% \noindent (A1) there exists a solution $u_t$ such that $a_{(N)} > 0$.\\
\begin{theorem}\label{thm:adjust-wass}
	Under Assumption \ref{assump1},
	% For any solution $u_t$, assume that\\
	% Given a solution $u_t$ and $\alpha_t$, assume that\\
	% \noindent (A1) $a_{(N)} > 0$.\\
	%  \noindent (A2) there exist index $j^*$ such that $a_{(j^*)} > 0 \ge a_{(j^*-1)}$.\\
	% the DR chance constraint \eqref{eq:drcc-cc} under the Wasserstein ambiguity set $\mathcal D_t = \mathcal D_t^2$ is feasible if and only if the following MILP constraints are feasible
	% there exist feasible values of $\Delta_{jk}, \varepsilon_{jk}, o_{\ell j k}$ such that 
	the DR chance constraint \eqref{eq:adjustable-drcc-cc} under the Wasserstein ambiguity set $\mathcal D_t = \mathcal D_t^2$ is   feasible if and only if the following  MILP constraints with auxiliary variables $\Delta_{jk}\in \{0,1\}, \varepsilon_{jk}\in\mathbb R, \tau_{\ell jk} \in \mathbb R, o_{\ell j k} \in \mathbb R $ are feasible. 
	\small
	\begin{subequations}
		\begin{eqnarray}
			\label{eq:adjustable-milp1}&& \sum_{j=1}^{N}\sum_{k = j-1}^{N-1}  \left[ -\frac{1}{N} \sum_{i=j}^k \left(P_{\mathrm{total},t}^{(k+1)}  - P_{\mathrm{total},t}^{(i)}\right)\Delta_{jk} - \sum_{\ell = 1}^{N_{\text{HVAC}}} P_{\ell} \left(o_{\ell jk} - \frac{j-1}{N} \tau_{\ell jk}\right) + P_{\mathrm{total},t}^{(k+1)} \left(\varepsilon_{jk} - \frac{j-1}{N}\Delta_{jk}\right) \right]
			% + \frac{j-1}{N}\sum_{\ell = 1}^{N_{\text{HVAC}}} P_\ell \tau_{\ell jk} - \frac{j-1}{N} P_{\text{total},t}^{(k+1)}\Delta_{jk}\right)
			\le -\delta_t \hspace{9mm}\\
			&& \sum_{j=1}^{N}\sum_{k = j-1}^{N-1} \Delta_{jk} = 1\\
			\label{eq:z-milp-alpha-con1}&& \sum_{j=1}^N \sum_{k=j-1}^{N-1} k \Delta_{jk} \le    \sum_{j=1}^N \alpha_t N \le \sum_{j=1}^N \sum_{k=j-1}^{N-1} (k+1) \Delta_{jk}\\
			%  \alpha_t N \le  \sum_{j=1}^{N}\sum_{k = j-1}^{N-1} (k+1)\Delta_{jk}\\
			% \label{eq:z-milp-alpha-con2} && \alpha_t N \ge \sum_{j=1}^{N}\sum_{k = j-1}^{N-1} k\Delta_{jk}\\
			\label{eq:z-milp-j-con1}&& \sum_{j=1}^{N}\sum_{k = j-1}^{N-1}  P_{\text{total},t}^{(j)} \Delta_{jk} \le \sum_{\ell = 1}^{N_{\text{HVAC}}} P_\ell u_{t,\ell} \le \sum_{j=1}^{N}\sum_{k = j-1}^{N-1}  P_{\text{total},t}^{(j-1)} \Delta_{jk}\\
			%  \label{eq:z-milp-j-con2}&& \sum_{\ell = 1}^{N_{\text{HVAC}}} P_\ell u_{t,\ell} \ge \sum_{j=1}^{N}\sum_{k = j-1}^{N-1}  P_{\text{total},t}^{(j)} \Delta_{jk}\\
			\label{eq:z-milp-mccormick1}&& \varepsilon_{jk} \le \Delta_{jk}, \ \varepsilon_{jk} \le \alpha_t, \ \varepsilon_{jk} \ge \alpha_t + \Delta_{jk} - 1, \ \varepsilon_{jk} \ge 0,\ 0\le j-1 \le k \le N-1\\
			&& o_{\ell jk} \le \varepsilon_{jk}, \ o_{\ell j k} \le u_{t,\ell}, \ o_{\ell jk} \ge \varepsilon_{jk} + u_{\ell} - 1, \ o_{\ell j k} \ge 0, \ 0\le j-1 \le k \le N-1, \ 1\le \ell \le N_{\text{HVAC}}\\
			\label{eq:z-milp-mccormick3}&& \tau_{\ell jk} \le \Delta_{jk}, \ \tau_{\ell j k} \le u_{t,\ell}, \ \tau_{\ell jk} \ge \Delta_{jk} + u_{\ell} - 1, \ \tau_{\ell j k} \ge 0, \ 0\le j-1 \le k \le N-1, \ 1 \le \ell \le N_{\text{HVAC}}%\hspace{7mm}
			\\
			\label{eq:adjustable-milp8}&& \Delta_{jk} \in \{0,1\}, \ 0\le j-1 \le k \le N-1.
		\end{eqnarray}
	\end{subequations}
	\normalsize
\end{theorem}
\proof{Proof: }

See Appendix \ref{sec:proof-adjust-w}  in the online supplement.
\Halmos
\endproof%\qed
The big-M free MILP  reformulation of the adjustable DRCC formulation \eqref{eq:adjustable-drcc} under the $\mathcal D_t^2$ is
% 
%\small
\begin{equation}\noindent\textbf{MILP4: }\label{eq:lp-ad-wass}
	\min_{u_t,\beta_{t,\ell}, x_{t,\ell}, \alpha_t}\left\{ c_\text{sys}\sum_{\ell=1}^{N_{\text{HVAC}}}\beta_{t,\ell} + c_\text{switch}\sum_{\ell=1}^{N_{\text{HVAC}}}u_{t,\ell} + c_t \alpha_t:\ \eqref{eq:x-function},\  \eqref{eq:absolute}-\eqref{eq:determ-binary} ,\ \eqref{eq:adjustable-milp1}-\eqref{eq:adjustable-milp8}, \ \eqref{eq:adjustable-alpha-bounds}\right\}.
	\nonumber
\end{equation}
\normalsize
\begin{remark}
	\textcolor{black}{Theorem \ref{thm:adjust-wass} also applies to general adjustable individual DR binary chance constraints with RHS uncertainty.
		The MILP3 reformulation (with big-M coefficients) for the adjustable chance constraint yields $3N$ binary variables and $3N+ 3N_\text{HVAC} + 1$ constraints. MILP4 has $N(N-1)/2$ binary variables and $3N(N-1)/2 + 3N(N-1)N_\text{HVAC} + 4$ constraints.}
\end{remark}

\section{Computation Setup}
\label{sec:comp-setup}

\begin{figure}[htb]
	\centering
	\includegraphics[width=0.4\linewidth]{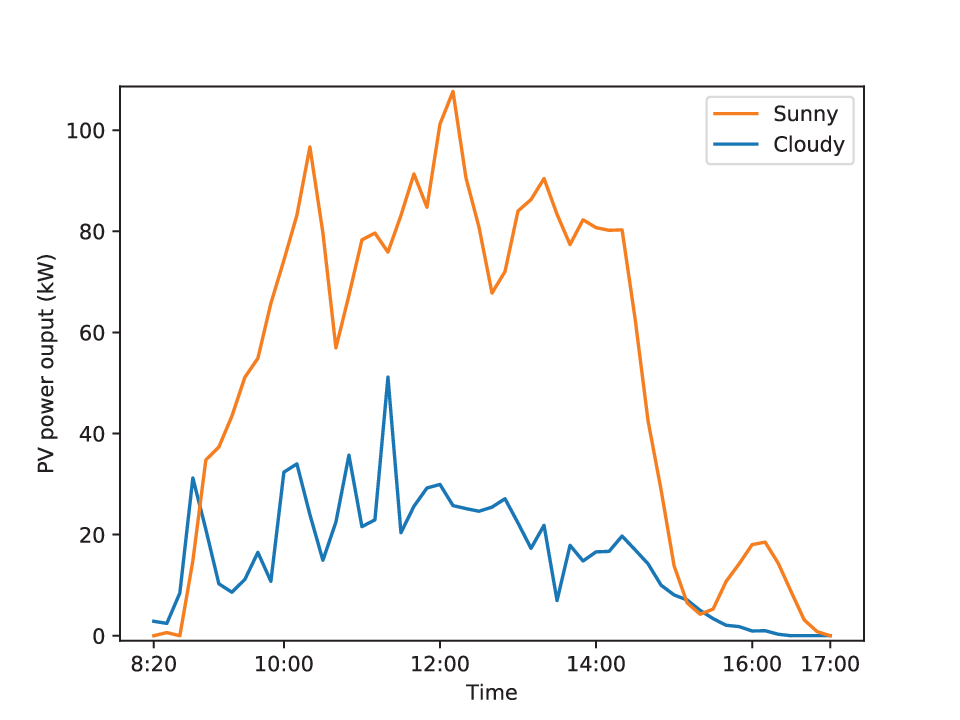}
	\caption{(Color online) PV profile}
	\label{fig:PV}
\end{figure}

We consider a fleet of $N_{\text{HVAC}} = 100$ identical buildings and $N_\text{PV} = 1$ PV panel for $N_p = 53$ periods, every 10 minutes from 8:20 am to 5:00 pm over a day.
We consider two typical weather conditions: a sunny day and a cloudy day. 
The PV power output data {\textcolor{blue}{(available at \url{https://drive.google.com/drive/folders/1ERACqKeP2yYcwzbsvgTbQ13uRqFQD4Gi?usp=sharing}.)}}
of $P_{\text{PV},t}$, as shown in Figure \ref{fig:PV}, is collected from a 13 kW PV panel located on the rooftop of the Distributed Energy Communication \& Control (DECC) laboratory at Oak Ridge National Laboratory (ORNL) in Tennessee.  We scale the PV output to be compatible with the aggregate of 100 residential HVAC systems (connected via a same step-down transformer).
For each building, a random initial room temperature is generated by a uniform distribution between 23.10\textdegree{C} and  23.15\textdegree{C}. The set-point $x_\text{ref}$ is $23.0\text{\textdegree{C}}$ and the comfort band $[x_\text{min}, x_\text{max}]$ is $[21.5, 24.5]$\textdegree{C}. {The building parameters are set as $A_\ell = 0.9914$, $B_\ell = -0.6767$, $G_\ell = (4.3\mathrm{e}{-5}, 
	%	4.3e-5, 
	0.0086)^\top$ for $\ell=1,\ldots,N_{\text{HVAC}}$}. Each HVAC system consumes 3.5 kW when it is ON. The objective costs are set as follows: $c_\text{sys} = 1.0$, and $c_\text{switch} = 1.0$. 
In our studies we consider   the two DRCC formulations, assuming unknown distributional information of the PV output generation, and the benchmark chance-constrained formulation, which are referred to as: 1. CC: stochastic chance-constrained formulation;
2. DRCC-M: DRCC  under the moment-based ambiguity set;
3. DRCC-W: DRCC  under the Wasserstein ambiguity set. 
%\begin{enumerate}
%    \item CC: stochastic chance-constrained formulation;
%    \item DRCC-M: DRCC formulation under the moment-based ambiguity set;
%    \item DRCC-W: DRCC formulation under the Wasserstein ambiguity set. 
%\end{enumerate}

We generate   $N=100$ i.i.d. samples (i.e., in-sample data) of $\widetilde{P}_{\text{PV},t}$   following uniform distributions  \citep{gaunt2017voltage} with the mean $P_{\text{PV},t}$ (shown in Figure \ref{fig:PV}) and half range of $0.15P_{\text{PV},t}$ for $t=1,\ldots,N_p$. We optimize the CC model and construct the ambiguity set of the DRCC-W model by using all $N=100$ samples, and for DRCC-M, only 10   samples are randomly picked from the $N$ samples to calculate the empirical mean and covariance. 
%  $P_{\text{PV},t}^n, \ n =1,\ldots,N$ following uniform distributions (\textcolor{red}{need some reference}) with the mean $P_{\text{PV},t}$ (shown in Figure \ref{fig:PV}) and half range $0.15P_{\text{PV},t}$ for $t=1,\ldots,N_p$.
With the optimal schedules obtained by solving different models, we generate 10 sets of $N^\prime = 1000$ i.i.d. samples (i.e., out-of-sample data) from the same uniform distribution to evaluate the out-of-sample performance of each schedule. All  models are computed
in Python 3.7.5 using Gurobi 9.0.0. The computations are performed on a Windows 10 Pro machine
with Intel(R) Core(TM) i7-8700 CPU 3.20 GHz and 16 GB memory. 

\section{Studies on the DRCC Models}
\label{sec:comp}
%We first compute solutions for all three models with risk level $\alpha_t \in \{0.1,0.2\}, \ t=1,\ldots,N_p$. 
In the DRCC-M model, we set the parameters of the moment-based ambiguity set $(\gamma_1,\gamma_2) = (0,1)$; in the DRCC-W model, we set the radius parameter $\delta_t = 0.02$.
In particular, we solve the DRCC-W model by using the two MILP formulations derived in Section \ref{sec:drcc-ref-wasserstein}. The comparison of the computation time and optimality gaps  is presented in Section \ref{sec:W-CPU}.  In Section \ref{sec:In-sample}, we present the solution details of the three models, including the tracking performance and resulting room temperatures. In Section \ref{sec:out-of-sample}, we present the out-of-sample performance of optimal solutions obtained. Furthermore,
we study the sensitivity of the out-of-sample performance on the in-sample data size.
% the impact of the in-sample test sample size on the out-of-sample performance. To demonstrate the performance sensitivity of the sample size used for computation, we show the performance comparison by solving the three models with only 10 samples and with all 100 samples.
% the chance-constrained models using the benchmark chance-constrained formulation and two DRCC formulations, assuming unknown distributional information of the PV output generation. We refer the three models as:
% \begin{enumerate}
%     \item CC: stochastic chance-constrained formulation;
%     \item DRCC-M: DRCC formulation under moment-based ambiguity set;
%     \item DRCC-W: DRCC formulation under Wasserstein ambiguity set. 
% \end{enumerate}

% We optimize CC and DRCC-W models by using all $N$ samples generated
\subsection{CPU Time and Optimality Gaps}
\label{sec:W-CPU}
%In this section, we solve the DRCC-W model using both MILP1 and MILP2 formulations proposed in Section \ref{sec:drcc-ref-wasserstein}. 
\textcolor{black}{
	%As later shown in Section \ref{sec:out-of-sample},  increasing the number of samples improves the performance. 
	In this section, we show the computational performance of  the two MILPs (MILP1 and MILP2 in Section \ref{sec:drcc-ref-wasserstein}) and the MILP reformulation (denoted as ``MILP-H'') proposed in \citet{ho2020distributionally}  under different sample sizes. The MILP2 formulation is strengthened by the techniques in Proposition \ref{prop:strengthen}.}
In particular, we generate 10 in-sample sets of the sunny weather following the description in Section \ref{sec:comp-setup}, with the size of $N = 100$ and $N=500$, respectively. \textcolor{black}{For further comparing MILP2 and MILP-H, we generate 10 more in-sample sets with larger sample size $N=3000$.} Each instance contains $N_p=53$ periods.  Given an initial room temperature of the first period, for each instance, we sequentially solve the remaining periods by using the resulting room temperature from previous periods as an initial room temperature. The CPU time limit is  100 seconds for each period. We test various risks $1-\alpha_t \in \{80\%, 90\%\}$ and Wasserstein radii $\delta_t \in \{ 0.02, 0.2\}$.

\textcolor{black}{In Table \ref{tab:Wasserstein}, for each instance, we report the total CPU time of solving all 53 periods. If any period cannot be solved to optimality within the time limit, in the parentheses after the CPU times, we report the number of periods that cannot be solved within the time limit and the average optimality gaps of them. Except for $1-\alpha_t = 90\%, \ \delta=0.2$ and $N=100$, for all other cases, the basic MILP1 has instances of which some periods cannot be solved with the time limit. In contrast, MILP-H and MILP2 solve all instances within much shorter time and (or) much smaller optimality gaps. For example, when $1-\alpha_t = 80\%,\ \delta = 0.02,\ N=500$, MILP1 solves instances using 2493.35 seconds and terminates with a 14.83\% optimality gap, on average. While both MILP-H and MILP2 solve the instances optimally within 5 seconds. MILP2 is more effective when the Wasserstein ball's radius $\delta$ is small. When the radius is larger $\delta = 0.2$ with high $1-\alpha_t = 90\%$, MILP-H yields shorter CPU times and (or) smaller gaps. If the sample size becomes larger (N=3,000), MILP-H starts to outperform HILP2 with smaller $1-\alpha_t = 80\%$.}
%
% Table generated by Excel2LaTeX from sheet 'Sheet3'
%
\begin{table}[htbp]
	\centering
	\textcolor{black}{
		\caption{Comparison of CPU time (in seconds) and optimality gaps of high risk requirement $1-\alpha_t$}
		\resizebox{.91\textwidth}{!}{%
			\begin{tabular}{ccc|lll|lll|ll}
				\hline
				\multirow{2}[2]{*}{$1-\alpha$} & \multirow{2}[2]{*}{$\delta$} &       & \multicolumn{3}{c|}{N = 100} & \multicolumn{3}{c|}{N = 500} & \multicolumn{2}{c}{N = 3000} \\
				&       & Instance & \multicolumn{1}{c}{MILP1} & \multicolumn{1}{l}{MILP-H} & \multicolumn{1}{l|}{MILP2} & \multicolumn{1}{c}{MILP1} & \multicolumn{1}{l}{MILP-H} & \multicolumn{1}{l|}{MILP2} & \multicolumn{1}{c}{MILP-H} & \multicolumn{1}{l}{MILP2} \\
				\hline
				\multirow{11}[4]{*}{80\%} & \multirow{11}[4]{*}{0.02} & 1     & 211.34 & 0.83  & \textbf{0.42} & 2522.87 (13, 14.69\%) & 4.43  & \textbf{1.76} & 128.38 (1, 2.92\%) & \textbf{32.85} \\
				&       & 2     & 211.42 (1, 1.30\%) & 5.31  & \textbf{0.37} & 2528.31 (12, 11.16\%) & 2.33  & \textbf{1.67} & 88.51 & \textbf{33.20} \\
				&       & 3     & 109.8 & 0.99  & \textbf{0.42} & 2390.65 (10, 15.68\%) & 2.84  & \textbf{1.60} & 143.16 & \textbf{47.51} \\
				&       & 4     & 146.08 & 100.75 (1, 1.25\%) & \textbf{0.38} & 2387.41 (11, 9.34\%) & 2.56  & \textbf{1.67} & 147.98 (1, 2.99\%) & \textbf{32.76} \\
				&       & 5     & 60.72 & 5.19  & \textbf{0.42} & 2365.25 (9, 18.86\%) & 2.82  & \textbf{1.68} & 89.44 & \textbf{33.37} \\
				&       & 6     & 93.74 & 0.78  & \textbf{0.41} & 2560.81 (10, 15.53\%) & 4.45  & \textbf{1.64} & 56.76 & \textbf{34.13} \\
				&       & 7     & 69.26 & 1.15  & \textbf{0.42} & 2653.68 (12, 21.59\%) & 3.29  & \textbf{1.70} & 93.66 & \textbf{33.95} \\
				&       & 8     & 53.91 & 100.70 (1, 1.33\%) & \textbf{0.37} & 2411.36 (14, 8.36\%) & 3.01  & \textbf{1.60} & 72.63 & \textbf{34.38} \\
				&       & 9     & 64.61 & 0.73  & \textbf{0.35} & 2676.71 (16, 12.29\%) & 2.92  & \textbf{1.71} & 82.45 & \textbf{33.53} \\
				&       & 10    & 153.73 (1, 0.45\%) & 1.28  & \textbf{0.43} & 2436.43 (10, 20.84\%) & 2.55  & \textbf{1.65} & 153.79 & \textbf{33.66} \\
				\cline{3-11}          &       & Avg.  & 117.46 (0.2, 0.18\%) & 21.77 (0.2, 0.25\%) & 0.40  & 2493.35 (11.7, 14.83\%) & 3.12  & 1.67  & 105.68 (0.2, 0.59\%) & 34.93 \\
				\hline
				\multirow{11}[4]{*}{80\%} & \multirow{11}[4]{*}{0.2} & 1     & 846.05 (4, 1.60\%) & 0.87  & \textbf{0.60} & 2833.75 (16, 3.35\%) & 3.18  & \textbf{2.35} & \textbf{47.32} & 51.66 \\
				&       & 2     & 906.48 (4, 1.67\%) & 0.78  & \textbf{0.53} & 3005.26 (16, 3.50\%) & 3.98  & \textbf{2.36} & 102.68 & \textbf{75.78} \\
				&       & 3     & 866.27 (5, 1.28\%) & 0.86  & \textbf{0.56} & 3109.46 (17, 3.32\%) & 10.45 & \textbf{2.14} & \textbf{45.73} & 86.92 \\
				&       & 4     & 830.88 (3, 1.12\%) & 0.99  & \textbf{0.53} & 2809.71 (16, 5.76\%) & 3.72  & \textbf{2.24} & \textbf{52.69} & 82.53 \\
				&       & 5     & 631.79 (3, 1.13\%) & 0.76  & \textbf{0.56} & 3168.74 (18, 2.44\%) & 6.59  & \textbf{2.43} & \textbf{66.68} & 89.73 \\
				&       & 6     & 843.26 (3, 1.02\%) & 0.83  & \textbf{0.57} & 2771.79 (14, 2.49\%) & 3.03  & \textbf{2.17} & \textbf{51.69} & 84.65 \\
				&       & 7     & 874.16 (4, 1.24\%) & 0.92  & \textbf{0.57} & 2938.23 (16, 2.98\%) & 2.83  & \textbf{2.37} & \textbf{46.80} & 98.53 \\
				&       & 8     & 835.96 (5, 1.15\%) & 0.98  & \textbf{0.58} & 2992.62 (20, 3.71\%) & 16.07 & \textbf{2.11} & \textbf{48.68} & 102.53 \\
				&       & 9     & 823.88 (4, 1.65\%) & 0.85  & \textbf{0.54} & 3040.97 (17, 3.89\%) & 3.23  & \textbf{2.33} & \textbf{51.64} & 102.78 \\
				&       & 10    & 688.03 (3, 1.20\%) & 0.98  & \textbf{0.62} & 3049.83 (18, 3.17\%) & 10.43 & \textbf{2.33} & \textbf{53.14} & 84.12 \\
				\cline{3-11}          &       & Avg.  & 814.68 (3.8, 1.31\%) & 0.88  & 0.57  & 2972.04 (16.8, 3.46\%) & 6.35  & 2.28  & 56.71 & 85.92 \\
				\hline
				\multirow{11}[4]{*}{90\%} & \multirow{11}[4]{*}{0.02} & 1     & 390.77 (2, 1.12\%) & 0.68  & \textbf{0.35} & 2174.93 (13, 2.15\%) & 1.21  & \textbf{0.75} & 13.72 & \textbf{9.15} \\
				&       & 2     & 321.15 (2, 1.56\%) & 0.64  & \textbf{0.34} & 2525.91 (18, 1.81\%) & 1.35  & \textbf{0.75} & 15.64 & \textbf{9.28} \\
				&       & 3     & 342.84 (2, 1.75\%) & 0.58  & \textbf{0.35} & 2193.18 (12, 2.96\%) & 1.16  & \textbf{0.76} & \textbf{8.35} & 9.27 \\
				&       & 4     & 312.02 (1, 1.11\%) & 0.61  & \textbf{0.31} & 2428.54 (18, 2.55\%) & 1.43  & \textbf{0.73} & 10.80 & \textbf{8.82} \\
				&       & 5     & 380.92 (3, 1.22\%) & 0.57  & \textbf{0.38} & 2311.06 (14, 2.17\%) & 1.37  & \textbf{0.68} & 13.55 & \textbf{9.23} \\
				&       & 6     & 347.04 (2, 1.62\%) & 0.60  & \textbf{0.33} & 2144.32 (13, 1.88\%) & 1.19  & \textbf{0.75} & \textbf{8.61} & 8.80 \\
				&       & 7     & 491.16 (1, 1.19\%) & 0.64  & \textbf{0.35} & 1958.54 (11, 1.66\%) & 1.07  & \textbf{0.74} & 10.01 & \textbf{8.67} \\
				&       & 8     & 349.96 (1, 1.14\%) & 0.56  & \textbf{0.32} & 2304.81 (15, 2.79\%) & 1.33  & \textbf{0.65} & 12.45 & \textbf{9.13} \\
				&       & 9     & 431.95 (3, 1.39\%) & 0.54  & \textbf{0.35} & 2098.43 (12, 1.98\%) & 1.21  & \textbf{0.77} & 9.58  & \textbf{8.64} \\
				&       & 10    & 332.82 (2, 1.63\%) & 0.48  & \textbf{0.35} & 2355.13 (16, 1.82\%) & 1.28  & \textbf{0.80} & 18.37 & \textbf{8.49} \\
				\cline{3-11}          &       & Avg.  & 370.06 (1.9, 1.37\%) & 0.59  & 0.34  & 2249.49 (14.2, 2.18\%) & 1.26  & 0.74  & 12.11 & 8.95 \\
				\hline
				\multirow{11}[4]{*}{90\%} & \multirow{11}[4]{*}{0.2} & 1     & 216.17 & 0.44  & 0.44  & 1962.85 (8, 2.61\%) & 1.10  & \textbf{0.89} & \textbf{9.84} & 10.32 \\
				&       & 2     & 220.72 & \textbf{0.42} & 0.44  & 1956.90 (7, 2.37\%) & 1.04  & \textbf{0.94} & \textbf{10.50} & 11.05 \\
				&       & 3     & 214.09 & \textbf{0.50} & 0.51  & 2013.63 (10, 2.07\%) & \textbf{0.95} & 0.96  & \textbf{10.14} & 12.42 \\
				&       & 4     & 199.31 & \textbf{0.45} & 0.47  & 1966.34 (10, 2.29\%) & 1.14  & \textbf{0.87} & \textbf{7.11} & 10.89 \\
				&       & 5     & 239.24 & 0.45  & \textbf{0.43} & 1867.35 (9, 2.49\%) & 1.08  & \textbf{0.93} & \textbf{6.52} & 11.15 \\
				&       & 6     & 285.80 & 0.47  & \textbf{0.45} & 1947.05 (11, 2.19\%) & \textbf{0.87} & 0.98  & \textbf{9.47} & 10.63 \\
				&       & 7     & 297.24 & \textbf{0.43} & 0.55  & 1873.74 (7, 2.51\%) & \textbf{0.96} & 0.98  & \textbf{9.26} & 11.96 \\
				&       & 8     & 236.04 & \textbf{0.42} & 0.49  & 1856.42 (8, 2.64\%) & 1.00  & \textbf{0.99} & \textbf{6.22} & 10.56 \\
				&       & 9     & 201.77 & \textbf{0.42} & 0.47  & 1947.70 (10, 1.95\%) & \textbf{0.92} & 100.93 (1, 1.42\%) & \textbf{9.92} & 10.86 \\
				&       & 10    & 254.37 & \textbf{0.42} & 0.45  & 1906.98 (10, 2.37\%) & \textbf{0.94} & 0.99  & \textbf{6.62} & 10.72 \\
				\cline{3-11}          &       & Avg.  & 236.48 & 0.44  & 0.47  & 1929.90 (9, 2.35\%) & 1.00  & 10.95 (0.1, 0.14\%) & 8.56  & 11.06 \\
				\hline
			\end{tabular}%
		}
		\label{tab:Wasserstein}
}\end{table}%

We also provide the average computational performance of the three MILPs when $1-\alpha_t = \{10\%, 20\%\}$, which may be of little interest in practice, but for a fair comparison. 
In Table \ref{tab:cpu-w-small}, MILP2 scales well when $1-\alpha_t$ is small. While MILP1 and MILP-H require much longer time.  For example, when $1-\alpha_t = 10\%,\ \delta = 0.02,\  N=100$, MILP2 finishes under 5 seconds, while MILP1 uses 2379.21 seconds with a  10.41\% gap and MILP-H of 195.82 seconds with a 0.09\% gap. More details of each instance can be found in the online appendix \ref{sec:cpu-w-89}. For the moment-based ambiguity set, the computational details are in the online appendix \ref{sec:cpu-moment}.
%
% Table generated by Excel2LaTeX from sheet 'cpu_wasserstein-89'
\begin{table}[htbp]
	\centering
	\textcolor{black}{
		\caption{Average CPU time (in seconds) and optimality gaps of low risk requirement $1-\alpha_t$}
		\resizebox{\textwidth}{!}{%
			\begin{tabular}{cc|llc|lll|ll}
				\hline
				\multirow{2}[2]{*}{$1-\alpha$} & \multirow{2}[2]{*}{$\delta$} & \multicolumn{3}{c|}{N = 100} & \multicolumn{3}{c|}{N = 500} & \multicolumn{2}{c}{N = 3000} \\
				&       & \multicolumn{1}{c}{MILP1} & MILP-H & MILP2 & \multicolumn{1}{c}{MILP1} & MILP-H & \multicolumn{1}{c|}{MILP2} & \multicolumn{1}{c}{MILP-H} & \multicolumn{1}{c}{MILP2} \\
				\hline
				\multirow{2}[1]{*}{10\%} & 0.02  & 80.31 (0.5, 0.17\%) & 8.60  & \textbf{0.84} & 2379.21 (6.9, 10.41\%) & 195.82 (0.2, 0.09\%) & \textbf{4.42} & 2205.50 (17.6, 4.54\%) & \textbf{53.67 (0.1, 0.04\%)} \\
				& 0.2   & 309.28 (1.8, 0.38\%) & 187.39 (1.5, 0.35\%) & \textbf{3.97} & 3199.71 (11.9, 20.31\%) & 387.83 (2, 0.39\%) & \textbf{10.01} & 2029.92 (14.9, 4.58\%) & \textbf{247.49 (1.3, 0.27\%)} \\
				\hline
				\multirow{2}[1]{*}{20\%} & 0.02  & 113.34 (0.7, 0.16\%) & 29.36 (0.2, 0.06\%) & \textbf{0.86} & 2979.62 (13.9, 26.00\%) & 186.22 (0.2, 0.06\%) & \textbf{14.70 (0.1, 0.01\%)} & 2029.92 (14.9, 4.58\%) & \textbf{35.54} \\
				& 0.2   & 345.01 (1.1, 0.36\%) & 25.9  & \textbf{1.58} & 3443.04 (12.9, 25.22\%) & 268.21 (0.6, 0.22\%) & \textbf{53.04 (0.3, 0.06\%)} & 3388.08 (27.3, 4.06\%) & \textbf{212.76 (1.2, 0.34\%)} \\
				\hline
			\end{tabular}	\label{tab:cpu-w-small}%
		}
	}
\end{table}%
\subsection{Tracking Performance and Room Temperatures}
\label{sec:In-sample}

In this section, we present the results under the sunny weather condition and the results for the cloudy weather condition are in the online appendix \ref{supplement-sec:tracking-temperature-cloudy}.
For the sunny weather, we solve three models  for all $N_p$ periods under  the sunny weather condition with $N_\text{HVAC} = 100$ HVAC units. For the cloudy weather,  as the PV generation is relatively lower than that of the sunny day (as shown in Figure \ref{fig:PV}), 
%The results under cloudy weather are included in  the online appendix \ref{supplement-sec:tracking-temperature-cloudy}.
fewer number of HVAC units are enrolled  to keep the total consumption compatible with the magnitude of the local PV generation \citep{dong2018model}. 
We remark that our model is flexible to incorporate fleet sizing decisions based on the nameplate capacity of HVAC devices and local solar PV generation scales. The details are in the online appendix \ref{sec:new-model}.

%For all three models, we solve them  for all $N_p$ periods under both the sunny weather and cloudy weather conditions, respectively. 
%{We note that the PV generation of the cloudy day is relatively lower than that of the sunny day (see in Figure \ref{fig:PV}).  Consequently, when solving instances of the cloudy weather condition, to keep the total number of enrolled HVAC devices compatible with the magnitude of local solar PV generation \citep{dong2018model}, we consider enrolling a subset of 35 residential HVAC units. 
%	That is, $N_{\text{HVAC}} = 35$ under the cloudy weather. While, we consider enrolling all the 100 HVAC units for the sunny day. \textcolor{black}{We remark that our model is flexible to incorporate the decision of the fleet size, which can be decided based on the nameplate capacity of HVAC devices and local solar PV generation scales.} The details are in the online appendix \ref{sec:new-model}.   }  

\begin{figure*} 
	\centering
	\begin{subfigure}{0.4\textwidth}
		\includegraphics[width=.9\linewidth]{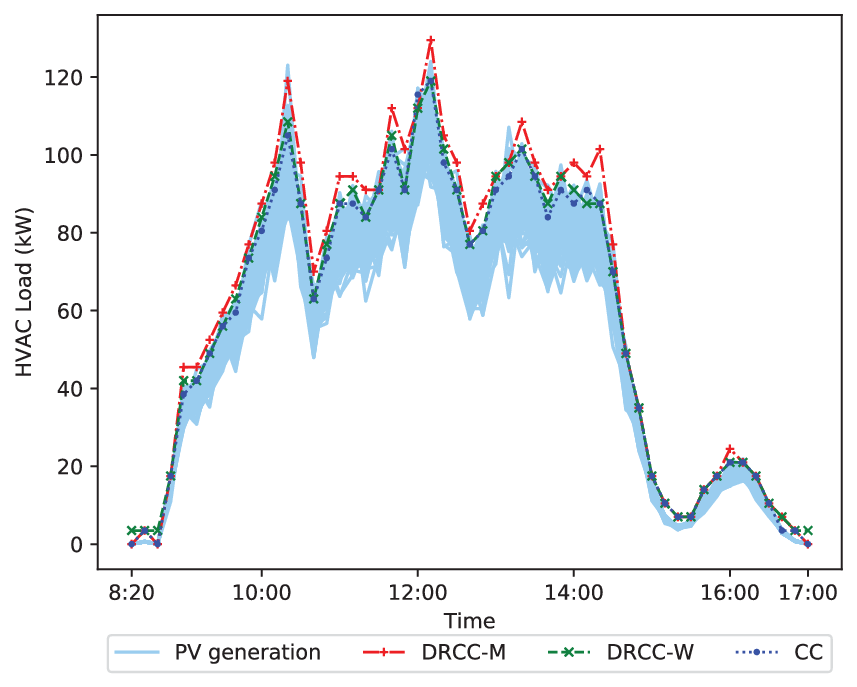} \caption{PV profile tracking} 
		\label{fig:PVtracking-sunny}
	\end{subfigure} 
	\begin{subfigure}{0.4\textwidth}
		\includegraphics[width=\linewidth]{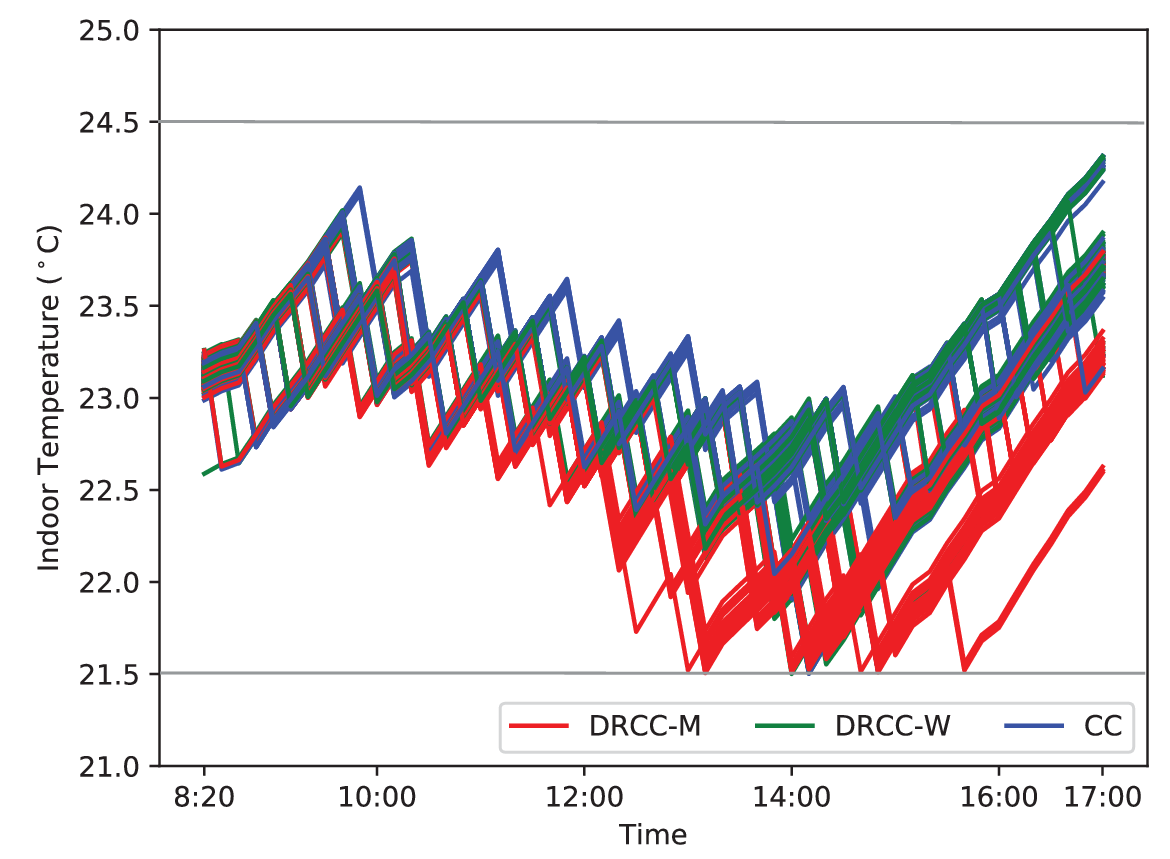} \caption{Room temperature} 
		\label{fig:temp-sunny}
	\end{subfigure} 
	\caption{(Color online) PV profile tracking and room temperatures of 100 buildings under sunny  weather}
	\label{fig:sunny-tracking-temp}
\end{figure*}

We present the overall tracking performance under sunny weather using 100 ON/OFF HVAC devices, i.e., $\sum_{j=1}^{N_{\text{HVAC}}} P_ju_{t,j}, \ t =1,\ldots,N_p$ of the three models in Figure \ref{fig:PVtracking-sunny}. 
The shaded blue areas in the background is the plot of 100 PV generation samples used for solving the models. Overall, all three models track the PV generation well.
%In Figure \ref{fig:PVtracking-cloudy} under the cloudy weather, in addition to the three models, we also present a ``Base'' model, where we solve a control problem without considering PV generation at all. That is, we solve a deterministic model of \eqref{eq:determ} by removing the last term of signal deviation of PV signal from the objective function. 
%
%Under both weather conditions, all three models track the PV generation well. 
Most of the periods, the two DRCC models provide higher HVAC loads than the stochastic CC model given that the DRCC models take into account ambiguous probability distributions and thus the solutions are more conservative.
%In Figure \ref{fig:PVtracking-sunny}, under the sunny weather, all three models track the PV generation well. The two DRCC models provide higher HVAC loads than the stochastic CC model as the DRCC models are more robust by considering ambiguous probability distributions. 
The DRCC-M model yields higher HVAC loads as the DRCC-M model is generally more conservative. \textcolor{black}{In Figure \ref{fig:temp-sunny}, we present the resulting room temperatures of all $N_{\text{HVAC}} = 100$ buildings over $N_p$ periods for all three models. All the indoor temperatures are maintained within the desired comfort band $[21.5, \ 24.5]$\textdegree{C}. The DRCC models provide cooler room temperatures for most buildings which is an immediate result of turning on more HVAC units as shown in Figure \ref{fig:PVtracking-sunny}.}

%In Figure \ref{fig:PVtracking-cloudy} under the cloudy weather, we observe the similar patterns. 
%The three models track the PV generation well at the beginning before 10 pm. For the rest of the day, there are several periodic jumps, which shows the similar pattern as in the Base model. In particular, when the PV generation output is small, the optimal schedule of all three models are dominated by the optimal schedule without considering PV generation most of the time.
%
%

\subsection{Out-of-Sample Performance}
\label{sec:out-of-sample}
After solving all the models and obtaining the optimal schedules, we fix them in ten out-of-sample data sets, each consisting of $N^\prime = 1000$ samples. For each data set, an out-of-sample probability is calculated as the ratio of the number of scenarios, where the PV generation is consumed locally (i.e., the total HVAC load is more than the PV generation), to the total number of scenarios $N^\prime$. The performance is measured by the 95th percentile of  probabilities of the ten out-of-sample sets.

\begin{figure*} 
	\centering
	\begin{subfigure}{0.45\textwidth}
		\includegraphics[width=\linewidth]{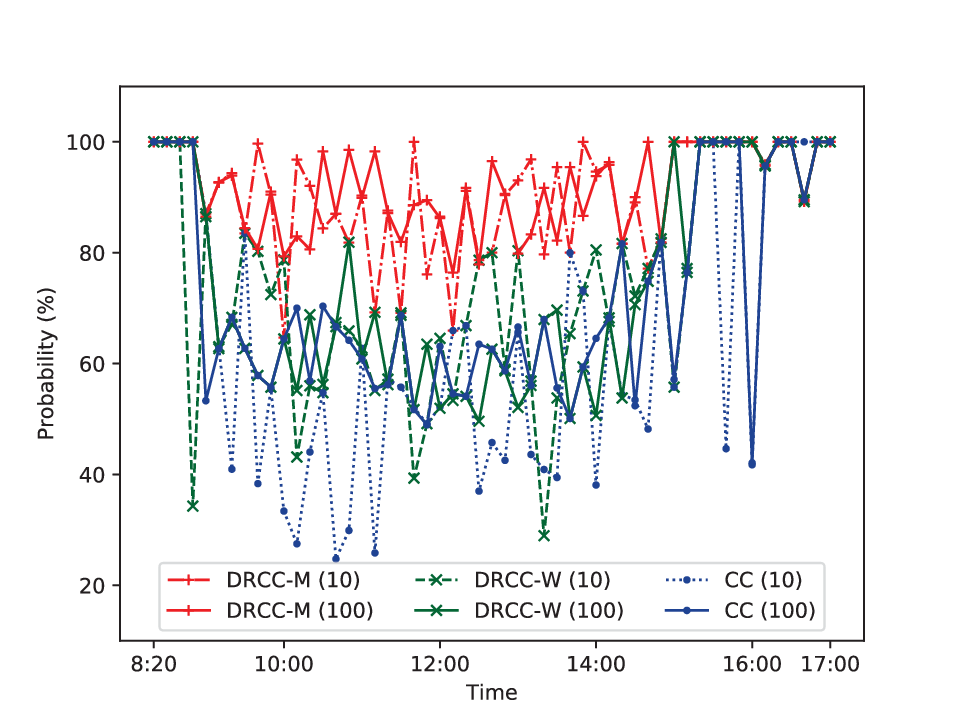} \caption{$1-\alpha_t = 50\%$} 
		\label{fig:alpha5}
	\end{subfigure} 
	\begin{subfigure}{0.45\textwidth}
		\includegraphics[width=\linewidth]{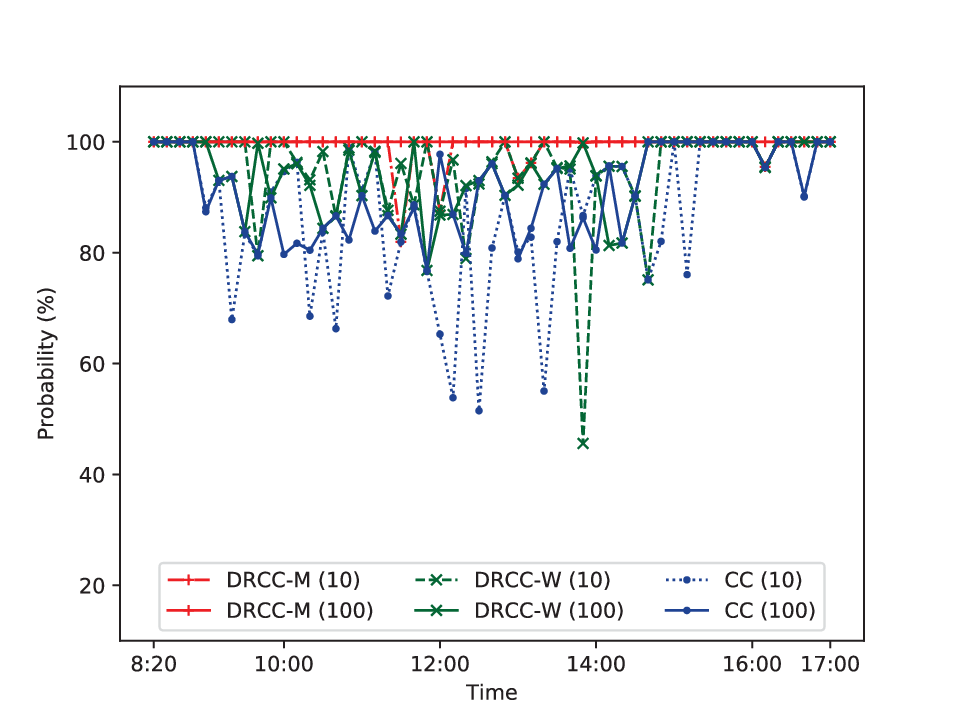} \caption{$1-\alpha_t = 80\%$} 
		\label{fig:alpha8}
	\end{subfigure}
	\caption{(Color online) Probability of locally consuming PV generation under sunny weather}
	\label{fig:prob_alpha}
\end{figure*}

To study the impact of the number of the samples used for solving the models, we consider two choices of the risk parameter $1-\alpha_t$: 80\% and 50\% under the sunny weather condition.
We solve the three models using 10 samples and 100 samples, respectively, with the two risk levels. The out-of-sample performance is presented in Figure \ref{fig:prob_alpha}. When the sample size is small $N = 10$, except for the DRCC-M model, both the DRCC-W and CC models fail to achieve the required risk level. For example, in Figure \ref{fig:alpha5}, the DRCC-W and the CC models perform below the required risk level $1-\alpha_t = 50\%$  between 11 am and  12 pm when using only 10 samples. However, with more samples $N = 100$,  the DRCC-W model is always above $1-\alpha_t = 50\%$, while, the CC model sometimes still fails to achieve the required risk level. Overall, the DRCC models perform better than the CC model, which is consistent with the previous observations. The DRCC-W models are more sensitive to the number of samples  compared to the DRCC-M model. A comparison of sunny and cloudy weather is included in the online appendix \ref{supplement-sec:out-of-sample}.
%
%
%
%
% \begin{figure*} 
% \centering
%   \begin{subfigure}{0.45\textwidth}
%     \includegraphics[width=\linewidth]{out-of-sample-sunny_Gaussian.pdf} \caption{Sunny} 
%   \end{subfigure} 
%   \begin{subfigure}{0.45\textwidth}
%     \includegraphics[width=\linewidth]{out-of-sample-cloudy_Gaussian.pdf} \caption{Cloudy} 
%   \end{subfigure}
%   \caption{Probability of locally consuming PV generation under Gaussian noise}
%   \label{fig:prob_Gaussian}
% \end{figure*}

% \begin{figure*} 
% \centering
%   \begin{subfigure}{0.45\textwidth}
%     \includegraphics[width=\linewidth]{out-of-sample-sunny_GaussianV.pdf} \caption{Sunny} 
%   \end{subfigure} 
%   \begin{subfigure}{0.45\textwidth}
%     \includegraphics[width=\linewidth]{out-of-sample-cloudy_GaussianV.pdf} \caption{Cloudy} 
%   \end{subfigure}
%   \caption{Probability of locally consuming PV generation under VARIANT Gaussian noise}
%   \label{fig:prob_Gaussianv}
% \end{figure*}
%
%
%
%
\section{Studies on the Adjustable DRCC Models}
\label{sec:comp-adjustable}
In this section, we focus on the adjustable variants of the two DRCC models, where we consider the risk parameter $\alpha_t$ as a variable  than a known parameter. Specifically, we solve the DRCC-W models using the two MILP formulations proposed in Section \ref{sec:adjustable-wasserstein}. The computational details of the adjustable models under the Wassestein set and the moment-based set are in the online appendices \ref{sec:wass-adjustable-cpu} and  \ref{sec:cpu-moment}.  In Section \ref{sec:wass-adjustable-cpu}, the  comparison of CPU time and optimality gaps can be found. %The solution details and sensitivity analysis of the risk level are presented in Section \ref{sec:wass-adjustable-alpha}. 

%\subsection{Sensitivity Analysis of Risk-Level Coefficient Cost}
We now study the impact of the coefficient cost $c_t$ in the objective \eqref{eq:adjustable-drcc-obj} of risk level $\alpha_t$ for both DRCC-M and DRCC-W models. We vary the coefficient from 10 to 20 in increments of 2. For the adjustable models, we observe that the optimal risk levels for the periods between 9:30 am and 2:30 pm are relatively lower. That is, during these periods, it is harder to consume all PV generation, which is consistent with the observations in Section \ref{sec:out-of-sample} for the non-adjustable models. Therefore, we focus on these periods in between 9:30 am  and 2:30 pm of a day with an increment of one hour. 

\begin{figure*} 
	\centering
	\begin{subfigure}{0.4\textwidth}
		\includegraphics[width=\linewidth]{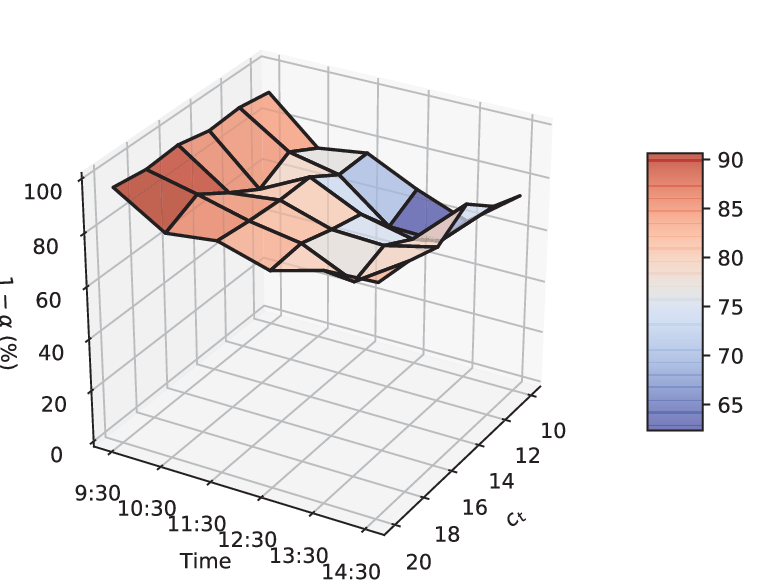} \caption{DRCC-M} 
		\label{fig:adjust_M}
	\end{subfigure} 
	\begin{subfigure}{0.4\textwidth}
		\includegraphics[width=\linewidth]{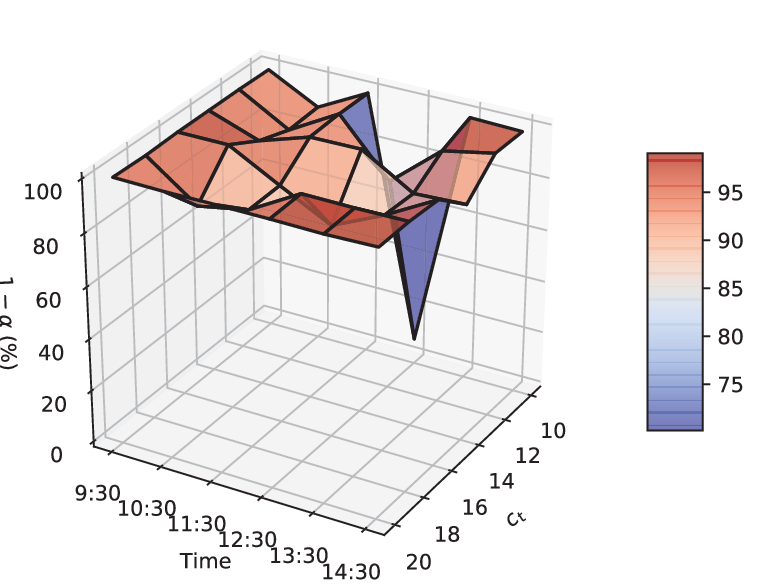} \caption{DRCC-W} 
		\label{fig:adjust_W}
	\end{subfigure}
	\caption{(Color online) The impact of the coefficient cost $c_t$ on the risk level }
	\label{fig:adjust_alpha}
\end{figure*}

In Figure \ref{fig:adjust_alpha}, we show the risk levels under different choices of coefficient $c_t$ for the DRCC-M and DRCC-W models.
In Figures \ref{fig:adjust_M} and \ref{fig:adjust_W}, the value $c_t$ of the risk level increases as larger coefficient is set. The graph in Figure \ref{fig:adjust_M} of the DRCC-M model is flatter as $c_t$ changes, because the performance of the DRCC-M model is less sensitive to the choice of $c_t$. While, as shown in Figure \ref{fig:adjust_W}, the DRCC-W model has a  low risk level  when $c_t$ is small and a higher risk level when $c_t$ is large. 

\begin{figure*} 
	\centering
	\begin{subfigure}{0.4\textwidth}
		\includegraphics[width=\linewidth]{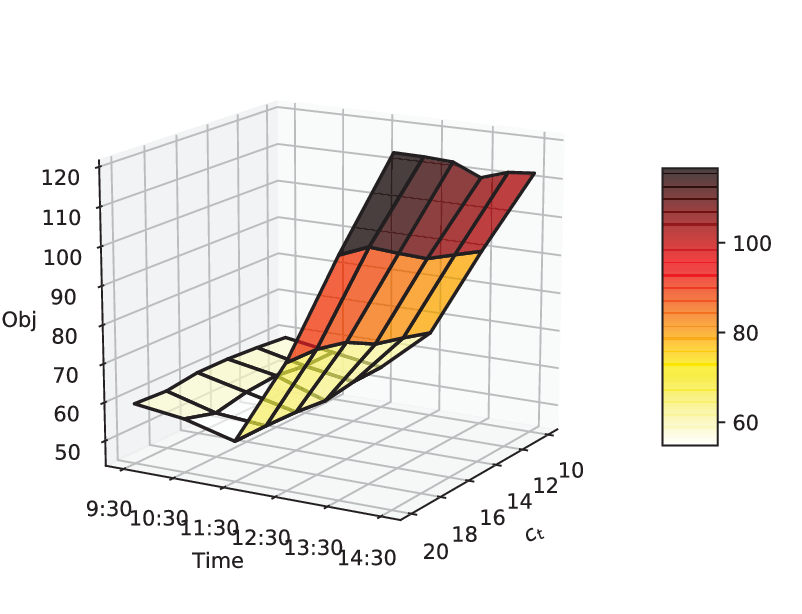} \caption{DRCC-M} 
		\label{fig:adjust_M_cost}
	\end{subfigure} 
	\begin{subfigure}{0.4\textwidth}
		\includegraphics[width=\linewidth]{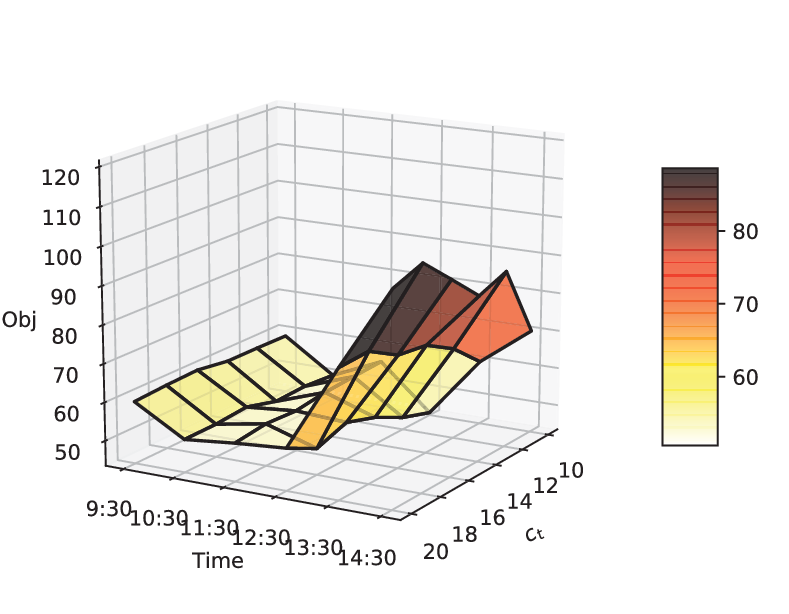} \caption{DRCC-W} 
		\label{fig:adjust_W_cost}
	\end{subfigure}
	\caption{(Color online) The impact of the coefficient cost $c_t$ on the objective cost}
	\label{fig:adjust_cost}
\end{figure*}

In Figure \ref{fig:adjust_cost}, we show the objective costs for the  various coefficient costs $c_t$. We see that both models yield higher objective costs as $c_t$ increases.
We also see that the DRCC-M model  yields higher objective costs under large $c_t$ even when the risk level is set lower than the DRCC-W model. For example, at 2:30 pm under $c_t = 20$,  DRCC-M sets the risk level around 85\% with an objective cost  around 120, while DRCC-W achieves a higher risk level above 95\% with an objective cost less than 90. This is expected as the DRCC-M model is more conservative and requires more HVAC units to run to achieve a similar risk level.
% 
%
%
%\label{sec:wass-adjustable-alpha}
%
%\begin{figure*} 
%	\centering
%	\begin{subfigure}{0.4\textwidth}
%		\includegraphics[width=\linewidth]{figures/Adjustable_alpha_contourM.pdf} \caption{DRCC-M} 
%		\label{fig:adjust_M}
%	\end{subfigure} 
%	\begin{subfigure}{0.4\textwidth}
%		\includegraphics[width=\linewidth]{figures/Adjustable_alpha_contourW.pdf} \caption{DRCC-W} 
%		\label{fig:adjust_W}
%	\end{subfigure}
%	\caption{(Color online) The impact of the coefficient cost $c_t$ on the risk level }
%	\label{fig:adjust_alpha}
%\end{figure*}

\label{sec:wass-adjustable-alpha}
%\begin{figure*} 
%	\centering
%	\begin{subfigure}{0.4\textwidth}
%		\includegraphics[width=\linewidth]{figures/Adjustable_cost_contourM.pdf} \caption{DRCC-M} 
%		\label{fig:adjust_M_cost}
%	\end{subfigure} 
%	\begin{subfigure}{0.4\textwidth}
%		\includegraphics[width=\linewidth]{figures/Adjustable_cost_contourW.pdf} \caption{DRCC-W} 
%		\label{fig:adjust_W_cost}
%	\end{subfigure}
%	\caption{(Color online) The impact of the coefficient cost $c_t$ on the objective cost}
%	\label{fig:adjust_cost}
%\end{figure*}
%
%
%
\section{Conclusions}
\label{sec:conclusions}

In this paper, we formulated a single-period BLC problem as a DRCC problem with uncertain PV generation under both the moment-based  and Wasserstein ambiguity sets as DRCC-M and DRCC-W, respectively. For the moment-based ambiguity set, we reformulated the DRCC problem as an MILP reformulation. For the Wasserstein ambiguity set, we provided an MILP1 reformulation and a more compact MILP2 reformulation by exploiting the CVaR interpretation. The results for DRCC apply to general individual DR chance constraints with RHS uncertainty. 

By considering the risk level as a decision variable,  we also proposed adjustable DRCC formulations that determine the optimal risk level of the chance constraint to balance the total cost and overall performance. For the moment-based ambiguity set, we developed an exact solution approach by solving two SOCP problems. For the Wasserstein ambiguity set, we derived a big-M MILP reformulation and   a big-M-free MILP reformulation. The results of the adjustable variants also apply to general individual \textcolor{black}{binary} chance constraints with RHS uncertainty. 

Extensive computational studies were conducted  on the non-adjustable and adjustable DRCC models under the two ambiguity sets. We found that the DRCC models achieve better out-sample performance while maintaining the indoor temperature within a desired comfort band. Specifically, the DRCC-M model requires fewer samples to achieve the required risk level. The DRCC-W model performs well when using enough many samples and only requires modest CPU time when solving the more compact reformulation MILP2. Furthermore, for the adjustable DRCC problem, we find that the DRCC-M is less sensitive to the choice of $c_t$ while it may require higher objective costs.  % In summary, the DRCC models 
\bibliographystyle{informs2014}
\bibliography{YilingRef}

\newpage

\begin{appendices}

	\section{Proof of Theorem \ref{thm:socp1}}\label{sec:proof-thm-socp1}
	\proof{Proof of Theorem \ref{thm:socp1}:}
	
	\textcolor{black}{We first rearrange constraint  \eqref{eq:linear-ref1} as 	
		\begin{equation}\label{eq:linear-ref1-rearrange}
			\sum_{\ell=1}^{N_s}P_\ell u_{t,\ell} - \theta_t - \sigma_t\sqrt{\gamma_1} \ge  \left( \sqrt{\frac{1-\alpha_t}{\alpha_t} (\gamma_2 - \gamma_1)}\right) \sigma_t. \end{equation}
		Given  the nonnegative RHS of \eqref{eq:linear-ref1-rearrange},  the left-hand side is implicitly enforced nonnegative, which results in constraint \eqref{eq:thm-soc7}. }By squaring both side, we obtain
	\begin{equation}\label{eq:thm1-inter1}(\gamma_2 - \gamma_1)\sigma_t^2 \le \alpha_t \left[ \left(\sum_{\ell = 1}^{N_{\text{HVAC}}} P_\ell u_{t,\ell} - \theta_t - \sigma_t\sqrt{\gamma_1}\right)^2 + (\gamma_2 - \gamma_1)\sigma_t^2  \right].\end{equation}
	The RHS of \eqref{eq:thm1-inter1} is nonlinear.
	We let $d =  \left(\sum_{\ell = 1}^{N_{\text{HVAC}}} P_\ell u_{t,\ell} - \theta_t - \sigma_t\sqrt{\gamma_1}\right)^2 + (\gamma_2 - \gamma_1)\sigma_t^2$. Then  inequality \eqref{eq:thm1-inter1} is equivalent to $(\gamma_2 - \gamma_1) \sigma_t^2 \le \alpha_t d$ which is equivalent to  \eqref{eq:thm-soc1}. To linearize $d$, we define $g_{ij} := u_{t,i}u_{t,j}$ by McCormick inequalities \eqref{eq:thm-soc3}. We conclude the proof.
	\Halmos
	\endproof
	
	\section{Proof of Theorem \ref{thm:socp2}} \label{sec:proof-thm-socp2}
	\proof{Proof of Theorem \ref{thm:socp2}:}
	To show the equivalence, we need to show that (i) constraint \eqref{eq:linear-ref2} implies constraints \eqref{eqn:soc-1}--\eqref{eqn:soc-4} and (ii) constraints \eqref{eqn:soc-1}--\eqref{eqn:soc-4} imply constraint \eqref{eq:linear-ref2}.
	
	\begin{enumerate}
		\item[(i)] \eqref{eq:linear-ref2} $\rightarrow$ \eqref{eqn:soc-1}--\eqref{eqn:soc-4}.
		
		Given a solution $(u_t^*, \alpha_t^*)$ that satisfies \eqref{eq:linear-ref2}, we let $\phi^* = \sqrt{1/\alpha_t^*}$, $w^* = \sqrt{\phi^*}$, $q^* = 1/w^*$. Then $(u_t^*, \alpha_t^*, \phi^*, w^*, q^*)$ is a solution to 
		\eqref{eqn:soc-1}--\eqref{eqn:soc-4}. 
		
		\item[(ii)] \eqref{eqn:soc-1}--\eqref{eqn:soc-4} $\rightarrow$ \eqref{eq:linear-ref2}.
		
		We notice that \eqref{eqn:soc-2} is equivalent to 
		\begin{equation}
			\label{eqn-ref-5}   \alpha_t \phi \ge q^2, \ \phi \ge 0;
		\end{equation}
		and \eqref{eqn:soc-4} can be rewritten as 
		\begin{equation}
			\label{eqn-ref-6}   q \ge \frac{1}{w}, \ w \ge 0.
		\end{equation}
		Combining \eqref{eqn-ref-5}, \eqref{eqn-ref-6}, and  \eqref{eqn:soc-3}, we have $\alpha_t \phi \ge 1/\phi $, which is further equivalent to 
		\begin{equation}\label{eqn-ref-7}
			% \alpha_t v \ge \frac{1}{v}, %\Leftrightarrow 
			v \ge \sqrt{\frac{1}{\alpha_t}}.
		\end{equation} 
		Combining \eqref{eqn:soc-1} and \eqref{eqn-ref-7}, we conclude that constraint \eqref{eqn:soc-1} implies \eqref{eq:linear-ref2}.
		\Halmos 
		\endproof
	\end{enumerate}
	
	\section{Proof of Theorem \ref{thm:moment-adjust}}\label{sec:proof-thm-moment-adjust}
	\proof{Proof of Theorem \ref{thm:moment-adjust}:}
	In the adjustable DRCC model \eqref{eq:adjustable-drcc}, we
	replace the adjustable DR chance constraint \eqref{eq:adjustable-drcc-cc} with the following convex reformulation.
	% Inspired by Theorem 14 in \cite{xie2019optimized}, 
	% when ${\gamma_1}/{\gamma_2} \le \alpha_t \le 0.75$, we can derive a more compact reformulation, which incorporates only one auxiliary variable $r$, by replacing \eqref{eq:adjustable-drcc-cc} with the following convex reformulation
	\begin{eqnarray}
		\label{eq:thm-qc1}&& \sum_{\ell=1}^{N_{\text{HVAC}}} P_\ell u_{t,\ell} - \theta_t - \sigma_t\sqrt{\gamma_1} \ge \sigma_t \sqrt{\gamma_2 - \gamma_1} r\\
		\label{eq:thm-qc2}&& r \ge \sqrt{\frac{1-\alpha_t}{\alpha_t}}.
	\end{eqnarray}
	The RHS of constraint \eqref{eq:thm-qc2} is convex when $0 \le  \gamma_1 / \gamma_2 \le \alpha_t \le 0.75$ as the second-order derivative $(3-4\alpha_t)(1-\alpha_t)^{-3/2}\alpha_t^{-5/2}/4$ is non-negative.  
	We can further construct an outer approximation of the reformulation \eqref{eq:thm-qc1}--\eqref{eq:thm-qc2} by replacing \eqref{eq:thm-qc2} with  
	\begin{equation}\label{eq:thm-outer-approx} 2r \ge \sqrt{1\over{\alpha_t}}. \end{equation} 
	Constraint \eqref{eq:thm-outer-approx} is implied by \eqref{eq:thm-qc2} when $\alpha_t \le 0.75$. Using a similar proof of Theorem \ref{thm:socp2}, we can show that \eqref{eq:thm-outer-approx} is equivalent to constraints \eqref{eq:thm-outer-socp1} -- \eqref{eq:thm-outer-socp2}.
	\Halmos
	\endproof
	
	\section{Proof of Proposition \ref{prop:adjust}}\label{sec:proof-prop-adjust}	
	\proof{Proof of Proposition \ref{prop:adjust}:}
	To show that $Z = Z_1$, we need to show that $Z\subseteq Z_1$ and $Z_1 \subseteq Z$.
	\begin{enumerate}
		\item[(i)]  $Z\subseteq Z_1$.
		
		Given $u_t\in Z$, there exists $\gamma \ge  0$ such that $(u_t, \gamma)$ satisfies \eqref{eq:z-1} and \eqref{eq:z-2}. If $\gamma > 0$, let $\lambda = 1/\gamma$. It is easy to see that $(u_t, \lambda)$ satisfies \eqref{eq:z1-1} and \eqref{eq:z1-2}. For the case $\gamma = 0$, \eqref{eq:z-1} is equivalent to 
		\begin{equation*}\label{eq:z1-gamma-0}
			\left\{ u_t: \ \delta_t \le  \frac{1}{N}\sum_{n = 1}^N \min 
			\left\{0,\ \max \left[\sum_{\ell=1}^{N_{\text{HVAC}}}P_\ell u_{t,\ell} - P_{total,t}^n , 0 \right]\right\}\right\} =  \left\{ u_t: \ \delta_t \le  0\right\}.
		\end{equation*}
		%  \begin{subequations}
		%  \begin{eqnarray}
		%     && \left\{ u_t: \ \delta_t \le  \frac{1}{N}\sum_{n = 1}^N \min 
		%     \left\{0,\ \max \left[\sum_{j=1}^{N_{\text{HVAC}}}P_j u_{t,j} - P_{total,t}^n , 0 \right]\right\}\right\}\\
		%     \label{eq:z1-gamma-0}&=& \left\{ u_t: \ \delta_t \le  0\right\}.
		%  \end{eqnarray}
		%  \end{subequations}
		Since $\delta_t > 0$, the left-hand side of \eqref{eq:z1-gamma-0} is equivalent to an empty set.
		\item[(ii)]  $Z_1 \subseteq Z$.
		
		Given $u_t \in Z_1$, there exists $\lambda \ge 0$ such that $(u_t, \lambda)$ satisfies \eqref{eq:z1-1} and \eqref{eq:z1-2}. Similarly, if $\lambda > 0$, we let $\gamma = 1/\lambda $, which satisfies \eqref{eq:z-1} and \eqref{eq:z-2}. In the case $\lambda = 0$, \eqref{eq:z1} is equivalent to 
		\begin{equation*}
			\left\{ u_t: \ \alpha_t \ge 1\right\}
		\end{equation*}
		which is empty.
		\Halmos
		\endproof	
	\end{enumerate}
	
	\section{Proof of Theorem \ref{thm:adjust-wass}}
	\label{sec:proof-adjust-w}
	
	\proof{Proof:}
	According to Theorem \ref{thm:Wasserstein-LP} and Proposition \ref{prop:strengthen}, for a given pair of ${u}_t$ and ${\alpha}_t$ which is feasible for the adjustable DR chance constraint \eqref{eq:adjustable-drcc-cc}, there exists a $({j}, {k})$ pair such that 
	\begin{eqnarray}
		\label{eq:kj-cond1}&& {{k}}/{N} \le {\alpha}_t < ({k}+1)/N \\
		\label{eq:kj-cond2}&& P_{\text{total},t}^{({j})} < \sum_{\ell = 1}^{N_{\text{HVAC}}} P_\ell {u}_{t,\ell} \le P_{\text{total},t}^{({j}-1)}\\
		\label{eq:kj-cond3}&&   -\frac{1}{N} \sum_{n = {j}}^{{k}} \left(P_{\text{total},t}^{({k}+1)} - P_{\text{total},t}^{(n)}\right) - \left({\alpha}_t - \frac{{j}-1}{N}\right) \left(\sum_{\ell=1}^{N_{\text{HVAC}}} P_\ell {u}_{t,\ell} - P_{\text{total},t}^{({k}+1)}\right) \le -\delta.
	\end{eqnarray}
	% \begin{equation}  {{k}}/{N} \le {\alpha}_t < ({k}+1)/N \text{, and } P_{\text{total},t}^{({j})} < \sum_{\ell = 1}^{N_{\text{HVAC}}} P_\ell {u}_{t,\ell} \le P_{\text{total},t}^{({j}-1)},\end{equation}
	% and 
	% \begin{equation}
	%     -\frac{1}{N} \sum_{n = {j}}^{{k}} \left(P_{\text{total},t}^{({k}+1)} - P_{\text{total},t}^{(n)}\right) - \left({\alpha}_t - \frac{{j}-1}{N}\right) \left(\sum_{\ell=1}^{N_{\text{HVAC}}} P_\ell {u}_{t,\ell} - P_{\text{total},t}^{({k}+1)}\right) \le -\delta.
	% \end{equation}
	% For a given pair of $\hat{u}_t$ and $\hat{\alpha}_t$ 
	%  that are feasible to the DR chance constraint \eqref{eq:drcc-cc} with the Wasserstein ambiguity set $\mathcal D_t^2$, 
	% there exists a $(\hat{j}, \hat{k})$ pair such that (i) \eqref{eq:alpha-con} and \eqref{eq:j-con}  hold; 
	% % and the optimal value of the dual problem \eqref{eq:z-opt-dual-obj}--\eqref{eq:z-opt-dual-3} is \eqref{eq:z-opt-dual-opt-val}.  
	% and (ii) the optimal value of the dual problem \eqref{eq:z-opt-dual-obj}--\eqref{eq:z-opt-dual-3}, or equivalently \eqref{eq:z-opt-dual-opt-val}, is no more than $-\delta$.
	%
	We denote $\Delta_{jk} \in \{0,1\}$ for all $0\le j-1 \le k \le N-1$ such that $\Delta_{jk} = 1$ if we select $j$ and $k$ as the critical index pair;
	% if  $(j,k)$ corresponds to a solution of $(u_t, \alpha_t)$; 
	$\Delta_{jk} = 0$, otherwise. 
	To impose constraint \eqref{eq:kj-cond3}, we require
	% 
	% For all the $(j,k)$ pairs which correspond to a feasible solution of ${u}_t$ and ${\alpha}_t$, we require  
	% $-\frac{1}{N} \sum_{n = {j}}^{{k}} \left(P_{\text{total},t}^{({k}+1)} - P_{\text{total},t}^{(n)}\right) - \left({{\alpha}}_t - \frac{{j}-1}{N}\right) \left(\sum_{\ell=1}^{N_{\text{HVAC}}} P_\ell {u}_{t,\ell} - P_{\text{total},t}^{({k}+1)}\right) \le -\delta$, which is equivalent to letting $   \sum_{j = 1}^{N} \sum_{k=j-1}^{N-1} \Delta_{jk} = 1$ 
	% 
	\begin{equation}\label{eq:z-milp-trilinear}
		\sum_{j = 1}^{N} \sum_{k=j-1}^{N-1}\left[-\frac{1}{N} \sum_{n = {j}}^{{k}} \left(P_{\text{total},t}^{({k}+1)} - P_{\text{total},t}^{(n)}\right) - \left({\alpha}_t - \frac{{j}-1}{N}\right) \left(\sum_{\ell=1}^{N_{\text{HVAC}}} P_\ell {u}_{t,\ell} - P_{\text{total},t}^{({k}+1)}\right)\right] \Delta_{jk} \le -\delta.
	\end{equation}
	Constraint \eqref{eq:z-milp-trilinear} is nonlinear due to  two bilinear terms, i.e., $\alpha_t \Delta_{jk}$ and $u_{t,\ell}\Delta_{jk}$, and one trilinear term,  $\alpha_t u_{t,\ell} \Delta_{jk}$. To linearize them, we introduce $\varepsilon_{jk} = \alpha_t \Delta_{jk}$, $\tau_{\ell jk} = u_{t,\ell}\Delta_{ji}$, $o_{\ell jk} = \alpha_t u_{t,\ell} \Delta_{jk}$ for  $0\le j-1 \le k \le N-1, 1\le \ell \le N_{\text{HVAC}}$, and the  McCormick inequalities \eqref{eq:z-milp-mccormick1}--\eqref{eq:z-milp-mccormick3}.
	
	To ensure the feasibility of the solution $(u_t,\alpha_t)$ associated with a $(j,k)$ pair (there can be multiple solutions associated with one $(j,k)$ pair), we need to further satisfy
	\eqref{eq:kj-cond1} and \eqref{eq:kj-cond2}, which is equivalent to  \eqref{eq:z-milp-alpha-con1} and \eqref{eq:z-milp-j-con1}.
	% we require \eqref{eq:alpha-con} and \eqref{eq:j-con} to be satisfied, which is equivalent to \eqref{eq:z-milp-alpha-con1}--\eqref{eq:z-milp-j-con2}. 
	Therefore, we conclude the proof.
	%  $\sum_{j=1}^{N_{\text{HVAC}}}P_j u_{t,j} - P_{\text{total},t}^{(j-1)} \le 0 < \sum_{j=1}^{N_{\text{HVAC}}}P_j u_{t,j} - P_{\text{total},t}^{(j)}$ 
	\Halmos
	\endproof%\qed
	
	\section{DRCC Model with the Decision of Fleet Size}
	\label{sec:new-model}
	
	In this section, we present a DRCC model that incorporates the decision of the fleet size of residential HVAC units, $N_{\text{HVAC}}$, which is a given parameter in previous models. We associate $N_{\text{HVAC}}$ with a unit penalty cost $c_{N_{\text{HVAC}}}$ in the objective coefficients.
	%Now, we set $N_{\text{HVAC}}$ to be a variable rather than a given parameters. 
	We denote $N_{\text{HVAC}}^\text{U}$  the maximum number of the HVAC units we can deploy to consume the PV generation. For HVAC unit $j$, $j=1,\ldots,N_{\text{HVAC}}^\text{U}$, we introduce a logical binary variable $\zeta_j$ such that $\zeta_j = 1$, if unit $j$ belongs to the fleet, and 0 otherwise. The DRCC model  is formulated as follows.
	%
	%
	%\small
	\begin{subequations}\label{eq:drcc-fleet-size}
		\begin{eqnarray}
			\min_{N_{\text{HVAC}},u_t,\beta_{t,\ell}, x_{t,\ell}t=1,\ldots,N_p} && \sum_{t=1}^{N_p}\left[ c_\text{sys} \sum_{\ell=1}^{N_{\text{HVAC}}^\text{U}}\beta_{t,\ell}  + c_\text{switch} \sum_{\ell=1}^{N_{\text{HVAC}}^\text{U}} u_{t,\ell}\right] + c_{N_{\text{HVAC}}} N_{\text{HVAC}}\\
			\mbox{s.t.} && %\eqref{eq:x-function}, \
			\eqref{eq:absolute}-\eqref{eq:determ-binary} \text{ for $t=1,\ldots,N_p$} \nonumber\\
			&& \inf_{f\in\mathcal D_t} \mathbb P \left( \sum_{\ell=1}^{N_{\text{HVAC}}^\text{U}} P_\ell u_{t,\ell} - \sum_{i=1}^{N_\text{PV}} \widetilde{P}_{\text{PV},t,i} \ge 0 \right) \ge 1-\alpha_t \text{ for $t=1,\ldots,N_p$}\\
			\label{eq:drcc-fleet-size-impose-zero} && \zeta_\ell \le \sum_{t=1}^{N_p} u_{t,\ell} \le N_p\zeta_\ell, \ \ell=1,\ldots, N_{\text{HVAC}}^\text{U}\\
			\label{eq:drcc-fleet-size-capacity}	&& \sum_{\ell=1}^{N_{\text{HVAC}}^\text{U}}\zeta_{\ell} \le N_{\text{HVAC}}\\
			\label{eq:drcc-fleet-size-temp1}	&& x_{1,\ell} = A_\ell x_{0,\ell}\zeta_\ell + B_\ell u_{1,\ell} + G_\ell v_\ell \zeta_\ell + (1-\zeta_\ell)x_\text{ref}, \ \ell = 1,\ldots,N_{\text{HVAC}}^\text{U}\\
			\label{eq:drcc-fleet-size-temp2}	&& x_{t,\ell} = A_\ell(x_{t-1,\ell}-x_\text{ref} + x_\text{ref}\zeta_\ell)  + B_\ell u_{t,\ell} + G_\ell v_\ell \zeta_\ell + (1-\zeta_j) x_\text{ref}, \nonumber\\
			&&  t = 2,\ldots,N_p, \ \ell = 1,\ldots,N_{\text{HVAC}}^\text{U}\\
			\label{eq:drcc-fleet-size-max}	&& 0\le N_{\text{HVAC}}  \le N_{\text{HVAC}}^\text{U}\\
			&& \zeta \in \{0,1\}^{N_{\text{HVAC}}^\text{U}}.
		\end{eqnarray}
	\end{subequations}
	\normalsize
	Constraint \eqref{eq:drcc-fleet-size-impose-zero} requires all $u_{t,\ell}$s' being zeros if $\zeta_\ell = 0$ and thus the HVAC unit $\ell$ is not in the fleet. Constraints \eqref{eq:drcc-fleet-size-temp1} and \eqref{eq:drcc-fleet-size-temp2} ensure that for HVAC unit $\ell$ not in the fleet, i.e., $\zeta_\ell = 0$, the indoor temperatures $x_{t,\ell}$ over all $N_p$ periods are imposed to be $x_\text{ref}$ and thus contribute zero to the objective value. We remark that \eqref{eq:drcc-fleet-size} is a multi-period model over all $N_p$ periods with individual DR chance constraints to guarantee the utilization of the PV generation for each period. All the solution methods and modeling techniques in Sections \ref{sec:DRCC} and \ref{sec:adjustable-cc} can still be applied to \eqref{eq:drcc-fleet-size}.
	\textcolor{black}{\section{CPU Time and Optimality Gaps for Wasserstein Set $\mathcal D_t^2$: Low Risk Requirement $1-\alpha_t$}
		\label{sec:cpu-w-89}}
	\textcolor{black}{
		See Table \ref{tab:cpu-w-low-details}.
	}
	% Table generated by Excel2LaTeX from sheet 'cpu_wasserstein-89'
	\begin{table}[htbp]
		\centering
		\caption{Comparison of CPU time (in seconds) and optimality gaps of low risk requirement $1-\alpha_t$}
		\resizebox{\textwidth}{!}{%
			\textcolor{black}{
				\begin{tabular}{ccc|llc|lll|ll}
					\hline
					\multirow{2}[2]{*}{$1-\alpha$} & \multirow{2}[2]{*}{$\delta$} &       & \multicolumn{3}{c|}{N = 100} & \multicolumn{3}{c|}{N = 500} & \multicolumn{2}{c}{N = 3000} \\
					&       & Instance & \multicolumn{1}{c}{MILP1} & MILP-H & MILP2 & \multicolumn{1}{c}{MILP1} & MILP-H & MILP2 & \multicolumn{1}{c}{MILP-H} & MILP2 \\
					\hline
					\multirow{11}[3]{*}{10\%} & \multirow{11}[3]{*}{0.02} & 1     & 122.45 (1, 0.21\%) & 6.22  & \textbf{0.84} & 2517.31 (9, 13.48\%) & 165.35 & \textbf{3.88} & 2065.91 (14, 5.47\%) & \textbf{38.27} \\
					&       & 2     & 122.61 (1, 0.34\%) & 12.49 & \textbf{0.81} & 2177.88 (6, 6.02\%) & 226.12 & \textbf{4.48} & 2116.55 (16, 3.36\%) & \textbf{47.35} \\
					&       & 3     & 125.81 (1, 0.29\%) & 14.49 & \textbf{0.87} & 2328.41 (5, 5.90\%) & 165.42 & \textbf{4.69} & 2382.27 (20, 4.25\%) & \textbf{45.83} \\
					&       & 4     & 35.72 & 10.88 & \textbf{0.74} & 2186.71 (3, 0.43\%) & 179.87 & \textbf{3.95} & 2447.35 (20, 6.34\%) & \textbf{46.88} \\
					&       & 5     & 65.05 & 4.02  & \textbf{0.86} & 2369.65 (6, 23.48\%) & 202.62 & \textbf{4.17} & 1973.34 (16, 4.01\%) & \textbf{42.08} \\
					&       & 6     & 24.46 & 8.07  & \textbf{0.80} & 2780.93 (15, 19.66\%) & 205.48 & \textbf{4.51} & 2347.32 (19, 4.27\%) & \textbf{142.53 (1, 0.38\%)} \\
					&       & 7     & 123.33 (1, 0.46\%) & 11.19 & \textbf{0.89} & 2528.08 (9, 8.65\%) & 247.28 (1, 0.44\%) & \textbf{4.30} & 2417.75 (21, 4.41\%) & \textbf{46.65} \\
					&       & 8     & 27.46 & 5.32  & \textbf{0.89} & 2339.40 (6, 9.09\%) & 211.23 (1, 0.43\%) & \textbf{4.73} & 2126.13 (17, 4.83\%) & \textbf{42.37} \\
					&       & 9     & 123.79 (1, 0.40\%) & 8.55  & \textbf{0.82} & 2201.85 (5, 12.77\%) & 202.42 & \textbf{5.02} & 2128.63 (17, 3.04\%) & \textbf{39.75} \\
					&       & 10    & 32.43 & 4.76  & \textbf{0.87} & 2361.87 (5, 4.60\%) & 152.37 & \textbf{4.44} & 2049.74 (16, 5.46\%) & \textbf{44.98} \\
					\cline{3-11}          &       & Avg. & 80.31 (0.5, 0.17\%) & 8.60  & 0.84  & 2379.21 (6.9, 10.41\%) & 195.82 (0.2, 0.09\%) & 4.42  & 2205.50 (17.6, 4.54\%) & 53.67 (0.1, 0.04\%) \\
					\hline
					\multirow{11}[3]{*}{10\%} & \multirow{11}[3]{*}{0.2} & 1     & 373.90 (2, 0.47\%) & 138.07 (1, 0.32\%) & \textbf{1.49} & 3113.95 (11, 14.50\%) & 489.31 (2, 0.39\%) & \textbf{9.77} & 3148.18 (25, 6.60\%) & \textbf{219.03 (1, 0.32\%)} \\
					&       & 2     & 312.83 (2, 0.27\%) & 120.61 (1, 0.30\%) & \textbf{4.02} & 3182.16 (12, 27.09\%) & 314.02 (2, 0.39\%) & \textbf{10.51} & 2694.06 (19, 5.38\%) & \textbf{295.30 (2, 0.26\%)} \\
					&       & 3     & 346.11 (2, 0.50\%) & 209.01 (1, 0.32\%) & \textbf{5.26} & 3146.71 (12, 19.78\%) & 328.88 (2, 0.39\%) & \textbf{10.07} & 2854.40 (23, 5.97\%) & \textbf{232.00 (1, 0.30\%)} \\
					&       & 4     & 316.45 (2, 0.39\%) & 221.34 (2, 0.38\%) & \textbf{8.86} & 3171.33 (8, 22.25\%) & 315.03 (2, 0.39\%) & \textbf{9.52} & 2388.50 (17, 4.98\%) & \textbf{221.34 (1, 0.29\%)} \\
					&       & 5     & 317.72 (2, 0.41\%) & 180.47 (1, 0.31\%) & \textbf{7.84} & 3119.98 (12, 26.54\%) & 412.70 (2, 0.39\%) & \textbf{9.97} & 2847.15 (20, 4.53\%) & \textbf{215.50 (1, 0.28\%)} \\
					&       & 6     & 246.86 (1, 0.43\%) & 110.10 (1, 0.42\%) & \textbf{1.55} & 3447.02 (15, 20.60\%) & 450.19 (2, 0.39\%) & \textbf{10.15} & 2885.27 (26, 5.42\%) & \textbf{298.43 (2, 0.22\%)} \\
					&       & 7     & 319.89 (2, 0.36\%) & 215.01 (2, 0.36\%) & \textbf{3.06} & 3215.30 (10, 16.71\%) & 442.92 (2, 0.38\%) & \textbf{9.92} & 2594.05 (20, 6.33\%) & \textbf{246.34 (1, 0.04\%)} \\
					&       & 8     & 224.64 (1, 0.33\%) & 256.21 (2, 0.38\%) & \textbf{1.59} & 3172.53 (14, 12.03\%) & 367.06 (2, 0.39\%) & \textbf{9.59} & 2781.21 (18, 6.75\%) & \textbf{303.22 (2, 0.38\%)} \\
					&       & 9     & 300.75 (2, 0.29\%) & 209.17 (2, 0.34\%) & \textbf{4.61} & 3134.99 (13, 28.85\%) & 396.29 (2, 0.39\%) & \textbf{10.15} & 3099.77 (26, 5.98\%) & \textbf{214.39 (1, 0.28\%)} \\
					&       & 10    & 333.62 (2, 0.39\%) & 213.95 (2, 0.37\%) & \textbf{1.38} & 3293.18 (12, 14.70\%) & 361.87 (2, 0.39\%) & \textbf{10.44} & 2686.10 (21, 5.77\%) & \textbf{229.30 (1, 0.32\%)} \\
					\cline{3-11}          &       & Avg. & 309.28 (1.8, 0.38\%) & 187.39 (1.5, 0.35\%) & 3.97  & 3199.71 (11.9, 20.31\%) & 387.83 (2, 0.39\%) & 10.01 & 2797.87 (21.5, 5.77\%) & 247.49 (1.3, 0.27\%) \\
					\hline
					\multirow{11}[3]{*}{20\%} & \multirow{11}[3]{*}{0.02} & 1     & 34.97 & 5.71  & \textbf{0.87} & 3278.54 (19, 27.73\%) & 184.14 & \textbf{4.70} & 2125.15 (14, 4.82\%) & \textbf{43.55} \\
					&       & 2     & 39.86 & 9.98  & \textbf{0.83} & 3010.49 (17, 24.43\%) & 143.21 & \textbf{4.38} & 2203.64 (16, 3.62\%) & \textbf{35.64} \\
					&       & 3     & 135.69 (1, 0.29\%) & 8.97  & \textbf{0.99} & 3100.17 (15, 34.20\%) & 155.43 & \textbf{4.72} & 1984.68 (15, 3.85\%) & \textbf{35.80} \\
					&       & 4     & 137.47 (1, 0.38\%) & 6.40  & \textbf{0.85} & 3061.05 (15, 28.30\%) & 131.96 & \textbf{4.23} & 2248.16 (17, 5.22\%) & \textbf{36.90} \\
					&       & 5     & 138.02 (1, 0.17\%) & 17.55 & \textbf{0.77} & 2870.46 (16, 31.43\%) & 208.30 & \textbf{4.23} & 2131.41 (17, 3.90\%) & \textbf{35.71} \\
					&       & 6     & 79.53 & 10.39 & \textbf{0.73} & 3068.02 (14, 24.83\%) & 177.49 & \textbf{5.52} & 1643.78 (11, 2.93\%) & \textbf{34.87} \\
					&       & 7     & 64.93 & 8.82  & \textbf{0.80} & 2823.78 (9, 15.16\%) & 317.74 (1, 0.32\%) & \textbf{104.82 (1, 0.10\%)} & 1981.36 (15, 5.53\%) & \textbf{29.77} \\
					&       & 8     & 231.45 (2, 0.19\%) & 105.15 (1, 0.31\%) & \textbf{0.80} & 2850.97 (14, 35.39\%) & 236.45 & \textbf{5.70} & 1799.90 (10, 4.67\%) & \textbf{32.08} \\
					&       & 9     & 134.99 (1, 0.36\%) & 108.33 (1, 0.31\%) & \textbf{1.15} & 3017.16 (11, 25.05\%) & 217.58 (1, 0.31\%) & \textbf{4.48} & 1996.12 (16, 5.39\%) & \textbf{33.86} \\
					&       & 10    & 136.47 (1, 0.22\%) & 12.35 & \textbf{0.79} & 2715.52 (9, 13.52\%) & 89.89 & \textbf{4.19} & 2185.03 (18, 5.88\%) & \textbf{37.21} \\
					\cline{3-11}          &       & Avg. & 113.34 (0.7, 0.16\%) & 29.36 (0.2, 0.06\%) & 0.86  & 2979.62 (13.9, 26.00\%) & 186.22 (0.2, 0.06\%) & 14.70 (0.1, 0.01\%) & 2029.92 (14.9, 4.58\%) & 35.54 \\
					\hline
					\multirow{11}[3]{*}{20\%} & \multirow{11}[3]{*}{0.2} & 1     & 388.77 (1, 0.31\%) & 43.8  & \textbf{1.47} & 3505.28 (11, 24.82\%) & 200.31 (1, 0.36\%) & \textbf{10.38} & 3277.78 (27, 3.97\%) & \textbf{196.14 (1, 0.29\%)} \\
					&       & 2     & 296.35 (1, 0.39\%) & 28.7  & \textbf{1.37} & 3547.98 (15, 18.34\%) & 268.09 & \textbf{40.14} & 3356.53 (28, 3.76\%) & \textbf{188.56 (1, 0.36\%)} \\
					&       & 3     & 309.24 (1, 0.34\%) & 37.9  & \textbf{1.40} & 3495.98 (16, 22.81\%) & 244.00 & \textbf{11.34} & 3468.97 (26, 3.80\%) & \textbf{290.79 (2, 0.28\%)} \\
					&       & 4     & 399.73 (2, 0.56\%) & 81.0  & \textbf{2.23} & 3354.51 (13, 27.32\%) & 226.25 (1, 0.36\%) & \textbf{15.48} & 3365.82 (27, 4.06\%) & \textbf{193.38 (1, 0.34\%)} \\
					&       & 5     & 286.19 (1, 0.37\%) & 11.4  & \textbf{1.51} & 3376.60 (13, 29.87\%) & 368.06 (1, 0.37\%) & \textbf{10.82} & 3527.95 (27, 4.97\%) & \textbf{183.83 (1, 0.27\%)} \\
					&       & 6     & 362.52 (1, 0.44\%) & 13.8  & \textbf{1.26} & 3384.83 (9, 33.56\%) & 298.02 (1, 0.36\%) & \textbf{110.27 (0.045\%)} & 3354.24 (26, 4.17\%) & \textbf{189.52 (1, 0.27\%)} \\
					&       & 7     & 400.97 (1, 0.38\%) & 9.1   & \textbf{2.12} & 3474.31 (15, 23.25\%) & 337.58 & \textbf{83.60} & 3191.83 (27, 4.10\%) & \textbf{222.98 (1, 0.36\%)} \\
					&       & 8     & 424.35 (2, 0.42\%) & 9.9   & \textbf{1.56} & 3412.10 (11, 26.36\%) & 206.94 (1, 0.36\%) & \textbf{110.55 (1, 0.23\%)} & 3447.36 (29, 4.11\%) & \textbf{187.14 (1, 0.23\%)} \\
					&       & 9     & 255.55 & 11.9  & \textbf{1.47} & 3412.99 (12, 26.23\%) & 283.64 & \textbf{28.68} & 3464.42 (29, 3.64\%) & \textbf{92.86} \\
					&       & 10    & 326.42 (1, 0.34\%) & 11.0  & \textbf{1.38} & 3465.78 (14, 19.64\%) & 249.26 (1, 0.35\%) & \textbf{109.15 (1, 0.23\%} & 3425.87 (27, 4.04\%) & \textbf{382.44 (3, 1.04\%)} \\
					\cline{3-11}          &       & Avg. & 345.01 (1.1, 0.36\%) & 25.9  & 1.58  & 3443.04 (12.9, 25.22\%) & 268.21 (0.6, 0.22\%) & 53.04 (0.3, 0.06\%) & 3388.08 (27.3, 4.06\%) & 212.76 (1.2, 0.34\%) \\
					\hline
				\end{tabular}%
				\label{tab:cpu-w-low-details}%
		}}
	\end{table}%

	\section{CPU Time and Optimality Gaps for Adjustable DRCC under Wasserstein Set $\mathcal D_t^2$}
	\label{sec:wass-adjustable-cpu}
	Following the in-sample data generation procedure in Section \ref{sec:comp-setup},
	we generate 10 instances (each of 10 samples) under the sunny weather condition and solve them using DRCC-W models by MILP3 and MILP4 reformulations, respectively.
	The CPU time limit for each period is  100 seconds. In Table \ref{tab:adjustable-wass}, MILP4 solves all instances faster with an average of 1154.29 seconds than 3449.39 seconds of MILP3. Among the 53 periods solved, MILP3 has more periods not  solved optimally. The average gap of the unsolved periods is up to 12.38\%. While MILP4 only yields a gap of 0.73\% as the MILP4 formulation provides a tighter linear relaxation than the MILP3. 
	% Table generated by Excel2LaTeX from sheet 'Adjustable'
	\begin{table}[htbp]
		\centering
		\caption{Comparison of CPU time (in seconds) and optimality gaps}
		\resizebox{0.5\textwidth}{!}{%
			\begin{tabular}{ccccccc}
				\hline
				\multirow{2}[2]{*}{Instance} & \multicolumn{3}{c}{MILP3} & \multicolumn{3}{c}{MILP4} \\
				& CPU  & \# Limit & Gap   & CPU & \# Limit & Gap \\
				\hline
				1     & 3410.35 & 33    & 12.35\% & 1165.00 & 8     & 0.70\% \\
				2     & 3456.28 & 34    & 12.44\% & 1069.39 & 7     & 0.41\% \\
				3     & 3460.71 & 33    & 12.84\% & 1557.45 & 13    & 0.68\% \\
				4     & 3453.83 & 34    & 11.68\% & 1099.37 & 8     & 0.88\% \\
				5     & 3486.60 & 34    & 12.13\% & 1307.64 & 11    & 0.77\% \\
				6     & 3459.33 & 34    & 12.47\% & 1048.38 & 8     & 1.03\% \\
				7     & 3461.45 & 34    & 12.43\% & 1380.35 & 11    & 0.76\% \\
				8     & 3457.13 & 33    & 11.88\% & 1194.67 & 9     & 0.71\% \\
				9     & 3392.60 & 33    & 12.57\% & 904.65 & 4     & 0.43\% \\
				10    & 3455.67 & 33    & 13.00\% & 816.00 & 5     & 0.91\% \\
				\hline
				Avg. & 3449.39 & 34    & 12.38\% & 1154.29 & 8     & 0.73\% \\
				\hline
			\end{tabular}%
		}
		\label{tab:adjustable-wass}%
	\end{table}%

	\textcolor{black}{\section{CPU Time and Optimality Gaps for Moment-based Ambiguity Set}
		\label{sec:cpu-moment}}
	
	\textcolor{black}{In Table \ref{tab:cpu-moment}, the total CPU time  of solving all 53 periods are reported for the DRCC and adjustable variants under the moment-based set $\mathcal D_t^1$. We set $\gamma_1 = 0$, $\gamma_2 = 1$, and $1-\alpha_t = 80\%$.  Column ``\# Limit'' indicates the number of
		periods that cannot be solved  when  the time limit is reached. The average optimality gaps of these (unsolved) periods are presented in the next column ``Gap.'' For the adjustable DRCC, as $\gamma_1/\gamma_2 = 0$, only SOCP1 or SOCP3 is required. Note that in our problem settings, $\alpha_t$ does not exceed 75\%. Therefore SOCP1 and SOCP3 are equivalent. For all instances, SOCP3 yields shorter CPU times and smaller optimality gaps than SOCP1.
	}
	\\
	% Table generated by Excel2LaTeX from sheet 'cpu-adjustable-moment'
	%
	\begin{table}[htbp]
		\centering
		\textcolor{black}{
			\caption{Comparison of CPU time (in seconds) and optimality gaps for DRCC and adjustable variants under $\mathcal D_t^1$}
			\resizebox{.65\textwidth}{!}{%
				\begin{tabular}{c|rcl|rcl|rcl|}
					\hline
					\multirow{3}[2]{*}{Instance} & \multicolumn{3}{c|}{\multirow{2}[1]{*}{DRCC}} & \multicolumn{6}{c}{adjustable DRCC} \\
					& \multicolumn{3}{c|}{} & \multicolumn{3}{c|}{SOCP1} & \multicolumn{3}{c}{SOCP3} \\
					& \multicolumn{1}{c}{CPU} & \multicolumn{1}{c}{\# Limit} & \multicolumn{1}{c|}{Gap} & \multicolumn{1}{c}{CPU} & \multicolumn{1}{c}{\# Limit} & \multicolumn{1}{c|}{Gap} & \multicolumn{1}{c}{CPU} & \multicolumn{1}{c}{\# Limit} & \multicolumn{1}{c}{Gap} \\
					\hline
					1     & 0.31  & 0     & N/A   & 3566.69 & 31    & 5.19\% & 1243.26 & 12    & \multicolumn{1}{l}{0.06\%} \\
					2     & 0.26  & 0     & N/A   & 3586.60 & 32    & 4.73\% & 1437.35 & 13    & \multicolumn{1}{l}{0.06\%} \\
					3     & 0.19  & 0     & N/A   & 3525.74 & 29    & 5.22\% & 1554.54 & 14    & \multicolumn{1}{l}{0.09\%} \\
					4     & 0.32  & 0     & N/A   & 3624.84 & 30    & 5.02\% & 1028.21 & 10    & \multicolumn{1}{l}{0.09\%} \\
					5     & 0.28  & 0     & N/A   & 3536.61 & 30    & 5.17\% & 1455.35 & 13    & \multicolumn{1}{l}{0.07\%} \\
					6     & 0.27  & 0     & N/A   & 3510.96 & 31    & 4.65\% & 1747.06 & 16    & \multicolumn{1}{l}{0.09\%} \\
					7     & 0.24  & 0     & N/A   & 3589.61 & 31    & 4.79\% & 1352.62 & 13    & \multicolumn{1}{l}{0.09\%} \\
					8     & 0.32  & 0     & N/A   & 3432.83 & 30    & 5.05\% & 925.62 & 9     & \multicolumn{1}{l}{0.08\%} \\
					9     & 0.28  & 0     & N/A   & 3598.59 & 29    & 5.09\% & 1220.28 & 12    & \multicolumn{1}{l}{0.04\%} \\
					10    & 0.31  & 0     & N/A   & 3607.45 & 32    & 5.10\% & 1107.61 & 10    & \multicolumn{1}{l}{0.11\%} \\
					\hline
					Avg.  & 0.28  & 0     & N/A   & 3557.99 & 31    & 5.00\% & 1307.19 & 12    & 0.08\% \\
					\hline
				\end{tabular}%
				\label{tab:cpu-moment}%
		}}
		
	\end{table}%
	
	\section{Tracking Performance and Room Temperatures for Cloudy Weather}
	\label{supplement-sec:tracking-temperature-cloudy}
	
\begin{figure*} 
		\centering
		\begin{subfigure}{0.4\textwidth}
			\includegraphics[width=.9\linewidth]{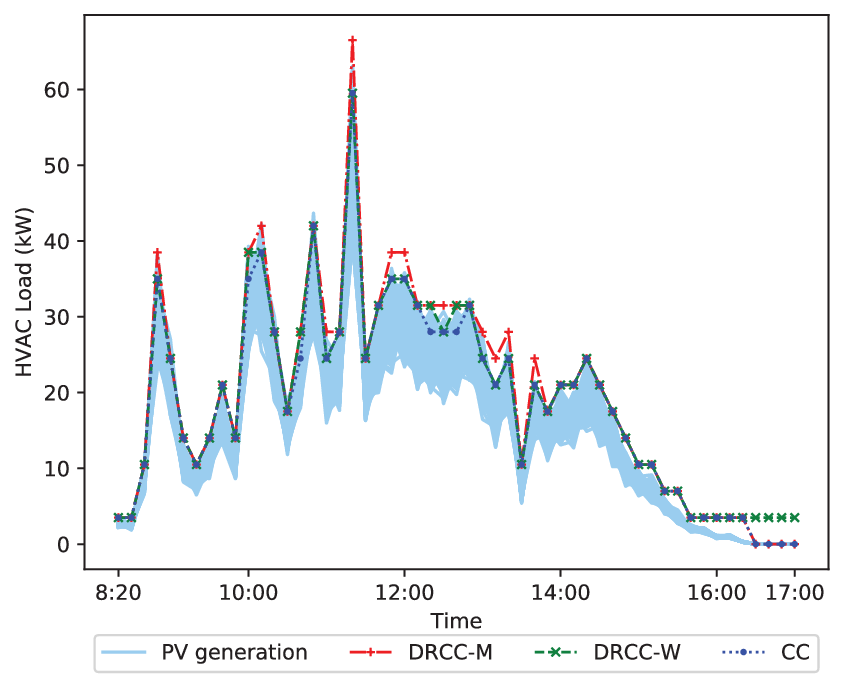} \caption{PV profile tracking} 
			\label{fig:PVtracking-cloudy}
		\end{subfigure}
		\begin{subfigure}{0.4\textwidth}
			\includegraphics[width=\linewidth]{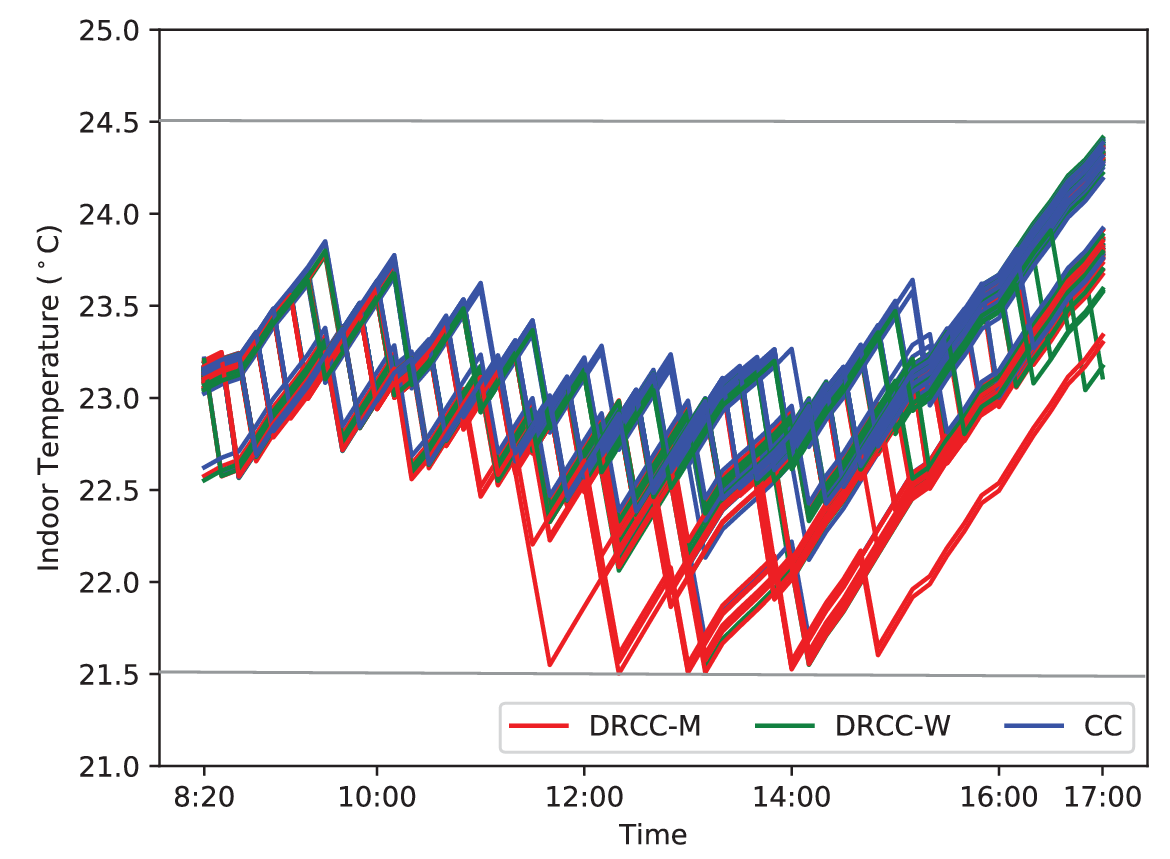} \caption{Room temperature} 
			\label{fig:temp-cloudy}
		\end{subfigure}
		\caption{(Color online)  PV profile tracking and room temperatures of 35 buildings under sunny  weather}
		\label{fig:cloudy-tracking-temp}
	\end{figure*}

	%For the cloudy weather, the number of HVAC units $N_\text{HVAC} = 35$.
	Figure \ref{fig:PVtracking-cloudy} shows the PV profile tracking under the cloudy weather. The three models track the PV profile well most of the time. Around 4:00 pm, when the PV generation is low, the optimal schedule of all three models do not track the PV generation as closely as before. In Figure \ref{fig:temp-cloudy}, again, all three models keep the room temperature within the comfort band and DRCC-M provides relatively lower temperature.

	%For all three models, we solve them  for all $N_p$ periods under both the sunny weather and cloudy weather conditions, respectively. 
	%{We note that the PV generation of the cloudy day is relatively lower than that of the sunny day (see in Figure \ref{fig:PV}).  Consequently, when solving instances of the cloudy weather condition, to keep the total number of enrolled HVAC devices compatible with the magnitude of local solar PV generation \citep{dong2018model}, we consider enrolling a subset of 35 residential HVAC units. 
	%	That is, $N_{\text{HVAC}} = 35$ under the cloudy weather. While, we consider enrolling all the 100 HVAC units for the sunny day. \textcolor{black}{We remark that our model is flexible to incorporate the decision of the fleet size, which can be decided based on the nameplate capacity of HVAC devices and local solar PV generation scales.} The details are in the online appendix \ref{sec:new-model}.   }  
	
	%In Figure \ref{fig:PVtracking-cloudy} under the cloudy weather, around 4:00 pm, when the PV generation is low, the optimal schedule of all three models do not track the PV generation as closely as before.
	
	\section{Out-of-sample Performance: Sunny vs. Cloudy }
	\label{supplement-sec:out-of-sample}
	
\begin{figure*} 
		\centering
		\begin{subfigure}{0.4\textwidth}
			\includegraphics[width=\linewidth]{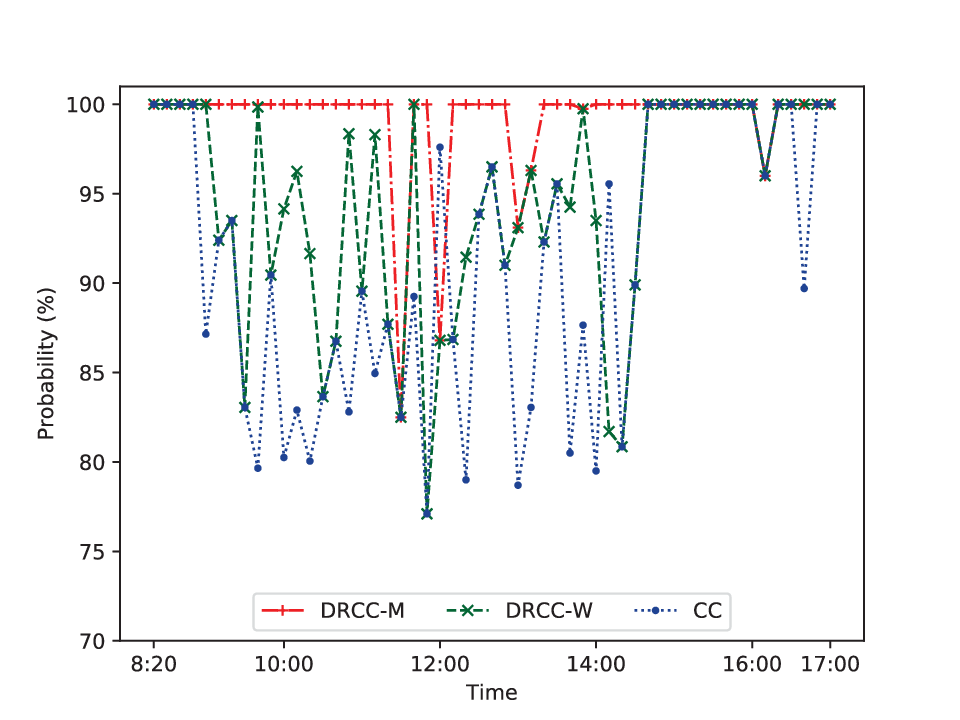} \caption{Sunny} 
		\end{subfigure} 
		\begin{subfigure}{0.4\textwidth}
			\includegraphics[width=\linewidth]{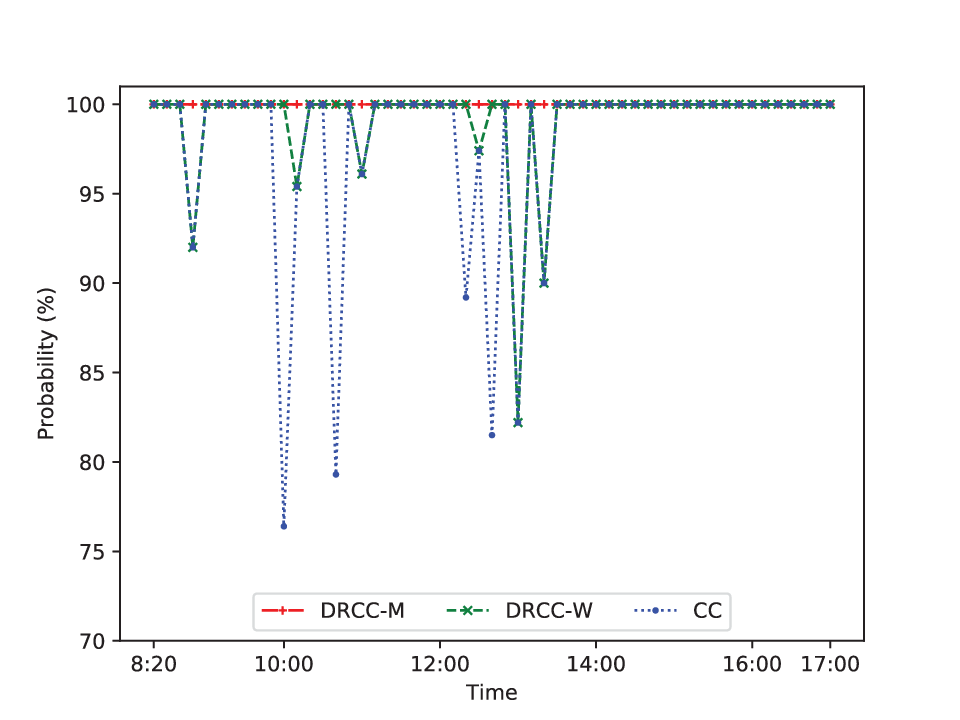} \caption{Cloudy} 
		\end{subfigure}
		\caption{(Color online) Probabilities of locally consuming PV generation under sunny and cloudy weather conditions}
		\label{fig:prob_uniform}
	\end{figure*}

	In Figure \ref{fig:prob_uniform}, the 95th percentile of probabilities  are shown for all three models under the sunny and cloudy weather conditions. In both plots, the two DRCC models perform better than the CC models. Again, as the DRCC-M model is more conservative,  the DRCC-M model  achieves higher probability than the DRCC-W model.  

\end{appendices}

\end{document}